\documentclass{article}
\usepackage{zibtitlepage}

\usepackage[utf8]{inputenc}
\usepackage[a4paper, margin={25mm}]{geometry}

\usepackage{amsmath,amssymb}
\usepackage{mathrsfs} % to use \mathscr on capital letters

\usepackage{mathtools} % for \mathmakebox
\usepackage{calc} % for \widthof

\usepackage{accents} % needed for commands
\usepackage{url}
\usepackage[table]{xcolor}

\usepackage{natbib}
\setcitestyle{authoryear,square,semicolon}
\usepackage{tikz}
\usetikzlibrary{circuits}
\usepackage{pgfgas}
\usetikzlibrary{fadings}
\usetikzlibrary{decorations.pathreplacing} % curly braces

\usepackage{graphicx}
\usepackage[format=hang]{caption} % multi-line captions do not start at left, but after "Figure X:"
\usepackage{subcaption} % for multiple pictures with

\usepackage[linesnumbered, ruled, vlined]{algorithm2e}

\usepackage{todonotes}
\input{commands}

\let\oldnl\nl% Store \nl in \oldnl
\newcommand{\nonl}{\renewcommand{\nl}{\let\nl\oldnl}}% Remove line number for one line
\newcommand{\dotcup}{\ensuremath{\mathbin{\dot{\cup}}}\xspace}
\newcommand{\dotbigcup}{\ensuremath{\mathop{\dot{\bigcup}}}\xspace}

\newcommand{\MINLP}{\ensuremath{\mathscr{P}}\xspace}
\newcommand{\inst}{\ensuremath{\Phi}\xspace} % problem instance
\newcommand{\instAgg}{\ensuremath{\inst^\mathrm{agg}}\xspace} % problem instances with aggregated time horizon
\newcommand{\instFixed}{\ensuremath{\inst^\mathrm{fixed}}\xspace} % problem instances with partially fixed binary variables
\newcommand{\paramMaxIt}{\ensuremath{n^\mathrm{it}}\xspace} % max number of iterations
\newcommand{\paramMaxItAgg}{\ensuremath{n^\mathrm{itAgg}}\xspace} % max number of iterations in the aggregated MIP
\newcommand{\paramMaxItFull}{\ensuremath{n^\mathrm{itFull}}\xspace} % max number of iterations in the whole time horizon MIP in the aggregated heuristic
\newcommand{\paramMaxRelTol}{\ensuremath{n^\mathrm{tolRel}}\xspace} % relative convergence tolerance
\newcommand{\paramMaxAbsTol}{\ensuremath{n^\mathrm{tolAbs}}\xspace} % absolute convergence tolerance
\newcommand{\paramRolHorizon}{\ensuremath{n^\mathrm{rolHor}}\xspace} % rolling horizon size in corresponding heuristic
\newcommand{\paramAggHorizon}{\ensuremath{n^\mathrm{aggHor}}\xspace} % time aggregated horizon size in corresponding heuristic

\newcommand{\paramNodeLimitMin}{\ensuremath{n^\mathrm{limitMin}}\xspace} % dynamic node limit: minimal node limit
\newcommand{\paramNodeLimitRel}{\ensuremath{n^\mathrm{limitRel}}\xspace} % dynamic node limit: new relative node limit (> 1.0)
\newcommand{\paramImprovement}{\ensuremath{n^\mathrm{improve}}\xspace} % minimal improvement to trigger a new node limit

\newcommand{\pressPrevI}[1]{\ensuremath{\press'_{#1}}\xspace}
\newcommand{\mFlowPrevI}[1]{\ensuremath{\mFlow'_{#1}}\xspace}

\begin{document}

\ZTPAuthor{
\ZTPHasOrcid{Felix Hennings}{0000-0001-6742-1983},
\ZTPHasOrcid{Kai Hoppmann-Baum}{0000-0001-9184-8215},
\ZTPHasOrcid{Janina Zittel}{0000-0002-0731-0314}}
\ZTPTitle{Optimizing transient gas network control for challenging real-world instances using MIP-based heuristics}
\ZTPNumber{22-08}
\ZTPMonth{May}
\ZTPYear{2022}
\ZTPInfo{}
%\ZTPInfo{The content of this report is also available as publication in <Journal>. Please always cite as:\\[1ex] <authors>. <title>. <journal>, <year>, \href{<link>}{\textcolor{blue}{DOI:~<DOI>}}}

\title{Optimizing transient gas network control for challenging real-world instances using MIP-based heuristics}
\author{
\ZTPHasOrcid{Felix Hennings}{0000-0001-6742-1983}\and\
\ZTPHasOrcid{Kai Hoppmann-Baum}{0000-0001-9184-8215}\and\
\ZTPHasOrcid{Janina Zittel}{0000-0002-0731-0314}}
\hypersetup{pdftitle={Optimizing transient gas network control for challenging real-world instances using MIP-based heuristics},
  pdfauthor={Felix Hennings, Kai Hoppmann-Baum, Janina Zittel}}

\zibtitlepage
\maketitle

\begin{abstract}
Optimizing the transient control of gas networks is a highly challenging task.
The corresponding model incorporates the combinatorial complexity of determining the settings for the many active elements as well as the non-linear and non-convex nature of the physical and technical principles of gas transport.
In this paper, %/talk
we present the latest improvements of our ongoing work to solve this problem for real-world, large-scale problem instances:
By adjusting our mixed-integer non-linear programming model regarding the gas compression capabilities in the network, we reflect the technical limits of the underlying units more accurately while maintaining a similar overall model size.
In addition, we introduce a new algorithmic approach that is based on splitting the complexity of the problem by first finding assignments for discrete variables and then determining the continuous variables as locally optimal solution of the corresponding non-linear program.
For the first task, we design multiple different heuristics based on concepts for general time-expanded optimization problems that find solutions by solving a sequence of sub-problems defined on reduced time horizons.
To demonstrate the competitiveness of our approach, we test our algorithm on particularly challenging historic demand scenarios.
The results show that high-quality solutions are obtained reliably within short solving times, making the algorithm well-suited to be applied at the core of time-critical industrial applications.
\end{abstract}
%
%
% define algorithm keyword break after ZIB title page (it uses some break as well ...)
\SetKw{break}{break}
\section{Introduction}
The natural gas supply of Europe and especially Germany is subject to a disruptive transformation.
With the possibility of cutting the supply from Russia and the recently solidified plans to build LNG terminals in Germany as published by the government~\citep{Fed2022}, new flow patterns are likely to emerge in the corresponding transport networks.
Moreover, the usage of fossil-based natural gas will soon be terminated in favor of more sustainable alternatives as part of the energy transformation towards massive reductions in CO$_2$ emissions.
In particular, the natural gas networks may be used to store synthesized methane generated via power-to-gas from excess renewable energy or could (partially) be repurposed for hydrogen transport, see, for example, the scenario framework for the German \emph{Gas Network Development Plan} as published by the German gas network operators~\citep{FNB2021}.
Especially when applying power-to-gas strategies, the intraday supply patterns are likely to become more volatile, complicating the control of the networks even further.

Mathematical optimization can assist the network operators in preparing for these new challenges and ensuring a secure and efficient network operation.
However, solving the corresponding optimization problems is hard.
On the one hand, the physics of gas transport in pipelines describing the interplay of gas pressure and flow is captured by non-convex equations.
On the other hand, the network contains many active elements to control the gas flow, like valves to dynamically change the network topology and compressors to increase the gas pressure.
Apart from the non-convex description of the maximum power available to compressors, these control options also introduce an immense combinatorial complexity to the problem.
Since the problem is both practically relevant and theoretically challenging, it has been intensively studied in the past, and we will give an overview split into two categories in the following.

First, there is the \emph{stationary} variant of the problem.
It defines an equilibrium state of the network by assuming a single infinitely long time step and can, for example, be used for planning purposes.
For a profound introduction to stationary gas network optimization, we recommend the book of~\citep{KocHilPfeSch2015} or the paper of~\citep{PfeFugGeiGei2015} for a summary of the corresponding findings.
They model the problem as a mixed-integer non-linear program (MINLP) and deal with the two sources of complexity separately by first finding solutions for the discrete decision variables and afterwards completing these to highly precise solutions by solving a corresponding non-linear programming (NLP) model.
A more recent study on the stationary problem is presented in~\citep{GonGonFer2019}.
They use a model based on simplified structures representing the network areas containing the compressors.
Algorithmically, the authors mainly focus on solving an NLP-model variant of the problem with a two-step sequential linear programming (SLP) approach, which is inspired by an early attempt at solving the problem reported in~\citep{PraWil1984}.
However, by replacing the first step with sequential mixed-integer programming, they can convert their algorithm into a heuristic for MINLP problems.
Test on real-world instances, which only have a limited number of binary variables, reveal a performance similar to state-of-the-art NLP and MINLP solvers.

Second, the \emph{transient} variant of the problem aims at finding a short-term control for the network over discrete future time steps.
An early publication on the topic was the thesis of~\citep{Mor2007}.
Here, the problem was formulated as a pure mixed-integer program (MIP) by approximating the non-linear pipe equations and compression power constraints using piece-wise linear functions.
To aid the corresponding solver, the author proposed a special branching scheme for the piece-wise linear functions and used a simulated annealing heuristic, which was also described in~\citep{MahMarMor2007}.
Shortly after, \citep{DomGeiKolLan2011} presented their approach for a similar problem formulation, which again used piece-wise linear functions to approximate the non-linear parts of the network.
They iteratively improved the linearization by first solving the corresponding MIP model, then using the binary solution values to form a non-linear programming model using the original not approximated equations, and afterwards refining the linearization by adding mesh-points corresponding to the found solutions.
The procedure is stopped when a predefined approximation threshold is satisfied.
In their computational results, they report problems in solving the single MIP models to optimality even after multiple hours.
In a more recent study, \citep{HahLeyZav2017} directly solve an MINLP formulation of the problem.
Although the problem is non-convex, they use a convex solver and apply a customized branching rule for the problem.
In their experiments, they design scenarios for small test networks of less than 10 nodes that include the reversal of the gas flow direction on single pipelines, which was not possible in the networks used in~\citep{Mor2007} and \citep{DomGeiKolLan2011}.
Despite the size of the networks, the different variants of their algorithm need, even in the best cases, hours to finish.
Another contribution was made by~\citep{BurEggGroMar2019}.
By introducing a new discretization for the differential equations describing the gas flow in pipes, they were able to reformulate the transient problem as a stationary one on a time-expanded graph and use a specialized algorithm designed for the stationary case to solve the problem.
It is based on iteratively solving mixed-integer programming relaxations of the problem, which they again obtained by applying a piece-wise linearization.
By refining the linearization in each iteration, their approach converges to the globally optimal solution.
To speed up the solving process, the intermediate MIP solutions were completed to valid primal solutions of the overall problem using an NLP solver.
For the one evaluated instance, the algorithm finds a solution with a gap of less than 10\% in less than 60 seconds and less than 8\% in half an hour.
However, after this point, even single iterations of the algorithm reach run times of over an hour and may not even improve the gap.

Finally, we mention the work of~\citep{HopHenLenKoc2021}, who presented a decision support system for network operators to ensure a secure and efficient network operation.
The system is split into two separate stages:
First, an MINLP problem formulation covering the complete network is solved, which approximates the capabilities of the technical elements located at large pipeline intersections, the \emph{network stations}.
In particular, each of these stations is replaced with a simplified graph model using artificial arcs.
Non-technical control measures available to the network operators are realized as slack variables to adjust the future supplies and demands in terms of pressure and inflow.
Their usage as well as changes in the technical control are penalized in the objective functions of a tri-level optimization model.
To derive a feasible solution, a series of linearized MIP models are solved, which use two types of primal heuristics to determine good initial solutions.
Afterwards, a linear programming (LP) model is solved to smooth the solution, followed by an SLP approach that transforms the solution into one respecting the original non-linear constraints of the MINLP problem.
In a computational study on subsequent real-world instances spanning multiple days, the algorithm's run time and the differences in the solutions are evaluated.
In the second stage of the decision support system, the actual control measures recommended to the network operators are determined by solving highly detailed models for each network station individually, as described in~\citep{HenAndHopTur2021}.
Note that since the overall system runs in regular intervals to provide up-to-date control recommendations, it must reliably satisfy strict time limit restrictions.

In this paper, we present several improvements to the simplified network-wide model of the first stage of the decision support system described in~\citep{HopHenLenKoc2021}:
We suggest an alternative model for the artificial network station arcs representing the gas compressors, which is closer to the compression capabilities of the highly detailed network station model used in the second stage without increasing the model's size.
Furthermore, we formulate a new algorithm to solve the MINLP problem of the first stage, which finds better solutions faster than the original approach.
In particular, we first determine sets of values for all the binary variables of the model, which we call \emph{binary assignments}, and afterwards complete these to full solutions of the MINLP by solving the remaining non-linear programming (NLP) model.
The binary assignments are obtained by combining a sequential mixed-integer programming (SMIP) routine with heuristics concepts for time-expanded problems: The rolling horizon and the aggregate horizon heuristic.
We additionally apply a dynamic branch-and-bound node limit to the single model solves to further reduce the maximum run times.
Finally, we thoroughly test the proposed algorithm's run time and its solution quality against valid lower bounds.
As a test set, we identify especially hard-to-solve real-world instances of the problem that feature many supply and demand changes.
Thereby we verify that our algorithm is suitable for time-critical industry applications and more volatile future supply patterns.

%%% motivation for our approach
Our proposed approach focuses on finding good primal solutions fast and does not provide a valid lower bound on the value of the problem's optimal solution.
This enables us to tackle larger and more challenging problem instances than those found in the literature while respecting strict time limits.
We follow the idea of dealing with the combinatorial and non-linear aspects of the problem separately and also use sequential mixed-integer programming for finding good binary assignments.
However, since we do not aim at finding lower bounds, we do not need to use an iteratively refined piece-wise linear formulation but only an updatable linearization, similar to what~\citep{GonGonFer2019} propose for the first level of their algorithm.
As a consequence, the general size of the MIP problems is smaller and does not increase during the solving process, which is one of the reasons for the reported slow lower bound convergence of the transient approaches listed above.

We note that the heuristic concepts we use for time-expanded problems are well known in the literature.
The rolling horizon approach was, for example, used for the optimization of product manufacturing in~\citep{TiaSae2012}, distributed power generation in~\citep{Cap2017}, or as one of the primal heuristics in~\citep{HopHenLenKoc2021}.
A theoretical analysis of rolling horizon approaches, as well as a more extensive list of applications, is given in \citep{GloLieRos2022}.
Similarly, heuristics using an aggregation of a time horizon were often described, as in the example of scheduling the iron ore production in \citep{NewKuc2007} or the extension of freight transport networks in \citep{BolErnKalRoc2013}.

The rest of the paper is structured as follows:
Section~\ref{sec:model} introduces the MINLP problem formulation, where we particularly highlight the differences with respect to the model used in~\citep{HopHenLenKoc2021}.
In the following Section~\ref{sec:mip_based_heuristics}, we present our heuristic algorithmic approach for solving this problem.
The computational results are the content of Section~\ref{sec:computational_experiments}, where we describe the creation of the instance sets, compare a base version of our algorithm against a global MINLP solver in terms of solution quality on smaller instances, and then compare the actually proposed algorithm against the base algorithm version on a more challenging instance set.
We close with concluding remarks in Section~\ref{sec:outlook}.
\section{Model formulation}\label{sec:model}
Our problem is based on the gas flow model used for the first stage of the decision support system presented in \citep{HopHenLenKoc2021}.
As in there, we assume that the gas network contains a set \setNetworkStations of \emph{network stations}.
Those are sub-networks that include the majority of all actively controllable network elements, among them all the network's compressor stations, and are mainly located at large pipeline intersections.
For the model, the inner network topology of each station is replaced by a set of artificial arcs representing an approximation of the station's gas routing, compression, and regulation capabilities.
This artificial model was created manually by industry experts who work for the network operator and have an excellent knowledge of the individual stations and their properties.

In the following, we summarize the differences between our problem formulation and the one of \citeauthor{HopHenLenKoc2021}.
Note that we refer to their model as the \emph{base model} for the remainder of the paper.
\begin{itemize}
\item As the most crucial change, we introduce in Section~\ref{sec:auxCompressor} a new model for the artificial compressor arcs in the network stations. It is based on the usage of configurations for each artificial arc, which we construct based on the properties of the compressor stations of the original network station topology. We show that the model is much closer to the highly detailed compressor station model used in~\citep{HenAndHopTur2021} while maintaining a similar model size with respect to the number of variables and constraints.
\item The second major difference is the introduction of artificial valves at the boundaries of the network stations in Section~\ref{sec:flow_dirs_and_fence_node_valves}. With these, it is possible to decouple the pressure inside and outside the station at these points. Note that this was initially suggested in \citep{Hop2022}.
\item In this paper, we do not use a tri-level model to reflect strict objective function priorities for using certain gas network control measures.
Instead, we use a single objective function representing the weighted sum of different terms to penalize them, see Section~\ref{sec:objective}.
This objective also includes \emph{smoothing} terms, which minimize pressure and inflow changes at the network station boundaries over time.
As a final change in the objective, we include a newly created incentive to reduce the usage of compressor units based on the new artificial compressor configurations to avoid unnecessary energy consumption.
\item We use the regulator model of \citep{HenAndHopTur2021}, which is slightly different from the base regulator model.
\item In the artificial arc models introduced in Section~\ref{sec:auxArcs}, we do not include \emph{combined arcs}, which have gas compression as well as regulation capabilities, as they do not appear in our real-world-based test instances, see Table~\ref{tab:stationStatistics} in Section~\ref{sec:instances}.
For the same reason, we did not include artificial \emph{bi-directional compressor arcs} or \emph{bi-directional combined arcs}.
\item Finally, we also do not introduce \emph{flow direction conditions} as they are not used in the network station description of our test instances.
\end{itemize}

When formulating our model below, we often use \emph{indicator constraints}, which have two parts: The first one is a logical condition in the form of a binary variable $y$ attaining a value $b\in\{0,1\}$ and the second one an implied linear constraint $a^T x \leq a_0$ of variables $x=\{x_1,\dots,x_n\}$ and coefficients $a\in\mathbb{R}^n$, $a_0\in\mathbb{R}$, which only holds if the logical condition is true.
We denote such an indicator constraint by
\begin{equation*}
    y = b \quad\implies\quad a^T x \leq a_0.
\end{equation*}
Note that we also use
\begin{equation*}
    y = b \quad\implies\quad a^T x = a_0
\end{equation*}
as an abbreviation of the two indicator constraints
\begin{align*}
    y = b \quad&\implies\quad  a^T x \leq a_0 & \text{and}& & y = b \quad&\implies\quad -a^T x \leq -a_0.
\end{align*}
Indicator constraints can be reformulated as linear constraints via the \emph{big-M} approach in case all the variables $x=\{x_1,\dots,x_n\}$ are bounded, see~\citep{BonLodTraWie2015} for a general introduction to the topic.
A list of all the variables used in our model, together with their domains, meanings, and units, is given in Table~\ref{tab:Variables} in the Appendix.

\subsection{Gas flow in networks}\label{sec:gas_flow_in_networks}
We model a gas network as a directed graph \generalGraph with \setVertices being the set of nodes or vertices and $\setArcs = \setPipes \dotcup \setValves \dotcup \setRegulators \dotcup \setArtificialLinks$ being the set of arcs representing the single element sets present in the network: Pipes \setPipes, valves \setValves, regulators \setRegulators, and artificial network station arcs \setArtificialLinks.
The latter is again composed of the different types of artificial arcs found in the network stations: Shortcuts \setShortcuts, regulating arcs \setRegulatorArcs, compressor arcs \setCompressorArcs, and bi-regulating arcs \setBiRegulatorArcs.
For station $i\in\setNetworkStations$, we denote by $\setArtificialLinksNaviStation{i}=\setShortcutsNaviStation{i}\dotcup\setRegulatorArcsNaviStation{i}\dotcup\setBiRegulatorArcsNaviStation{i}\dotcup\setCompressorArcsNaviStation{i}$ we denote its set of artificial arcs and by $\setFenceNodesNaviStation{i}\subseteq\setVertices$ the nodes forming the station's boundary, which we call \emph{fence nodes}.
Formally, the fence nodes are those nodes having incident arcs from within the station as well as from outside of the station, i.e.,
\begin{align*}
    \forall v\in\setFenceNodesNaviStation{i} : (\exists a^\mathrm{in}\in\setArtificialLinksNaviStation{i}: a^\mathrm{in}=(v,r) \lor a^\mathrm{in}=(\ell,v)) \land (\exists a^\mathrm{out}\in\setArcs\setminus\setArtificialLinksNaviStation{i}: a^\mathrm{out}=(v,r) \lor a^\mathrm{out}=(\ell,v)).
\end{align*}
Note that, in theory, a node can be the fence node of different stations.
However, this does not happen in practice.
Furthermore, there can be artificial non-fence nodes inside the station, which are needed for the manual station modeling.
Finally, we assume for each station $i\in\setNetworkStations$ that the graph induced by the station's arcs \setArtificialLinksNaviStation{i} is connected.

Furthermore, we consider a discrete time horizon \setTimestepsAll, where \setTimestepsNoZeroAll denotes the set of all future time steps.
We denote by \granularityI{t} the time difference in seconds from timestep $t\in\setTimesteps$ to the initial state time step 0.

The gas flow on each arc $a=(\ell,r)$ from node $\ell$ to node $r$ at some time $t\in\setTimesteps$ is characterized by the pressure \press and the mass flow \mFlow at each end of the arc, i.e., the incoming pressure \pressI{\ell,a,t}, the outgoing pressure \pressI{r,a,t}, the incoming mass flow \mFlowI{\ell,a,t}, and the outgoing mass flow \mFlowI{r,a,t}.
Here, negative flow values represent flow against the arc's orientation, while the pressure values are always positive.
Since we assume that only pipes have a nonzero volume, the inflow and outflow of all non-pipe arcs is identical.
Thus, we define the mass flow of an arc $a\in\setArcs\setminus\setPipes$ as $\mFlowI{a,t} := \mFlowI{\ell,a,t} = \mFlowI{r,a,t}$.

The quantities of the arcs are connected at the nodes of the network.
For all arcs incident to a node, the corresponding pressure values at that end node are equal.
Hence, we can reference this value as pressure $\pressI{v,t}$ of node $v\in\setVertices$ for time $t\in\setTimesteps$ and it holds:
\begin{align*}
    \pressI{v,t}&= \pressI{v,a,t} \quad \forall a\in\setArcs : a=(v,r) \lor a=(\ell,v).
\end{align*}
For the mass flow holds the conservation of mass, a Kirchhoff-type law which demands that the flow into and out of a node $v\in\setVertices$ must be balanced for all time steps $t\in\setTimestepsNoZero$, i.e.,
\begin{align}
    \sum_{a\in\setArcs : a=(v,r)} \mFlowI{v,a,t} - \sum_{a\in\setArcs : a=(\ell,v)} \mFlowI{v,a,t} &= \inflowI{v,t}. \label{eq:flow_balance}
\end{align}
This equation is also known as \emph{flow balance equation}.
Note that the inflow \inflowI{v,t} denotes gas injection into node $v$ at time $t$ from outside of the network for positive values and gas withdrawal for negative values.
Furthermore, the inflow value at time $t$ represents the average gas inflow rate into the node in the continuous time interval between $t$-1 and $t$.
We partition the nodes into three sets: The set of \emph{entries} \setEntries with $\inflowI{v,t}\geq 0$ for $v\in\setEntries$, the set of \emph{exits} \setExits with $\inflowI{v,t}\leq 0$ for $v\in\setExits$, and the set of inner nodes \setInnerNodes with $\inflowI{v,t}=0$ for $v\in\setInnerNodes$ for all $t\in\setTimestepsNoZero$.
We call the union of all entries and exits the \emph{boundary nodes} \setBoundaryNodes of the network.

For each pressure and flow variable, we are given lower and upper bounds for each $t\in\setTimestepsNoZero$, which are denoted by bars below and above the corresponding variable.

Apart from the network topology, we are also given the initial state of the network, which specifies for $t=0$ the node pressure \pressI{v,0} for nodes $v\in\setVertices$ and the arc flow values \mFlowI{l,a,0} and \mFlowI{r,a,0} for arcs $a=(l,r)\in\setArcs$.
Furthermore, the initial states of the active network elements outside the network stations, i.e., valves and regulators, are given.

%% single element models
\subsection{Pipes}\label{sec:pipes}
For pipes, we use the \emph{friction-dominated} model given as (ISO3) in \citep{DomHilLanMeh2021}. It is based on the \emph{Euler equations} modeling the one-dimensional gas flow in straight cylindric pipes and can be derived by assuming a time-constant temperature (isothermal conditions), gas mixture, and compressibility factor, as well as clearly subsonic gas velocities and a negligible inertia term, which both are usually given in real-world conditions, see \citep{Osi1996,DomHilLanMeh2021} and \citep{Hen2021}.
The usage of the \emph{implicit box} scheme of \citep{DomGeiKolLan2011,KolLanBal2010} as a discretization yields the model for a pipe $a=(\ell,r)\in\setPipes$ and two adjacent time points $t_1\in\setTimestepsNoZero$ and $t_0:=t_1-1$ given as
\begin{subequations}\label{eq:pipe_model}
\begin{align}
    \pressI{\ell,t_1} + \pressI{r,t_1} - \pressI{\ell,t_0} - \pressI{r,t_0} + \frac{2 \sGasConst \tempI{a} \zFactI{a} (\granularityI{t_1}-\granularityI{t_0})}{\lenI{a} \areaI{a}}\left(\mFlowI{r,a,t_1} - \mFlowI{\ell,a,t_1}\right) &=0 \label{eq:fd_discrete_continuity}\\
    \pressI{r,t_1} - \pressI{\ell,t_1} +
    \underbrace{\frac{\fricI{a} \sGasConst \tempI{a} \zFactI{a} \lenI{a}}{4 \areaI{a}^2 \diamI{a}} \left(\frac{|\mFlowI{\ell,a,t_1}| \mFlowI{\ell,a,t_1}}{\pressI{\ell,t_1}} + \frac{|\mFlowI{r,a,t_1}| \mFlowI{r,a,t_1}}{\pressI{r,t_1}}\right)}_{\text{friction term}}
    + \frac{\gravAcc \slopeI{a} \lenI{a}}{2 \sGasConst \tempI{a} \zFactI{a}} \left(\pressI{\ell,t_1} + \pressI{r,t_1}\right) &=0.\label{eq:fd_discrete_momentum}
\end{align}
\end{subequations}
The first equation \eqref{eq:fd_discrete_continuity} is called the \emph{Continuity} equation, while the second equation \eqref{eq:fd_discrete_momentum} is known as the \emph{Momentum} equation and contains the non-linear \emph{friction term}.
In the equations, the specific gas constant is given by \sGasConst, the gravitational acceleration by \gravAcc, the gas temperature of the pipe by \tempI{a}, its compressibility factor by \zFactI{a}, its length by \lenI{a}, its diameter by \diamI{a}, its area by \areaI{a} defined as $\areaI{a}=\diamI{a}^2\pi/4$, its Darcy friction coefficient by \fricI{a}, and finally its slope by \slopeI{a} defined as $\slopeI{a} = (\heightI{r}-\heightI{\ell})/\lenI{a}$ using the height or altitude \height of the pipes' end nodes.
The gas temperature depends on the individual pipe, as we determine it based on the given gas temperatures \tempI{v,0} in the initial state for each node $v$ as
\begin{equation*}
    \tempI{a}=\frac{\tempI{\ell,0}+\tempI{r,0}}{2}.
\end{equation*}
In a similar fashion, we determine the pipes' compressibility factor using the formula of P\'apay \citep{Pap1968,Tak1989} based on the initial pressure and temperature of each node as
\begin{equation*}
    \zFactI{a} = 0.5(\zFact(\pressI{\ell,0},\tempI{\ell,0}) + \zFact(\pressI{r,0},\tempI{r,0})).
\end{equation*}
Finally, the Darcy friction factor is determined using the formula of Nikuradse \citep{Nik1950}, which is an explicit approximation depending only on the pipe's diameter and its integral roughness \roughI{a}.

\subsection{Valves}\label{sec:valves}
Valves are active elements that are able to dynamically connect or disconnect two network nodes.
Their state is represented by a binary mode variable \modeOp{a,t}, which specifies if the valve $a$ is open ($\modeOp{a,t}=1$) at time $t$, connecting both end nodes and forcing their pressures variables to attain identical values, or closed ($\modeOp{a,t}=0$), disconnecting the nodes, decoupling the respective pressure values, and prohibiting gas flow exchange.
Thus, the corresponding model for a valve $a=(\ell,r)\in\setValves$ reads
\begin{subequations}\label{eq:valve_model}
\begin{align}
    \modeOp{a,t} = 1 \quad&\implies\quad \pressI{\ell,t} = \pressI{r,t} && \forall t\in\setTimestepsNoZero\\
    \modeOp{a,t} = 0 \quad&\implies\quad \mFlowI{a,t} = 0 && \forall t\in\setTimestepsNoZero\\
    \modeOp{a,t} &\in\{0,1\} && \forall t\in\setTimesteps. \notag
\end{align}
\end{subequations}

\subsection{Regulators}\label{sec:regulators}
Regulators are active elements that can be seen as an extension of valves.
In addition to changing the network's connectivity, they can reduce the pressure of the gas in the direction of the flow.
However, this is only possible along the arc's orientation.
We use the regulator model of \citep{HenAndHopTur2021} here, in which the state of a regulator $a$ at time $t$ is modeled by a set of three different binary mode variables:
The two valve state variables of an open regulator \modeOp{a,t} and a closed regulator \modeCl{a,t} as well as the active mode \modeAc{a,t}, which allows for pressure reduction.
The model for all $a=(\ell,r)\in\setRegulators$ can be stated as
\begin{subequations}\label{eq:regulator_model}
\begin{align}
    1 = \modeAc{a,t} &+ \modeOp{a,t} + \modeCl{a,t} \label{eq:regulator_mode_choice}  && \forall t\in\setTimestepsNoZero\\
    \modeOp{a,t} = 1 \quad&\implies\quad \pressI{\ell,t} \leq \pressI{r,t}  && \forall t\in\setTimestepsNoZero\\
    \modeAc{a,t} + \modeOp{a,t} = 1 \quad&\implies\quad \pressI{\ell,t} \geq \pressI{r,t} \label{eq:regulator_pl_ge_pr} && \forall t\in\setTimestepsNoZero\\
    \modeCl{a,t} = 1 \quad&\implies\quad \mFlowI{a,t} \leq 0  && \forall t\in\setTimestepsNoZero\\
    \mFlowI{a,t} &\geq 0  && \forall t\in\setTimestepsNoZero\\
    \modeOp{a,t},\modeCl{a,t},\modeAc{a,t} &\in\{0,1\} && \forall t\in\setTimesteps. \notag
\end{align}
\end{subequations}
Note that in the indicator condition of Equation~\eqref{eq:regulator_pl_ge_pr}, the term $\modeAc{a} + \modeOp{a}$ acts as a binary variable, as Equation~\eqref{eq:regulator_mode_choice} forces the choice of exactly one of the three modes.
Furthermore, we assume an internal check valve in the element, which prevents flow against the orientation even in the case of an open regulator.
For a more detailed description of the internal regulator behavior, we refer to \citep{HenPetStr2021}.

\subsection{Artificial station arcs}\label{sec:auxArcs}
As stated above, the original network topology inside the stations is replaced by artificial arcs of different types, which represent an approximation of the station's gas transport capabilities.
Each arc can be active or inactive, which we model by a binary variable
\begin{equation*}
    \activeArcI{a,t}\in\{0,1\} \qquad \forall a\in\setArtificialLinks \quad \forall t\in\setTimesteps,
\end{equation*}
where $\activeArcI{a,t} = 1$ represents the active state.
An inactive arc behaves like a closed valve, i.e., the flow is zero and the end nodes' pressures are decoupled, while the active behavior is specific to each artificial arc type, which we discuss in the following.

\subsubsection{Shortcuts}
Shortcuts connect its end nodes if active.
Hence, their model is equivalent to that of a valve and can be stated for a shortcut $a=(\ell,r)\in\setShortcuts$ and time $t\in\setTimestepsNoZero$ as:
\begin{subequations}\label{eq:artificial_shortcut_model}
\begin{align}
    \activeArcI{a,t} = 1 \quad&\implies\quad \pressI{\ell,t} = \pressI{r,t} \\
    \activeArcI{a,t} = 0 \quad&\implies\quad \mFlowI{a,t} = 0.
\end{align}
\end{subequations}

\subsubsection{Regulating arcs}
Artificial regulating arcs capture the general behavior of regulators.
Hence, they only allow flow along their orientation and offer the possibility of pressure reduction if they are active.
The corresponding model for a regulating arc $a=(\ell,r)\in\setRegulatorArcs$ and time $t\in\setTimestepsNoZero$ can be stated as
\begin{subequations}\label{eq:artificial_regulator_model}
\begin{align}
    \activeArcI{a,t} = 1 \quad&\implies\quad \pressI{\ell,t} \geq \pressI{r,t} \\
    \activeArcI{a,t} = 0 \quad&\implies\quad \mFlowI{a,t} \leq 0 \\
    \mFlowI{a,t} &\geq 0.
\end{align}
\end{subequations}
Note that the Equations~\eqref{eq:regulator_model} do not apply here.

\subsubsection{Bi-regulating arcs}\label{sec:auxBiDirRegulators}
To represent the behavior of a pair of anti-parallel regulators in the original network topology, we use a bi-regulating arc.
We can model its behavior by choosing an orientation of the element in active mode and applying the corresponding constraints of a regulating arc for each direction separately.
The model for a bi-regulating arc $a=(\ell,r)\in\setBiRegulatorArcs$ and time $t\in\setTimestepsNoZero$ is given as
\begin{subequations}\label{eq:artificial_biregulator_model}
\begin{align}
    \activeArcI{a,t} &= \activeFwdI{a,t} + \activeBwdI{a,t} \\
    \mFlowI{a,t} &= \mFlowFwdI{a,t} - \mFlowBwdI{a,t} \\
    \activeFwdI{a,t} = 1 \quad&\implies\quad \pressI{\ell,t} \geq \pressI{r,t} \\
    \activeFwdI{a,t} = 0 \quad&\implies\quad \mFlowFwdI{a,t} \leq 0 \\
    \mFlowFwdI{a,t} &\geq 0 \\
    \activeBwdI{a,t} = 1 \quad&\implies\quad \pressI{r,t} \geq \pressI{\ell,t} \\
    \activeBwdI{a,t} = 0 \quad&\implies\quad \mFlowBwdI{a,t} \leq 0 \\
    \mFlowBwdI{a,t} &\geq 0 \\
    \activeFwdI{a,t},\activeBwdI{a,t} &\in\{0,1\}. \notag
\end{align}
\end{subequations}
Here, \activeFwdI{a,t} represents the choice for the \emph{forward orientation} for arc $a$ and time $t$, in which flow and regulation happen along the arc's orientation, while \activeBwdI{a,t} represents the choice of the \emph{backward orientation}, in which the orientation is flipped.
In the same fashion, \mFlowFwdI{a,t} represents the forward flow while \mFlowBwdI{a,t} represents the flow in the backward direction.

\subsubsection{Compressor arcs}\label{sec:auxCompressor}
Finally, there are the compressor arcs, which are the only elements able to increase the gas pressure in the direction of the flow.
They approximate the capabilities of compressor station arcs in the original network station topology.
For this reason, we shortly introduce these in the following.

Each compressor station, which we assume to be modeled as an arc $a=(\ell,r)$ as well, represents a combination of single compressor units, which again are a combination of a compressor machine and a corresponding drive that provides the necessary power.
The compressor machine of a compressor unit $u$ can operate in its feasible operating range, defined as a polytope in the space $(\vFlow, \pRatio)$ of the volumetric flow \vFlow and the compression ratio \pRatio given as
\begin{subequations}
\begin{align}
    p &= \dens \sGasConst \tempI{a} \zFactI{a} \label{eq:state_equation_of_real_gases} \\
    \vFlow = \frac{\mFlow}{\densI{\ell}} \quad&\implies \quad \vFlow = \mFlow \frac{\sGasConst  \tempI{a} \zFactI{a}}{\pressI{\ell}} \label{eq:volumetric_flow_definition} \\
    \pRatio &= \frac{\pressI{r}}{\pressI{\ell}}.
\end{align}
\end{subequations}
Here, \dens denotes the density of the gas, \densI{\ell} the incoming gas density, while \zFactI{a} and \tempI{a} are defined in the same way as they were for pipes in Section~\ref{sec:pipes}.
Equation~\eqref{eq:state_equation_of_real_gases} is called the equation of state for real gases \citep{Osi1996,FugGeiGolMor2015}.

Furthermore, the compression is restricted by the maximal amount of power $\ub{P}_u$ the corresponding drive can provide.
The power needed for compressing a mass flow value of \mFlow from \pressI{\ell} to \pressI{r}, i.e., with $\pRatio = \pressI{r}/\pressI{\ell}$, is given as
\begin{equation}
    \power(\mFlow,\pressI{\ell},\pressI{r}) = \frac{\mFlow}{\adEff}\sGasConst \tempI{a} \zFactI{a} \frac{\isenExp}{\isenExp-1} \left[ \pRatio^{\frac{\isenExp-1}{\isenExp}} - 1 \right].
  \label{eq:power_definition}
\end{equation}
Here, \adEff is the adiabatic efficiency and is assumed to be a constant given per compressor unit.
Further, \isenExp is the isentropic exponent, which we approximate by the constant value of $1.296$, see, for example, \citep{FugGeiGolMor2015} for more details.
The single units can then be combined to predefined configurations $\setCompressorConfigurations{a}$ of the compressor station $a$, defined as serial connections of parallel unit arrangements.
For a more in-depth description, we refer to~\citep{FugGeiGolMor2015,WalHil2017,HenAndHopTur2021}.

In the base model of \citep{HopHenLenKoc2021}, the artificial pressure increasing arcs were modeled as a single compressor representing the combination of single units, which could dynamically be assigned to pressure-increasing arcs using binary assignment variables.
As a consequence of this dynamic combination of different units, some of the constraints describing the resulting compression capabilities were approximated quite generously to keep the overall model linear.

To avoid this problem, we develop a novel formulation, which is closer to the original compressor station model while the model size with respect to the number of variables and constraints is not increased.
Like for the compressor stations, we create for each artificial compressor arc $a=(\ell,r)\in\setCompressorArcs$ a set of configurations \setAuxCompressorConfigsI{a} to choose from.
These form different compression settings based on how the compressor units can be composed in the underlying compressor stations.
The feasible operating range of a chosen configuration $c\in\setAuxCompressorConfigsI{a}$ of an artificial compressor arc $a$ is restricted in terms of three parameters:
The maximum compression ratio $\pRatioI{c}$, the maximum volumetric flow $\vFlowI{c}$, and the maximum compression power $\powerI{c}$.
Furthermore, each configuration has an associated number $u_c$ of used compressor units.
From the parameters restricting the feasible operating range, we further find parameters $\alpha^\mathrm{\pressI{\ell}}_c$, $\alpha^\mathrm{\pressI{r}}_c$, and $\alpha^\mathrm{\mFlow}_c$, which represent a linearization of the compression power of the form
\begin{equation*}
    P \approx \alpha^\mathrm{\pressI{\ell}}_c \pressI{\ell} + \alpha^\mathrm{\pressI{r}}_c \pressI{r} + \alpha^\mathrm{\mFlow}_c \mFlow,
\end{equation*}
and are determined by a linear regression based on triples $(\pressI{\ell}^s,\pressI{r}^s,\mFlow^s)_{s\in\mathbb{N}}$ sampled from within the intersection of the corresponding variable bound intervals over time, which satisfy the constraints
\begin{align*}
    \mFlow^s \frac{\sGasConst  \tempI{a} \zFactI{a}}{\pressI{\ell}^s} &\leq \vFlowI{c} &
    \frac{\pressI{r}^s}{\pressI{\ell}^s} &\leq \pRatioI{c} &
    \power(\mFlow^s,\pressI{\ell}^s,\pressI{r}^s) &= \powerI{c}.
\end{align*}
Note that usually, $\alpha^\mathrm{\pressI{\ell}}_c < 0$, $\alpha^\mathrm{\pressI{r}}_c > 0$, and $\alpha^\mathrm{\mFlow}_c > 0$ hold.

The final component of the configuration-based model is the set of artificial compressor configuration conflicts \setAuxCompressorConfigConflicts.
Each conflict $z$ is defined as tuple $z=(c_1,c_2)$ of configurations of different artificial compressors from the same network station that cannot be used simultaneously.
Conflicts are needed to model, for example, compressor units that can be used in different compressor stations, but not at the same time.
As a consequence, configurations $c_1$ and $c_2$ of two corresponding artificial compressors cannot be used simultaneously if both can use the same unit.
Since conflicts are always defined for configurations of artificial compressors in the same network station, we can partition the conflicts per network station as $\setAuxCompressorConfigConflicts = \dotbigcup_{i\in\setNetworkStations} \setAuxCompressorConfigConflictsI{i}$.

Given this description, we state the configuration-based model for an artificial compressor $a=(\ell,r)\in\setCompressorArcs$ and time $t\in\setTimestepsNoZero$ as
\begin{subequations}\label{eq:artificial_compressor_model}
\begin{align}
    \activeArcI{a,t} &= \sum_{c\in\setAuxCompressorConfigsI{a}} \activeCfgI{c,t} \label{eq:one_cfg_per_active_compressor}\\
    \activeArcI{a,t} = 1 \quad&\implies\quad \pressI{\ell,t} \leq \pressI{r,t} \label{eq:auxCompressor_active_pressure_increase}\\
    \activeArcI{a,t} = 0 \quad&\implies\quad \mFlowI{a,t} \leq 0\\
    \mFlowI{a,t} &\geq 0, \label{eq:auxCompressor_positive_flow}
\end{align}
\end{subequations}
where we additionally demand for each configuration $c\in\setAuxCompressorConfigsI{a}$ of an artificial compressor $a=(\ell,r)\in\setCompressorArcs$ that
\begin{subequations}\label{eq:artificial_compressor_config_model}
\begin{align}
    \activeCfgI{c,t} = 1 \quad&\implies\quad \pressI{r,t} \leq \pressI{\ell,t} \pRatioI{c} && \forall t\in\setTimestepsNoZero\\
    \activeCfgI{c,t} = 1 \quad&\implies\quad \mFlowI{a,t}    \leq \pressI{\ell,t} \frac{\vFlowI{c}}{\sGasConst \tempI{a} \zFactI{a}} && \forall t\in\setTimestepsNoZero \\
    \activeCfgI{c,t} = 1 \quad&\implies\quad \alpha^\mathrm{\pressI{\ell}}_c \pressI{\ell,t} + \alpha^\mathrm{\pressI{r}}_c \pressI{r,t} + \alpha^\mathrm{\mFlow}_c \mFlowI{a,t} \leq \powerI{c} && \forall t\in\setTimestepsNoZero \label{eq:maxPowerOfConfig}\\
    \activeCfgI{c,t} &\in\{0,1\} && \forall t\in\setTimesteps. \notag
\end{align}
\end{subequations}
Here, the binary variable \activeCfgI{c,t} represents the decision to use (value=1) or not to use (value=0) the corresponding configuration $c\in\setAuxCompressorConfigsI{a}$ at time $t$.
In addition, we need the following constraint for each conflict $z=(c_1,c_2)\in\setAuxCompressorConfigConflicts$ and time $t\in\setTimestepsNoZero$:
\begin{equation}
    \activeCfgI{c_1,t} + \activeCfgI{c_2,t} \leq 1. \label{eq:compressor_conf_conflict}
\end{equation}
Note that the artificial compressor arcs are directed and therefore can only have positive flow as well as compress gas from the end node $\ell$ towards $r$, which is ensured via Equations~\eqref{eq:auxCompressor_positive_flow} and \eqref{eq:auxCompressor_active_pressure_increase}.

\paragraph{Artificial configuration set construction}
Since the artificial compressor configurations are not given as part of the input, we create them for each artificial compressor $a\in\setCompressorArcs$ based on the original compressor stations represented by it.
Here, we associate each artificial configuration $c$ with a tuple $(u^p_c, u^s_c, U_c)$, where $u^p_c$ is the number of parallel compressor units to use, $u^s_c$ is the number of serial units to use, and $U_c$ is the set of usable units.
We manually create these tuples based on the configurations of the original compressor stations as well as their topological position in the network.
If multiple original compressor compositions use similar units in the same layout, we represent them by the same artificial configuration if possible.
An example of the construction is given in Figure~\ref{fig:aux_configuration_construction}.

\begin{figure}[t]
  \newcommand{\nodeSquare}[2]{\draw[fill=white] (#1-0.05,#2-0.05) rectangle (#1+0.05,#2+0.05);}
  \centering
  \fbox{
  \begin{tikzpicture}
  [circuit,
   circuit symbol unit=2pt,
   scale=1.6, every node/.style={scale=1.0}]
  % some nice explanation of styles https://tex.stackexchange.com/a/52379
  \draw (0,0) to [cs] (1,0); %u1
  \draw (2,0) to [cs] (3,0); %u2

  % u1 || u3
  \draw (4.2, 0.2) to [cs] (4.8, 0.2);
  \draw (4.2,-0.2) to [cs] (4.8,-0.2);
  \draw (4.0, 0.0) to (4.2, 0.2);
  \draw (4.0, 0.0) to (4.2,-0.2);
  \draw (5.0, 0.0) to (4.8, 0.2);
  \draw (5.0, 0.0) to (4.8,-0.2);

  % u2 || u3
  \draw (6.2, 0.2) to [cs] (6.8, 0.2);
  \draw (6.2,-0.2) to [cs] (6.8,-0.2);
  \draw (6.0, 0.0) to (6.2, 0.2);
  \draw (6.0, 0.0) to (6.2,-0.2);
  \draw (7.0, 0.0) to (6.8, 0.2);
  \draw (7.0, 0.0) to (6.8,-0.2);

  \nodeSquare{0}{0}
  \nodeSquare{1}{0}
  \nodeSquare{2}{0}
  \nodeSquare{3}{0}
  \nodeSquare{4}{0}
  \nodeSquare{5}{0}
  \nodeSquare{6}{0}
  \nodeSquare{7}{0}
  % label
  \node at (0.5,-0.2) {$u_1$};
  \node at (2.5,-0.2) {$u_2$};
  \node at (4.5, 0.0) {$u_1$};
  \node at (4.5,-0.4) {$u_3$};
  \node at (6.5, 0.0) {$u_2$};
  \node at (6.5,-0.4) {$u_3$};
  % braces
  \draw [decorate,decoration={brace}] (-0.15,-0.3) --  (-0.15, 0.3);
  \draw [decorate,decoration={brace}] (7.15, 0.3) --  (7.15,-0.3);
  \draw [decorate,decoration={brace}] (2.5,-2.1) --  (2.5, -1.1);
  \draw [decorate,decoration={brace}] (4.2,-1.1) --  (4.2, -2.1);
  % text
  \node[align=left,anchor=south west] at (-1,0.4) {Assume a artificial compressor $a_{\mathrm{ac}}$ with a single associated compressor station $a_{\mathrm{cs}}$.\\ Let the configurations of $a_{\mathrm{cs}}$ be defined as};
  \node[align=left,anchor=south west] at (-1,-0.2) {$\mathcal{C}_{a_{\mathrm{cs}}} :=$};
  \node[align=left,anchor=south west] at (7.2, -0.2) {.};
  \node[align=left,anchor=south west] at (-1,-0.9) {Then the corresponding artificial configurations for $a_{\mathrm{ac}}$ are};
  \node[align=left,anchor=south west] at (1.65,-1.8) {$\mathcal{C}^{\mathrm{ar}}_{a_{\mathrm{ac}}} :=$};
  \node[align=left,anchor=south west] at (4.25,-1.8) {.};
  % configuration sets
  \node at (2.8,-1.0) {$u^p$};
  \node at (3.1,-1.0) {$u^s$};
  \node at (3.6,-1.0) {$U$};
  \node[align=left,anchor=south west] at (2.5,-1.5) {$(\,1\,,\,1\,,\{u_1,u_2\}),$};
  \node[align=left,anchor=south west] at (2.5,-1.8) {$(\,2\,,\,1\,,\{u_1,u_3\}),$};
  \node[align=left,anchor=south west] at (2.5,-2.1) {$(\,2\,,\,1\,,\{u_2,u_3\})$};
  % size
  \node at (3.5, 1.0) {};
  \node at (3.5,-2.1) {};
  \end{tikzpicture}}
  \caption{Example of the construction of the configurations for an artificial compressor $a_\mathrm{ac}$. The associated compressor station $a_\mathrm{cs}$ has four different configurations using three different compressor units: Two configurations using only a single unit, either unit $u_1$ or unit $u_2$, and two configurations in which each of the two units is used in parallel to unit $u_3$. There are three resulting artificial configurations for $a_\mathrm{ac}$, which all allow only one serial stage: One configuration with only one allowed parallel unit, in which $u_1$ or $u_2$ can be chosen, and two configurations allowing for two parallel units, one with the units set $\{u_1,u_3\}$ and one with the set $\{u_2,u_3\}$. While we can combine the first two original configurations into one artificial configuration, this is not possible for the parallel configurations. The corresponding unit set would have to be $\{u_1,u_2,u_3\}$, which is not feasible as there is no original configuration combining $u_1$ and $u_2$ in parallel.}
  \label{fig:aux_configuration_construction}
\end{figure}

\begin{figure}[t]
    \centering
    \includegraphics[width=0.5\textwidth]{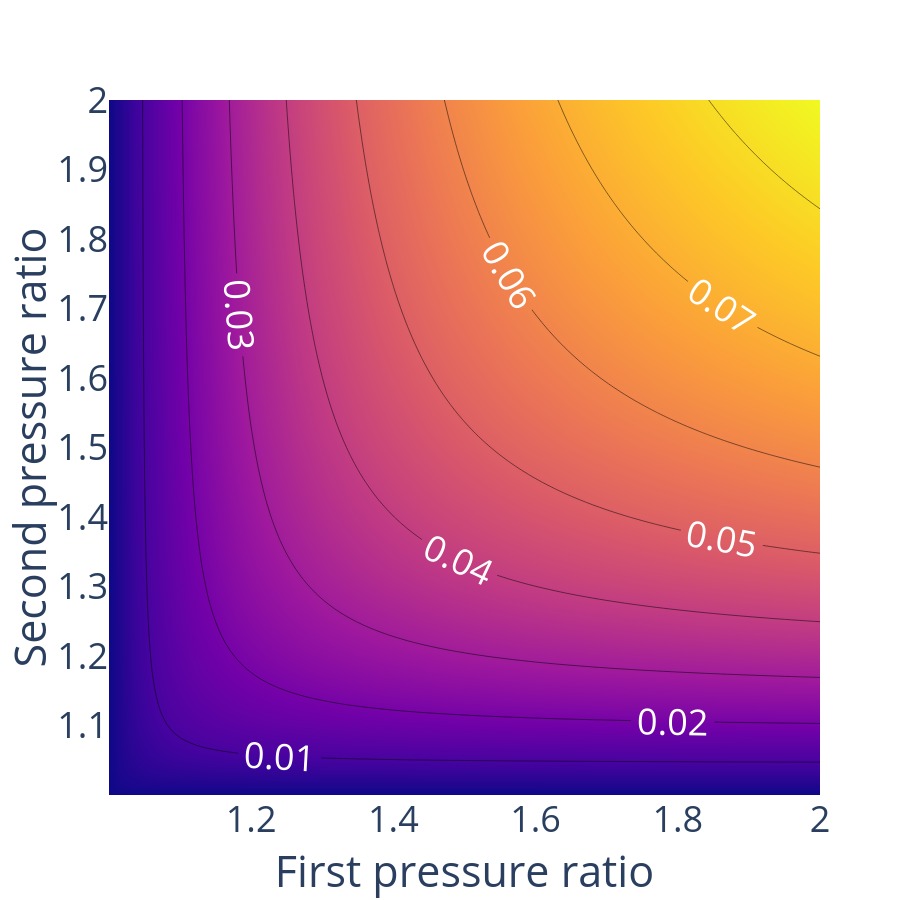}
    \caption{Relative error in compression power. We compare the power \powerI{\mathrm{sum}} needed for a serial compression using two separate compression processes with ratios \pRatioI{1} and \pRatioI{2} to the power \powerI{\mathrm{one}} needed by our approximation of using a single compression process with ratio $\pRatioI{1} \cdot \pRatioI{2}$. We plotted the relative error $(\powerI{\mathrm{one}}-\powerI{\mathrm{sum}})/\powerI{\mathrm{sum}}$ for compression ratios $\pRatioI{1}$ and $\pRatioI{2}$ up to $2.0$. Figure created with \citep{Plo2015}.}
    \label{fig:contour_powerSum_RelError}
\end{figure}

Given $(u^p_c, u^s_c, U_c)$ for each artificial configuration $c$, we now derive the corresponding parameters for the number of simultaneously used compressor units $u_c$, the maximum compression ratio $\pRatioI{c}$, the maximum volumetric flow $\vFlowI{c}$, and the maximum power $\powerI{c}$.
First, we present a procedure for a configuration $c$ for the special case of either $u^p_c$ or $u^s_c$ being equal to 1, i.e., for the case of a pure serial or pure parallel compression.
For each compressor unit $u\in U_c$, we are given the maximum power $\ub{\power}_u$ its associated drive can provide.
Furthermore, we can determine its maximum compression ratio $\ub{\pRatio}_u$ and maximum volumetric flow $\ub{\vFlow}_u$ by iterating the extreme points of its feasible operating range polytope.
Let then $(\tilde{\pRatio}_n)_{n=1}^{|U_c|}$ be the sequence of maximum compression ratios $\ub{\pRatio}_u$ for all $u\in U_c$ sorted in descending order, and analogously let $(\tilde{\vFlow}_n)_{n=1}^{|U_c|}$ and $(\tilde{\power}_n)_{n=1}^{|U_c|}$ be the corresponding descending sequences of maximum volumetric flows $\ub{\vFlow}_u$ and power values $\ub{\power}_u$.
We then define the parameters $u_c$, \pRatioI{c}, \vFlowI{c}, and \powerI{c} for $c\in\setAuxCompressorConfigsI{a}$ of $a\in\setCompressorArcs$ as
\begin{subequations}\label{eq:cfg_param_construction_pure_serial_or_parallel}
\begin{align}
    u_c &= u^p_c\cdot u^s_c \\
    \pRatioI{c} &=
    \begin{cases}
      \prod_{i=1}^{u^s_c} \tilde{\pRatio}_i & \text{ if } u^s_c > 1 \\
      \tilde{\pRatio}_{u^p_c} & \text{ if } u^s_c = 1,
    \end{cases}\label{eq:auxConfig_ratio_easy_case}\\
    \vFlowI{c} &=
    \begin{cases}
      \sum_{i=1}^{u^p_c} \tilde{Q}_i & \text{ if } u^p_c > 1 \\
      \tilde{Q}_{u^s_c} & \text{ if } u^p_c = 1,
    \end{cases}\label{eq:auxConfig_flow_easy_case}\\
    \powerI{c} &= \sum_{i=1}^{u_c} \tilde{\power}_i.
\end{align}
\end{subequations}
Note that for the serial compression, summing up the power values is not quite accurate and slightly underestimates the necessary compression power when comparing the original two serial compression processes with one combined compression process using the product of the single compression ratios, see Figure~\ref{fig:contour_powerSum_RelError}.
As a consequence, $P_c$ is potentially too small to allow for all compression processes possible in the original network.
However, in the data used for our computational experiments in Section~\ref{sec:computational_experiments}, the maximum serial compression number $u^s_c$ for a configuration $c$ is 2.
Furthermore, all involved units have compressor ratios smaller than 2.0.
Thus, the approximation error is guaranteed to be smaller than 9\%.

If a configuration features multiple serial and multiple parallel stages, the above-presented formulas for creating the configuration parameters cannot be applied directly.
Instead,
we create for each such arrangement and each concrete assignment of a compressor unit to a position in the arrangement a separate artificial configuration for the corresponding artificial compressor arc.
The total number of used compressor units $u_c$ is then equal to the number of active units in each specific arrangement, and the maximum power $\powerI{c}$ is determined as the sum of maximum power values of these units.
To determine the parameters \pRatioI{c} and \vFlowI{c}, we recursively apply the formulas used in Equations~\eqref{eq:auxConfig_ratio_easy_case} and \eqref{eq:auxConfig_flow_easy_case} for each of the concrete serial and parallel sub-compositions of units.
The approach of creating one configuration for each possible unit arrangement can, in theory, result in a lot of artificial configurations.
However, in our experience, this is not the case in practice, as such configurations do not occur very often.

\paragraph{Model comparison}
The main benefit of the new artificial compressor model in Equations~\eqref{eq:artificial_compressor_model}, \eqref{eq:artificial_compressor_config_model}, and \eqref{eq:compressor_conf_conflict} compared to the base model is that each compressor unit contributes either in serial or in parallel to the compression capabilities and does not use the benefits both at the same time.

Furthermore, we improved the description of the maximum compression ratio \pRatio and the maximum volumetric flow \vFlow.
For both, we do not the need to replace the incoming pressure of the arc \pressI{\ell} by a constant value, which was done explicitly in Equation (40) in \citep{HopHenLenKoc2021} and implicitly by assuming a maximum mass flow bound per compressor unit, which can only be derived from the inherent volumetric flow bound by assuming a constant incoming pressure.
Furthermore, we can correctly employ the product of compressor ratios in the computation of the maximum usable compression ratio of serial compression processes in
Equation~\eqref{eq:auxConfig_ratio_easy_case}
instead of the linear approximation formula used in (35) in \citep{HopHenLenKoc2021}.

\begin{figure}[thp]
  \centering

  \begin{subfigure}[c]{0.48\textwidth}
  % trim options: <left> <lower> <right> <upper>
  \includegraphics[trim={3.8cm 4.5cm 8.0cm 4.0cm},clip,width=0.9\textwidth]{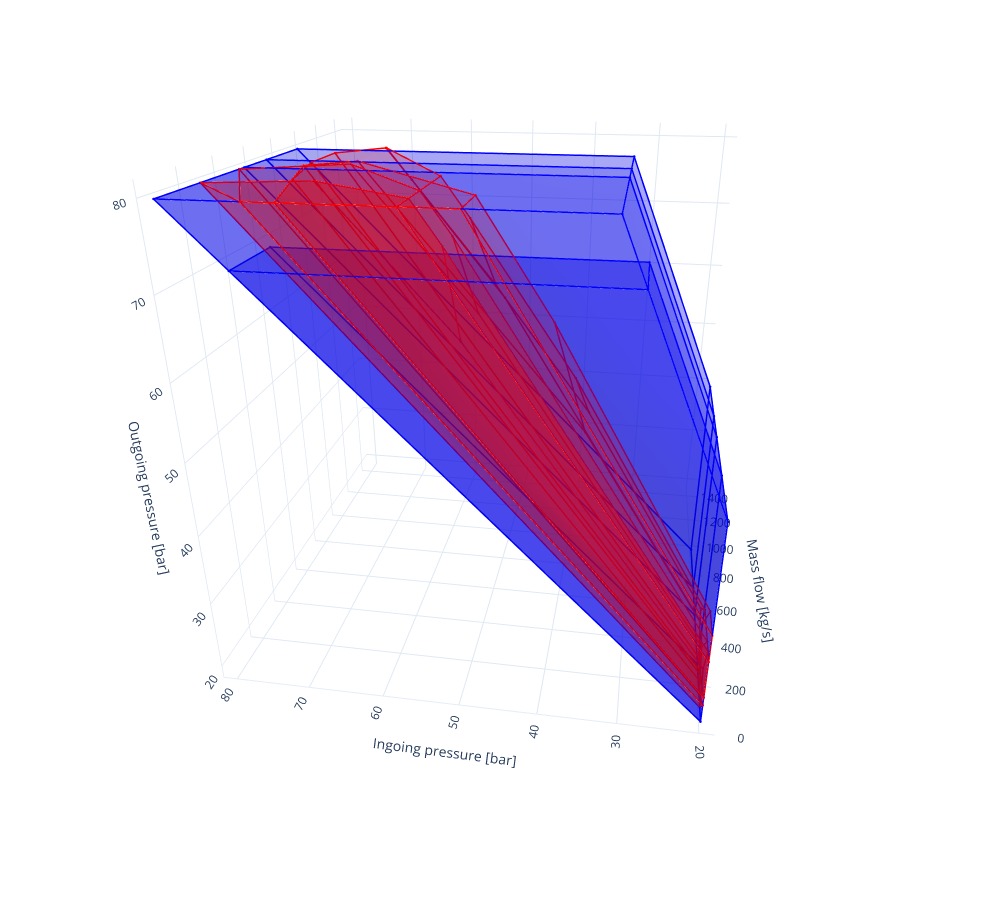}
  \subcaption{Base model, View 1}
  \label{fig:baseModel_cam1}
  \end{subfigure}
  \begin{subfigure}[c]{0.48\textwidth}
  \includegraphics[trim={3.8cm 4.5cm 8.0cm 4.0cm},clip,width=0.9\textwidth]{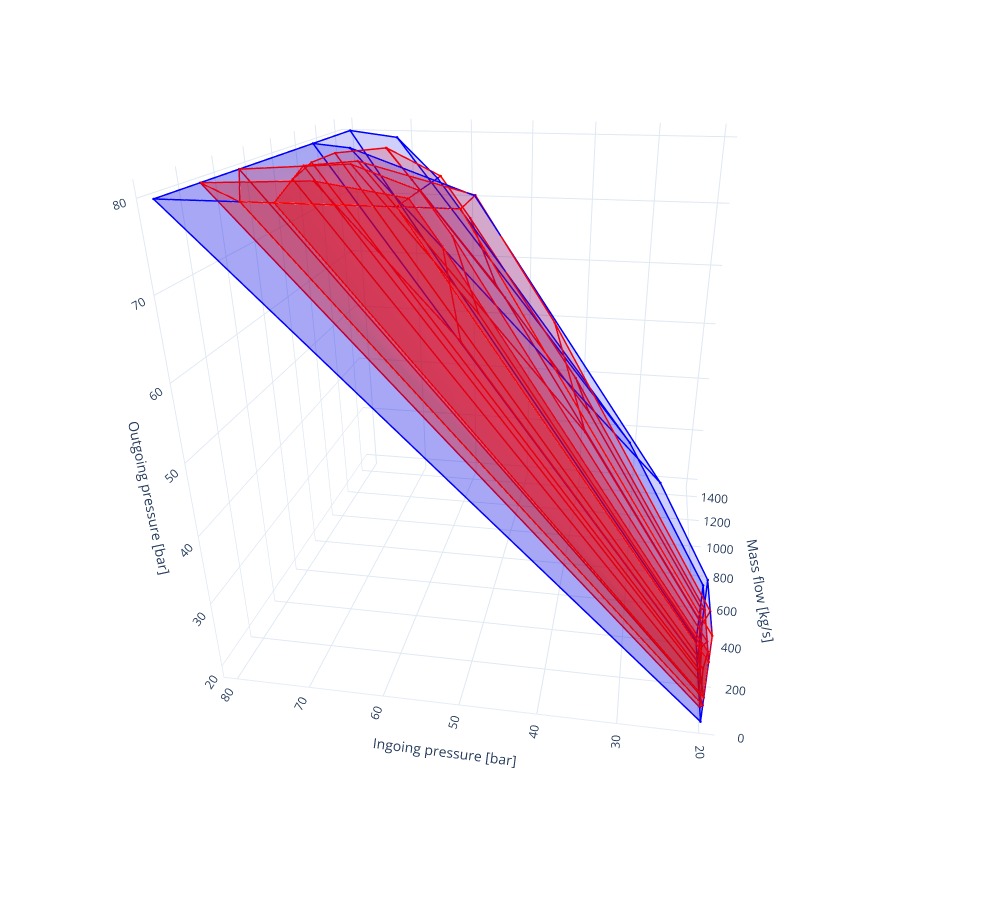}
  \subcaption{New model, View 1}
  \label{fig:newModel_cam1}
  \end{subfigure}

  \bigskip

  \begin{subfigure}[c]{0.48\textwidth}
  % trim options: <left> <lower> <right> <upper>
  \includegraphics[trim={6.0cm 5.0cm 6.0cm 5.0cm},clip,width=0.9\textwidth]{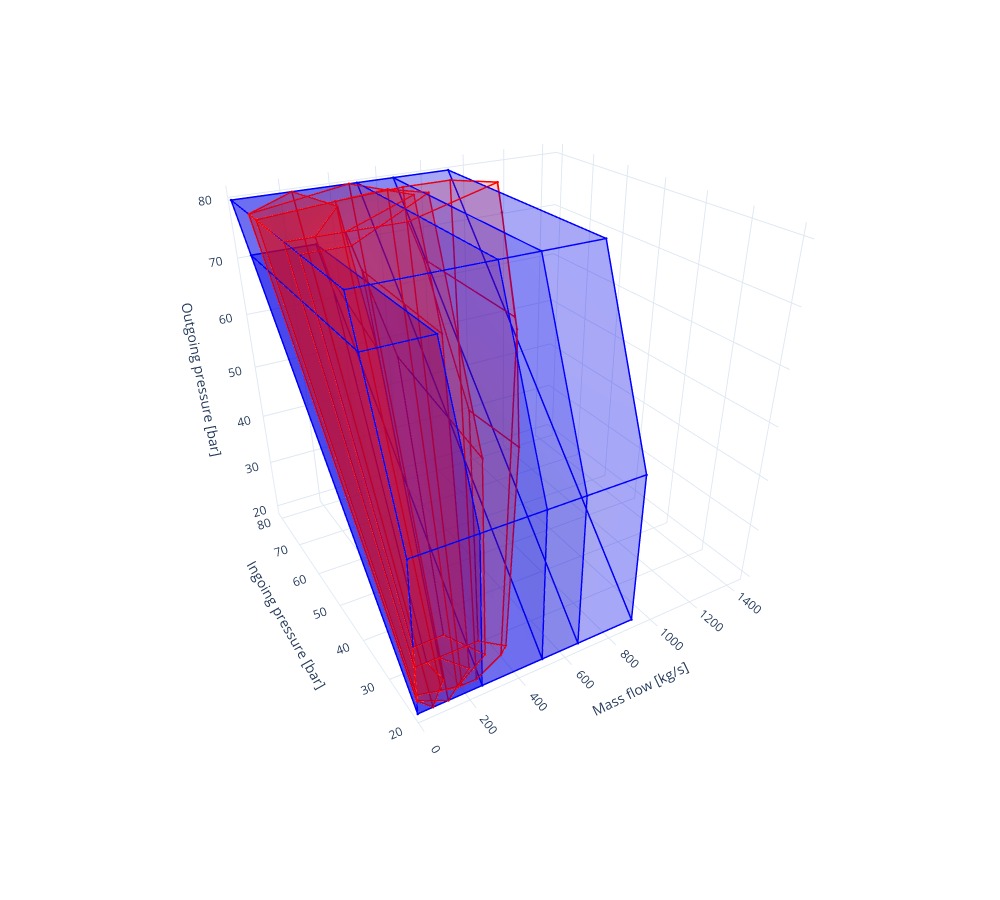}
  \subcaption{Base model, View 2}
  \label{fig:baseModel_cam2}
  \end{subfigure}
  \begin{subfigure}[c]{0.48\textwidth}
  \includegraphics[trim={6.0cm 5.0cm 6.0cm 5.0cm},clip,width=0.9\textwidth]{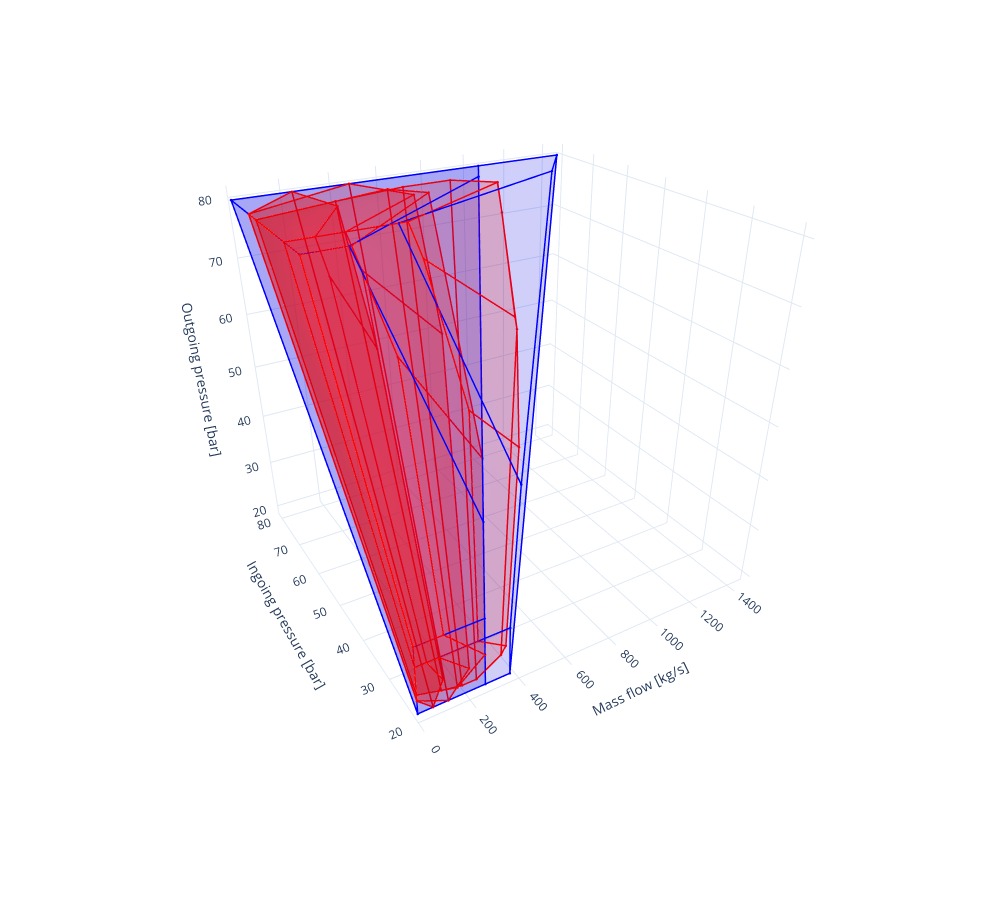}
  \subcaption{New model, View 2}
  \label{fig:newModel_cam2}
  \end{subfigure}
  \caption{Feasible region comparison of the artificial compressor base model against our new one at the example of one artificial compressor arc. As a reference solution, we plotted the feasible region polytopes of the configurations of the associated compressor station from the original network topology in red.
  The polytopes for the artificial compressor arc are plotted in blue, where Figures (a) and (c) show the base model polytopes and Figures (b) and (d) show the new model polytopes.
  For the new model, each polytope represents one artificial configuration to choose from.
  Since the base model does not feature configurations, we plotted one polytope for each feasible compressor unit combination on the artificial arc.
  The pressure bound interval used for both end nodes is $[20\,\text{bar}, 80\,\text{bar}]$.
  For the constant incoming pressure needed for the base model, we use the median value of $50\,\text{bar}$.
  For a better impression of the three-dimensional polytopes, we show them from two different camera positions in View 1, visible in Figures (a) and (b), and View 2, visible in Figures (c) and (d). Figures created with \citep{Plo2015}.}
  \label{fig:modelComparison_auxCompressor}
\end{figure}

As an example, we show the model improvements in comparison to \citep{HopHenLenKoc2021} for a single artificial compressor arc in Figure~\ref{fig:modelComparison_auxCompressor}, where we compare the feasible region polytopes produced by the base model, our new formulation, and the actual feasible region of the configuration polytopes of the compressor station in the original station topology.
To have a clean comparison, we choose an artificial compressor arc representing exactly one compressor station.
This compressor station features three compressor units, from which two units can be combined in parallel, but has no serial unit combination capabilities.
Note that we excluded the power bound from the feasible region polytopes since it is modeled for both models in the same way, i.e., as a linear approximation using a linear regression approach.

We clearly see the increased accuracy of the new model already for an artificial compressor arc, which is only able to combine two compressor units.
When combining more compressor units, the differences would further increase.

\subsection{Station simple states}
The stations' artificial arcs do not operate independently from each other.
Instead, each station $i\in\setNetworkStations$ is always in exactly one \emph{simple state}, where the set of possible simple states of the station is given as \setSimpleStatesI{i}, and the union of simple states over all stations is given as $\setSimpleStates = \bigcup_{i\in\setNetworkStations} \setSimpleStatesI{i}$.
Each simple state $s\in\setSimpleStatesI{i}$ restricts the usage of the artificial arcs by defining a set of arcs $\simpStOnArcsI{s}\subseteq\setArtificialLinksNaviStation{i}$, which need to be active if the station uses state $s$, and a set of arcs $\simpStOffArcsI{s}\subseteq\setArtificialLinksNaviStation{i}$, which cannot be active if the station uses state $s$.
For all other artificial arcs $\setArtificialLinksNaviStation{i}\setminus(\simpStOnArcsI{s} \cup \simpStOffArcsI{s})$, the activity is optional.
We represent the choice for a simple state $s$ at time $t$ by the binary variable \activeStateI{s,t}, where a value of 1 represents an active simple state.
Using this, we formulate the simple state model for each network station $i\in\setNetworkStations$ as
\begin{align}
    \sum_{s\in\setSimpleStatesI{i}} \activeStateI{s,t} &= 1 && \forall t\in\setTimestepsNoZero  \label{eq:one_simple_state} \\
    \activeStateI{s,t} &\in\{0,1\} && \forall s\in\setSimpleStatesI{i} \quad \forall t\in\setTimesteps \notag
\end{align}
and for all artificial arcs $a\in\setArtificialLinksNaviStation{i}$ of each network station $i\in\setNetworkStations$ and time $t\in\setTimestepsNoZero$ as
\begin{align}
    \sum_{s\in\setSimpleStatesI{i}: a\in\simpStOnArcsI{s}} \activeStateI{s,t} &\leq \activeArcI{a,t} &
    1 - \sum_{s\in\setSimpleStatesI{i}: a\in\simpStOffArcsI{s}} \activeStateI{s,t} &\geq \activeArcI{a,t}. \label{eq:auxArc_simple_state_connection_model}
\end{align}

\subsection{Station flow directions and fence node valves}\label{sec:flow_dirs_and_fence_node_valves}
Apart from the simple state, each station is also assigned exactly one \emph{flow direction} $f$ at each point in time $t\in\setTimesteps$, captured by the binary variable \activeFlowDirI{f,t}, where the value 1 again stands for an active flow direction.
The set of possible flow directions for station $i$ is given as \setFlowDirectionsI{i}.
A flow direction $f\in\setFlowDirectionsI{i}$ determines the connection between the interior and exterior of each station by defining two sets of fence nodes: Those fence nodes $\flowDirEntriesI{f}\subseteq\setFenceNodesNaviStation{i}$ allowing for flow into the station and those fence nodes $\flowDirExitsI{f}\subseteq\setFenceNodesNaviStation{i}$ allowing for flow out of the station.
All other fence nodes $\setFenceNodesNaviStation{i}\setminus(\flowDirEntriesI{f} \cup \flowDirExitsI{f})$ do not allow any exchange of flow between the station's interior and exterior.

\begin{figure}[t]
\centering

\newcommand{\nodeSquare}[2]{\draw[fill=white] (#1-0.05,#2-0.05) rectangle (#1+0.05,#2+0.05);}
\begin{tikzpicture}
  [circuit,
   circuit symbol unit=2pt,
   scale=1.6, every node/.style={scale=1.0}]
  % some nice explanation of styles https://tex.stackexchange.com/a/52379
  %%%%%%%%%%%%%%%%%%%%%%%%%% alpha paths %%%%%%%%%%%%%%%%%%%%%%%
  \draw (2,0) to [va] (3,0);
  \draw [dashed,thick,path fading=west] (2,0) -- (1, 1);
  \draw [dashed,thick,path fading=west] (2,0) -- (1,-1);
  \draw [dashed,thick,path fading=east] (3,0) -- (4, 1);
  \draw [dashed,thick,path fading=east] (3,0) -- (4,-1);
  % draw arc, see https://tex.stackexchange.com/a/66491
  % mid point: (6,0)
  % radius:    4
  % from angle 150 to 210 (180 +- 30)
  \draw [dotted] (6,0) ++(150:4) arc (150:210:4);
  \nodeSquare{2}{0}
  \nodeSquare{3}{0}
  % label
  \node at (0.8,1.5) {station exterior};
  \node at (3.8,1.5) {station interior};
  \node [align=center] at (1,0) {network\\ arcs};
  \node [align=center] at (4,0) {station\\ arcs};
  \draw [->, shorten >=10pt] (2.85,-1.3) to [bend left=10] (2.5, 0);
  \node [label={[align=center]below:fence node\\valve}] at (2.95,-1.3) {};
  \draw [->, shorten >=10pt] (1.6,-0.95) to [bend right=10] (2, 0);
  \node [label={[align=center]below:fence node}] at (1.5,-0.95) {};
  \draw [->, shorten >=10pt] (3.3,-0.95) to [bend left=10] (3, 0);
  \node [label={[align=center]below:copy of fence node}] at (4.0,-0.95) {};
  \node at (-1,0) {}; % add white space left equivalent to the white space right created by the circle mid point
\end{tikzpicture}
\caption{Schematic fence node valve construction. We create a copy of the fence node and add a new valve between the fence node and its copy, the fence node valve of the node. The artificial station arcs that have been connected to the fence node are instead now connected to the copied node.}
\label{fig:fenceNodeValve}
\end{figure}
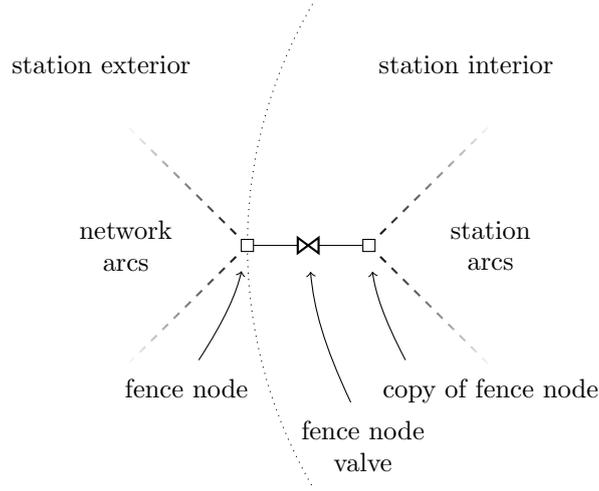

As an extension to the model given in \citep{HopHenLenKoc2021}, we also decouple the pressure at those fence nodes without flow exchange.
To do that, we adopt an idea initially mentioned in \citep{Hop2022}: We create a copy $v'$ of each fence node $v$, create a valve $a'$ in between the two nodes, and reconnect each artificial arc that was connected to the fence node $v$ to its copy $v'$ instead.
The construction is also visualized in Figure~\ref{fig:fenceNodeValve}.
Note that we add the copy $v'$ to \setVertices and the new valve $a$ to \setValves.
Thereby, we create all the corresponding variables and constraints.
We further denote by $a^\mathrm{vaFn}_v$ the \emph{fence node valve} of fence node $v$, by $\setValvesFN\subseteq\setValves$ the set of all fence node valves contained in the network, and by $\setValvesOrig=\setValves\setminus\setValvesFN$ the set of original non-fence-node valves of the network.

In addition to the flow restrictions, some fence nodes $v$ have an upper pressure bound \pressUBExitFG{v} in case flow leaves the station over $v$ for some flow direction, i.e. the fence node serves as an exit of the station.
For all other fence nodes $v$, we set $\pressUBExitFG{v}=\pressUB{v}$.

Finally, not all flow directions fit to all simple states.
Hence, we are given the set $\setSimpleStationValidPairs \subset \setFlowDirections \times \setSimpleStates$ of suitable combinations of flow directions and simple states.

For the model regarding flow directions and fence node valves, we demand for all network stations $i\in\setNetworkStations$ and all times $t\in\setTimestepsNoZero$ that
\begin{align}
    \sum_{f\in\setFlowDirectionsI{i}} \activeFlowDirI{f,t} &= 1 \label{eq:one_flow_dir} \\
    \activeFlowDirI{f,t} &\in\{0,1\} \qquad \forall f\in\setFlowDirectionsI{i}, \notag
\end{align}
for all simple states $s\in\setSimpleStatesI{i}$ of all stations $i\in\setNetworkStations$ and all times $t\in\setTimestepsNoZero$ that
\begin{equation}
    \sum_{(f,s)\in \setSimpleStationValidPairs} \activeFlowDirI{f,t} \geq \activeStateI{s,t}, \label{eq:simple_state_flow_dir_coupling}
\end{equation}
and finally, for all fence nodes $v\in\setFenceNodesNaviStation{i}$ of all stations $i\in\setNetworkStations$ and all times $t\in\setTimestepsNoZero$ that
\begin{subequations}\label{eq:fence_node_model}
\begin{align}
   \sum_{f\in\setFlowDirectionsI{i}: v\in(\flowDirEntriesI{f} \cup \flowDirExitsI{f})} \activeFlowDirI{f,t} &= \modeOp{a^\mathrm{vaFn}_v,t} \\
   \sum_{f\in\setFlowDirectionsI{i}: v\in\flowDirEntriesI{f}} \activeFlowDirI{f,t} = 1
   \quad&\implies\quad \mFlowI{a^\mathrm{vaFn}_v,t} \geq 0\\
   \sum_{f\in\setFlowDirectionsI{i}: v\in\flowDirExitsI{f}} \activeFlowDirI{f,t} = 1
   \quad&\implies\quad \mFlowI{a^\mathrm{vaFn}_v,t} \leq 0\\
   \sum_{f\in\setFlowDirectionsI{i}: v\in\flowDirExitsI{f}} \activeFlowDirI{f,t} = 1, \quad&\implies\quad \pressI{v,t} \leq \pressUBExitFG{v}.
\end{align}
\end{subequations}

\subsection{Demand scenario}\label{sec:futureDemands}
At each boundary node, we are given future requirements in terms of inflow and pressure, which we call the \emph{demands}.
A distinct set of demand values for all boundary nodes of the network is called a \emph{demand scenario} of the problem.
While for each future time step $t\in\setTimesteps$, there is a required inflow value \demandInflow{v,t} for each boundary node $v\in\setBoundaryNodes$, the prescribed pressure \demandPressure{v,t} only exists for entries $v\in\setEntries$ with non-zero inflow values at time $t$.
Furthermore, we do not have to meet them exactly but within a tolerance of $\varepsilon=1\,\text{bar}$ for each value.

However, there is no guarantee that meeting the demands is actually possible.
For these cases, we allow a deviation from the demands, which we measure and penalize in the objective function.
We call this deviation \emph{slack} and capture it in the variables \slackQPos{v,t} and \slackQNeg{v,t} for the positive and negative inflow deviation and \slackPPos{v,t} and \slackPNeg{v,t} for the positive and negative pressure deviation for a boundary node $v\in\setBoundaryNodes$ and time $t\in\setTimestepsNoZero$.
When denoting by \slackQUB the maximum amount of inflow slack per boundary node and time step as well as by \slackPUB the corresponding maximum amount of pressure slack, we can formulate the model of the inflow demands for all boundary nodes $v\in\setBoundaryNodes$ and a future time $t\in\setTimestepsNoZero$ as
\begin{subequations}\label{eq:inflow_demand_model}
\begin{align}
   \inflowI{v,t} &= \demandInflow{v,t} + \slackQPos{v,t} - \slackQNeg{v,t} \\
   0 &\leq \slackQPos{v,t} \leq \slackQUB \\
   0 &\leq \slackQNeg{v,t} \leq \slackQUB,
\end{align}
\end{subequations}
and the model of the pressure demands for all entry nodes $v\in\setEntries$ with non-zero inflow demands \demandInflow{v,t} and a future time $t\in\setTimestepsNoZero$ as
\begin{subequations}\label{eq:pressure_demand_model}
\begin{align}
   \pressI{v,t} &\leq \demandPressure{v,t} + \varepsilon + \slackPPos{v,t} \\
   \pressI{v,t} &\geq \demandPressure{v,t} - \varepsilon - \slackPNeg{v,t}  \\
   0 &\leq \slackPPos{v,t} \leq \slackPUB \\
   0 &\leq \slackPNeg{v,t} \leq \slackPUB.
\end{align}
\end{subequations}

\subsection{Initial state of artificial elements}\label{sec:station_initial_state}
As stated at the end of Section~\ref{sec:gas_flow_in_networks}, we are given an initial state for the network, specifying the pressures, flows, and modes of all nodes and elements of the original network.
However, since the network station model is based on an artificial topology, there are no initial state values for the non-fence nodes inside the station, the artificial arcs, the station's simple state, and its flow direction.

For the inner station nodes, we do not need the initial pressure values as variables in the model.
However, the initial pressures and temperatures are used as parameters in the artificial compressor model.
For these values, we use the corresponding average values per quantity over the set of fence nodes.
As the only exception, we use the initial quantities given for the original fence node for those nodes being fence node copies.
Regarding the flow values, we do not use the initial flow values of artificial arcs.
However, for the objective function, we need the initial flow on the fence node valve (see Section~\ref{sec:objective}), which can be determined by using the flow balance equation \eqref{eq:flow_balance} of the corresponding fence node.

Regarding the discrete states of the network stations as well as their contained artificial arcs, we aim to choose settings that fit the initial state of the fence nodes, fence node valves, and the compressor stations in the original network station topology.
Since there might be multiple suitable settings, we only sort out non-fitting ones and include the final decision into the optimization model.

First, we determine a list of possible initial simple states \setInitActiveSimpleStates.
This is done by checking for each simple state of each network station if there would be a corresponding feasible initial state for the flow direction choice, as well as all the elements in the network station which fits to the initial pressure and flow values of the fence nodes and corresponding fence node valve.
As a second step, we determine the set of initially running compressor units \setInitActiveCompressorUnits based on the active configurations of the compressors stations in the original network station topology.
Finally, we determine the set \setInitActiveAuxCompressorConfigs of those configurations of artificial compressors $c$, for which there are enough initially running compressor units contained in their set of usable compressor units $U_c$, such that the configuration might be active.

Having these sets, we can define the model for the initial time step.
From the above-defined constraints, we use \eqref{eq:one_cfg_per_active_compressor}, \eqref{eq:compressor_conf_conflict}, \eqref{eq:one_simple_state}, and \eqref{eq:auxArc_simple_state_connection_model} also for $t=0$.
In addition, we add the following new constraints:
\begin{align}
\activeStateI{s,0} &= 0 && \forall s\in\setSimpleStates\setminus\setInitActiveSimpleStates \label{eq:init_state_deactivate_simple_states}\\
\activeCfgI{c,0} &= 0 && \forall a\in\setCompressorArcs \quad \forall c\in\setAuxCompressorConfigsI{a}\setminus\setInitActiveAuxCompressorConfigs \label{eq:init_state_deactivate_configs}\\
\sum_{a\in\setCompressorArcs} \quad \sum_{c\in\setAuxCompressorConfigsI{a}: u\in U_c} \activeCfgI{c,0} &\geq 1 && \forall u\in\setInitActiveCompressorUnits. \label{eq:init_state_cover_running_units}
\end{align}
While the first two constraints ensure that only initially valid simple states and artificial configurations are chosen for the initial state, the last constraint demands for each initially running compressor unit $u$ that there is at least one initially active artificial configuration $c$ with $u\in U_c$.

\subsection{Objective function}\label{sec:objective}
The general goal of the model is to reduce the amount of changes necessary for fulfilling the given demands while trying to reduce unnecessary compressor unit usage.
In contrast to the base model of~\citep{HopHenLenKoc2021}, where a multi-level model was formulated to ensure each solution has minimal slack values, we instead used a single-level objective function representing the sum of differently weighted terms.
Another difference from the base model is that we do not perform a solution \emph{smoothing} in a post-processing step but instead include this already in the model, which increases its complexity.

To build our objective function, we need to introduce additional variables to track the corresponding changes of the single entities and connect these with the quantities describing the network state.
We start with the changes of valve modes.
For each non-fence-node valve $a\in\setValvesOrig$ and time $t_1\in\setTimestepsNoZero$, we capture a change in the binary variable $\valveModeChg{a,t_1}$ and specify the corresponding model as
\begin{subequations}\label{eq:orig_valve_change_model}
\begin{align}
    \modeOp{a,t_1} - \modeOp{a,t_0} &\leq \valveModeChg{a,t_1} \\
    \modeOp{a,t_0} - \modeOp{a,t_1} &\leq \valveModeChg{a,t_1} \\
    \valveModeChg{a,t_1} &\in\{0,1\}, \notag
\end{align}
\end{subequations}
where $t_0=t_1-1$ denotes the time step preceding $t_1$.

For the regulator changes, we use the model introduced in \citep{HenAndHopTur2021} tracking changes of the regulator's mode and, in case it is active, the regulator's operation point in terms of changes in the incoming pressure, outgoing pressure, and flow.
For a regulator $a=(\ell,r)\in\setRegulators$ and times $t_1\in\setTimestepsNoZero$ and $t_0=t_1-1$, the model reads as
\begin{subequations}\label{eq:regulator_change_model}
\begin{align}
    \regModeChg{a,t_1} &\geq \modeX{a,t_1} - \modeX{a,t_0} && \forall \text{x}\in\{\text{cl}, \text{op}, \text{ac}\}\\
    \regModeChg{a,t_1} &\leq 2 - \modeX{a,t_1} - \modeX{a,t_0} && \forall \text{x}\in\{\text{cl}, \text{op}, \text{ac}\}\\
    y_{t_1} - y_{t_0} &\leq \rgYChg{a,t_1} + (\modeOp{a,t_1} + \modeCl{a,t_1} + \regModeChg{a,t_1})(\ub{y}_{t_1}-\lb{y}_{t_0}) && \forall y\in\{\pressI{\ell}, \pressI{r}, \mFlowI{a}\}\\
    y_{t_0} - y_{t_1} &\leq \rgYChg{a,t_1} + (\modeOp{a,t_1} + \modeCl{a,t_1} + \regModeChg{a,t_1})(\ub{y}_{t_0}-\lb{y}_{t_1}) && \forall y\in\{\pressI{\ell}, \pressI{r}, \mFlowI{a}\}\\
    \regModeChg{a,t_1} &\in\{0,1\}, \notag
\end{align}
\end{subequations}
where \regModeChg{a,t_1} denotes the binary variable indicating a mode change at the regulator between times $t_0$ and $t_1$, and \rgPLChg{a,t_1}, \rgPRChg{a,t_1}, and \rgQChg{a,t_1} track the changes in the operation point for an active regulator for the incoming pressure, outgoing pressure, and flow through the element.

Next are the changes regarding the network stations and their corresponding elements.
First, we introduce binary variables \auxArcActivityChg{a,t}, representing a change in the activity of the artificial arc $a$ between time $t$ and the previous time step.
The corresponding constraints have the same structure as the valve constraints above and read for an artificial arc $a\in\setArtificialLinks$ and times $t_1\in\setTimestepsNoZero$ and $t_0=t_1-1$ as
\begin{subequations}\label{eq:aux_arc_change_model}
\begin{align}
    \activeArcI{a,t_1} - \activeArcI{a,t_0} &\leq \auxArcActivityChg{a,t_1} \\
    \activeArcI{a,t_0} - \activeArcI{a,t_1} &\leq \auxArcActivityChg{a,t_1} \\
    \auxArcActivityChg{a,t_1} &\in\{0,1\}. \notag
\end{align}
\end{subequations}
For artificial compressors, we also count configuration changes as changes in the arc's activity.
We only need to track if a configuration is newly activated, as turning it off either turns another configuration on or deactivates the whole compressor arc.
Hence, we formulate for each configuration $c\in\setAuxCompressorConfigsI{a}$ of each artificial compressor $a\in\setCompressorArcs$ and times $t_1\in\setTimestepsNoZero$ and $t_0=t_1-1$ the corresponding constraint as
\begin{align}
    \activeCfgI{c,t_1} - \activeCfgI{c,t_0} &\leq \auxArcActivityChg{a,t_1}. \label{eq:aux_comp_conf_change_model}
\end{align}

As in the base model, we track changes in the network stations' simple states but not in the flow directions.
Since each station has exactly one simple state at each point in time, we capture the activation of a simple state $s$ at time $t$ in the binary variable \simpleStateTurnOn{s,t} and give the model for a simple state $s\in\setSimpleStatesI{i}$ of a station $i\in\setNetworkStations$ and times $t_1\in\setTimestepsNoZero$ and $t_0=t_1-1$ as
\begin{align}
    \activeStateI{s,t_1} - \activeStateI{s,t_0} &\leq \simpleStateTurnOn{s,t_1} \label{eq:simple_state_turn_on}\\
    \simpleStateTurnOn{s,t_1} &\in\{0,1\} \notag.
\end{align}

Finally, we capture the pressure changes at each fence node $v$ between time $t$ and the previous time point in the variable \fnPChg{v,t} and the flow changes at the corresponding fence node valve by the variable \fnQChg{v,t}.
The corresponding model for a fence node $v\in\setFenceNodesNaviStation{i}$ of network station $i\in\setNetworkStations$ and times $t_1\in\setTimestepsNoZero$ and $t_0=t_1-1$ is given as
\begin{subequations}\label{eq:fence_node_p_and_q_change}
\begin{align}
    \pressI{v,t_1} - \pressI{v,t_0} &\leq \fnPChg{v,t_1} \\
    \pressI{v,t_0} - \pressI{v,t_1} &\leq \fnPChg{v,t_1} \\
    \mFlowI{a^\mathrm{vaFn}_v,t_1} - \mFlowI{a^\mathrm{vaFn}_v,t_0} &\leq \fnQChg{v,t_1} \\
    \mFlowI{a^\mathrm{vaFn}_v,t_0} - \mFlowI{a^\mathrm{vaFn}_v,t_1} &\leq \fnQChg{v,t_1}.
\end{align}
\end{subequations}

Having defined all the change tracking variables, we are now able to formulate the objective function as a weighted sum of different penalty terms.
Apart from the changes and slack values defined above, we also penalize the usage of compressor units by penalizing the usage of artificial configurations $c$ weighted by their number of used compressor units $u_c$.
For each of the terms in the sum, we have specific weights defining the penalty costs: The demand slack cost per hour \costSlackQ, the pressure slack cost per hour \costSlackP, the cost for each compressor unit running in an active configuration per hour \costUnitRunning, the non-fence-node valve change cost \costVaModeChg, the regulator mode change cost \costRgModeChg, the operation point change costs of active regulators for the incoming pressure \costRgPLChg, the outgoing pressure \costRgPRChg and the flow \costRgQChg, the artificial arc change costs \costArtificialLinkChange, the costs for activating a simple state \costSimpleStateChangeI{s}, which are specific values manually defined by the industry experts for each simple state $s\in\setSimpleStates$, and the costs for changes of the fence node pressure \costFenceNodePChg and flow changes on the corresponding fence node valves \costFenceNodeQChg.
The complete objective is then given as
\begin{align}
    \min \sum_{t\in\setTimestepsNoZero} \Bigg( \quad & \sum_{v\in\setBoundaryNodes} \frac{\granularityI{t}-\granularityI{t-1}}{1\,\text{h}} \left( \costSlackQ (\slackQPos{v,t} + \slackQNeg{v,t}) + \costSlackP (\slackPPos{v,t} + \slackPNeg{v,t}) \right) \notag \\
    +& \sum_{a\in\setCompressorArcs}\sum_{c\in\setAuxCompressorConfigsI{a}} \frac{\granularityI{t}-\granularityI{t-1}}{1\,\text{h}} \cdot u_c \cdot \costUnitRunning \cdot \activeCfgI{c,t} \notag \\
    +& \sum_{a\in\setValvesOrig} \costVaModeChg \cdot \valveModeChg{a,t} \notag \\
    +&  \sum_{a\in\setArtificialLinks} \costArtificialLinkChange \cdot \auxArcActivityChg{a,t} \label{eq:objective} \\
    +&  \sum_{i\in\setNetworkStations} \sum_{s\in\setSimpleStatesI{i}} \costSimpleStateChangeI{s} \cdot \simpleStateTurnOn{s,t} \notag \\
    +& \sum_{a\in\setRegulators} \left(\costRgModeChg \cdot \regModeChg{a,t} + \costRgPLChg \cdot \rgPLChg{a,t} + \costRgPRChg \cdot \rgPRChg{a,t} + \costRgQChg \cdot \rgQChg{a,t} \right) \notag \\
    +& \sum_{i\in\setNetworkStations} \sum_{v\in\setFenceNodesNaviStation{i}} \left( \costFenceNodePChg \cdot \fnPChg{v,t} + \costFenceNodeQChg \cdot \fnQChg{v,t} \right)
    \Bigg). \notag
\end{align}

\subsection{Complete model}
By combining all the above-defined variables and constraints with the objective function, we can formulate our gas flow problem \MINLP as the following MINLP model:
\begin{align*}
  \min\quad & \eqref{eq:objective} \\
  \text{s.t.} \quad \forall t\in\setTimestepsNoZero \qquad \qquad \\
  & \eqref{eq:flow_balance} \quad \forall v\in\setVertices
  && \eqref{eq:pipe_model} \quad \forall a\in\setPipes \\
  & \eqref{eq:valve_model} \quad \forall a\in\setValves
  && \eqref{eq:regulator_model},\eqref{eq:regulator_change_model} \quad \forall a\in\setRegulators \\
  & \eqref{eq:artificial_shortcut_model} \quad \forall a\in\setShortcuts
  && \eqref{eq:artificial_regulator_model} \quad \forall a\in\setRegulatorArcs \\
  & \eqref{eq:artificial_biregulator_model} \quad \forall a\in\setBiRegulatorArcs
  && \eqref{eq:artificial_compressor_model} \quad \forall a\in\setCompressorArcs \\
  & \eqref{eq:artificial_compressor_config_model},\eqref{eq:aux_comp_conf_change_model} \quad \forall a\in\setCompressorArcs \quad \forall c\in\setAuxCompressorConfigsI{a}
  && \eqref{eq:compressor_conf_conflict} \quad \forall z\in\setAuxCompressorConfigConflicts \\
  & \eqref{eq:one_simple_state},\eqref{eq:one_flow_dir} \quad \forall i\in\setNetworkStations
  && \eqref{eq:auxArc_simple_state_connection_model},\eqref{eq:aux_arc_change_model} \quad \forall i\in\setNetworkStations \quad \forall a\in\setArtificialLinksNaviStation{i}\\
  & \eqref{eq:simple_state_flow_dir_coupling},\eqref{eq:simple_state_turn_on} \quad \forall i\in\setNetworkStations \quad \forall s\in\setSimpleStatesI{i}
  && \eqref{eq:fence_node_model},\eqref{eq:fence_node_p_and_q_change} \quad \forall i\in\setNetworkStations \quad \forall v\in\setFenceNodesNaviStation{i}\\
  & \eqref{eq:inflow_demand_model} \quad \forall v\in\setBoundaryNodes
  && \eqref{eq:pressure_demand_model} \quad \forall v\in\setEntries: \demandInflow{v,t} \neq 0 \\
  & \eqref{eq:orig_valve_change_model} \quad \forall a\in\setValvesOrig \\
  % initial state
  \quad t=0 \qquad \qquad \\
  & \eqref{eq:one_cfg_per_active_compressor} \quad \forall a\in\setCompressorArcs
  && \eqref{eq:compressor_conf_conflict} \quad \forall z\in\setAuxCompressorConfigConflicts \\
  % start problem name
  \MINLP\rule{4cm}{0.0cm}
  % end problem name
  & \eqref{eq:one_simple_state} \quad \forall i\in\setNetworkStations
  && \eqref{eq:auxArc_simple_state_connection_model} \quad \forall i\in\setNetworkStations \quad \forall a\in\setArtificialLinksNaviStation{i} \\
  & \eqref{eq:init_state_deactivate_simple_states} \quad \forall s\in\setSimpleStates\setminus\setInitActiveSimpleStates
  && \eqref{eq:init_state_deactivate_configs} \quad \forall a\in\setCompressorArcs \quad \forall c\in\setAuxCompressorConfigsI{a}\setminus\setInitActiveAuxCompressorConfigs \\
  & \eqref{eq:init_state_cover_running_units} \quad \forall u\in\setInitActiveCompressorUnits \\
  % binaries
  \forall t\in\setTimesteps \qquad \qquad & \\
  \{0,1\}\ni \enskip
  & \modeOp{a,t} \quad \forall a\in\setValves,
  && \modeOp{a,t},\modeAc{a,t},\modeCl{a,t} \quad \forall a\in\setRegulators, \\
  & \activeArcI{a,t} \quad \forall a\in\setArtificialLinks,
  && \activeCfgI{c,t} \quad \forall a\in\setCompressorArcs \quad \forall c\in\setAuxCompressorConfigsI{a}, \\
  & \activeStateI{s,t} \quad \forall i\in\setNetworkStations \quad \forall s\in\setSimpleStatesI{i} \\
  % binaries for future time steps
  \forall t\in\setTimestepsNoZero \qquad \qquad & \\
  \{0,1\}\ni \enskip
  & \activeFwdI{a,t},\activeBwdI{a,t} \quad \forall a\in\setBiRegulatorArcs,
  && \activeFlowDirI{f,t} \quad \forall i\in\setNetworkStations \quad \forall f\in\setFlowDirectionsI{i}, \\
  & \valveModeChg{a,t} \quad \forall a\in\setValvesOrig,
  && \regModeChg{a,t} \quad \forall a\in\setRegulators, \\
  & \auxArcActivityChg{a,t} \quad \forall a\in\setArtificialLinks,
  && \simpleStateTurnOn{s,t} \quad \forall s\in\setSimpleStates.
\end{align*}
\section{MIP-based heuristic algorithms}\label{sec:mip_based_heuristics}
We aim at solving the transient gas flow problem presented in the previous section for challenging demand scenarios on large real-world-based instances.
Since this is currently out of reach for black-box MINLP solvers, we create a specialized solution approach given as Algorithm~\ref{algo:baseAlgo}.
It is based on the idea of dealing with the two sources of complexity, the non-continuity of binary variables and the non-convexity of the constraints, separately:
In the first step, we determine a promising and feasible set of solution values for all the binary variables, which we call a \emph{binary assignment} of the problem.
Afterwards, we transform the problem into an NLP by fixing these values and then solve this remaining problem with a black-box NLP solver.
Finally, we return the overall best solution out of all successful NLP runs.
If non of the runs found a feasible solution, the heuristic failed.
Note that the algorithm could be easily parallelized by solving the NLP model immediately after finding a new binary assignment and in parallel to the ongoing search for binary assignments.
Hence, solutions would be obtainable early during the execution, which is especially valuable in environments with strict time limit restrictions.

\begin{algorithm}[tbp]
  \SetKwInOut{Parameters}{Parameters}
  \DontPrintSemicolon
  \KwData{An instance \inst of problem \MINLP}
  \KwResult{A valid solution for \inst}
  % vars
  \SetKwData{baSet}{baSet}
  \SetKwData{solSet}{solSet}
  % functions
  \SetKwFunction{findBinaryAssignments}{findBinaryAssignments}
  \SetKwFunction{solveNLP}{solveNLP}
  \SetKwFunction{chooseBestSol}{chooseBestSol}
  % algo start
  \baSet $ \gets $ \findBinaryAssignments{\inst}\;
  \solSet $ \gets $ \solveNLP{\baSet, \inst}\;
  \KwRet{\chooseBestSol{\solSet}}\;
  \caption{Base algorithm}
  \label{algo:baseAlgo}
\end{algorithm}

Since we consider especially challenging and large instances, we purely focus on finding good feasible solutions fast, i.e., we do not aim for valid global lower bounds for the problem and thereby prove optimality of our solutions.
However, note that in our approach, the black-box NLP solver usually provides, at least locally, optimal solutions for the corresponding sub-problem for a given binary assignment.

In the remainder of this section, we focus on finding binary assignments that yield high-quality solutions for the overall problem.
Note that we sometimes claim that these binary-assignment-producing algorithms solve the overall problem, which we use as a short form of explaining that we solve the problem by using them as function \texttt{findBinaryAssignments} in Algorithm~\ref{algo:baseAlgo}.

\subsection{Sequential mixed-integer programming}
The core of our algorithm for finding binary assignments consists of a sequential algorithm solving linearized versions of the original MINLP problem, which we refer to as \texttt{SMIP} and which is given as Algorithm~\ref{algo:standardSMIP}.
After initializing a linearization for the non-convex Momentum equation~\eqref{eq:fd_discrete_momentum} in Line~\ref{line:SMIP_linearizationInit}, we solve the remaining MIP problem.
If successful, this yields a valid binary assignment in the sense that it respects all the constraints involving no continuous variables.
We then check if the binary assignment was already found before and if the linearization error is small enough to fulfill the convergence criteria.
If both are not the case and we have not reached the maximum number of iterations \paramMaxIt, we update our linearization based on the current solution and repeat the process by solving the updated MIP model.
The algorithm finally returns the set of found binary assignments to be completed to full solutions for the original problem \MINLP.
Note that given a valid binary assignment, the resulting NLP model might be infeasible and hence does not yield a feasible overall solution.

\begin{algorithm}[ht]
  \SetKwInOut{Parameters}{Parameters}
  \DontPrintSemicolon
  \KwData{An instance \inst of problem \MINLP}
  \Parameters{%
    Maximum number of iteration \paramMaxIt\\
    Absolute convergence tolerance \paramMaxAbsTol\\
    Relative convergence tolerance \paramMaxRelTol
  }
  \KwResult{A list of valid binary assignments for \inst}
  % vars
  \SetKwData{baSet}{baSet}
  \SetKwData{mip}{mip}
  \SetKwData{foundPrimalSolution}{foundPrimalSolution}
  \SetKwData{solutionPointer}{solutionPointer}
  \SetKwData{ba}{ba}
  % functions
  \SetKwFunction{createSet}{createSet}
  \SetKwFunction{applyInitialLinearization}{applyInitialLinearization}
  \SetKwFunction{solve}{solve}
  \SetKwFunction{getBinarySolutionValues}{getBinarySolutionValues}
  \SetKwFunction{add}{add}
  \SetKwFunction{isLinearizationErrorSmall}{isLinearizationErrorSmall}
  \SetKwFunction{linearizeBasedOnPreviousSolution}{linearizeBasedOnPreviousSolution}
  % algo start
  \baSet $\gets$ \createSet{}\;
  \mip $\gets$ \applyInitialLinearization{\inst} \label{line:SMIP_linearizationInit}\;
  \For{$i \leftarrow 1$ \KwTo \paramMaxIt}{
    \tcp{solve the linearized MIP problem}
    \foundPrimalSolution, \solutionPointer $ \gets $\solve{\mip} \label{line:SMIP_mipSolve}\;
    \eIf{\foundPrimalSolution}{
      \tcp{Found a feasible solution for the linearized MIP}
      \ba $\gets$ \getBinarySolutionValues{\solutionPointer}\;
      \If{\ba $\in$ \baSet}{
        \tcp{Found an already known binary assignment\\$\Rightarrow$ abort the iteration}
        \break \tcp*[h]{the loop}\;
      }
      \tcp{Found a new binary assignment}
      \baSet.\add{\ba}\;
      \If{\isLinearizationErrorSmall{\mip, \solutionPointer, \paramMaxAbsTol, \paramMaxRelTol}\label{line:SMIP_checkLinearizationError}}{
        \tcp{Given solution fulfills the convergence criterion\\$\Rightarrow$ abort the iteration}
        \break \tcp*[h]{the loop}\;
      }
      \tcp{update linearization}
      \mip $\gets$ \linearizeBasedOnPreviousSolution{\inst,\solutionPointer} \label{line:SMIP_linearizationUpdate}\;
    }
    {
      \tcp{linearized MIP no successful\\$\Rightarrow$ abort the iteration}
      \break \tcp*[h]{the loop}\;
    }
  }
  \KwRet{\baSet}\;
  \caption{General sequential mixed-integer programming algorithm (\texttt{SMIP})}
  \label{algo:standardSMIP}
\end{algorithm}

\subsubsection{Linearization and convergence}
As already mentioned above, the only non-linear term in the model is the friction term in the Momentum equation~\eqref{eq:fd_discrete_momentum} for pipes given as
\begin{equation*}
    c_a \cdot \left(\frac{|\mFlowI{\ell,a,t}| \mFlowI{\ell,a,t}}{\pressI{\ell,t}} + \frac{|\mFlowI{r,a,t}| \mFlowI{r,a,t}}{\pressI{r,t}}\right) \qquad \text{ with the constant } c_a := \frac{\fricI{a} \sGasConst \tempI{a} \zFactI{a} \lenI{a}}{4 \areaI{a}^2 \diamI{a}}.
\end{equation*}
Hence, we are given for each Momentum equation of pipe $a=(\ell,r)\in\setPipes$ and time $t\in\setTimestepsNoZero$ two non-linear functions of the form $c_a (|\mFlowI{x,a,t}| \mFlowI{x,a,t}/\pressI{x,t})$, one for each end node $x\in\{\ell,r\}$.
Both functions are approximated by a linear combination of the involved quantities as
\begin{equation*}
    c_a \frac{|\mFlowI{x,a,t}| \mFlowI{x,a,t}}{\pressI{x,t}} \approx \alpha^0_{x,a,t} + \alpha^p_{x,a,t} \cdot \pressI{x,t} + \alpha^q_{x,a,t} \cdot \mFlowI{x,a,t}.
\end{equation*}

To check if a solution of a linearized MIP model has converged against a non-linear solution, we use two tolerance parameters: The absolute convergence tolerance \paramMaxAbsTol and the relative convergence tolerance \paramMaxRelTol.
Given these, we define that a solution with pressure values $\press'$ and mass flow values $\mFlow'$ converged, if for all pipes and future time points the differences between the linearized friction function and the non-linear friction function with respect to $\press'$ and $\mFlow'$ are smaller than the tolerances.
Hence, we check in function \texttt{isLinearizationErrorSmall} used in Line~\ref{line:SMIP_checkLinearizationError} of Algorithm~\ref{algo:standardSMIP} if the following constraints are fulfilled:
\begin{align*}
    \forall t\in\setTimestepsNoZero \quad \forall a=(\ell,r)\in\setPipes: & \\
    && \Lambda^{\text{lin}} =& \quad \alpha^0_{\ell,a,t} + \alpha^p_{\ell,a,t} \cdot \pressPrevI{\ell,t} + \alpha^q_{\ell,a,t} \cdot \mFlowPrevI{\ell,a,t} \\
    && &+ \alpha^0_{r,a,t} + \alpha^p_{r,a,t} \cdot \pressPrevI{r,t} + \alpha^q_{r,a,t} \cdot \mFlowPrevI{r,a,t} \\
    &&\Lambda^{\text{orig}} =& \quad c_a \cdot \left(\frac{|\mFlowPrevI{\ell,a,t}| \mFlowPrevI{\ell,a,t}}{\pressPrevI{\ell,t}} + \frac{|\mFlowPrevI{r,a,t}| \mFlowPrevI{r,a,t}}{\pressPrevI{r,t}}\right) \\
    && \paramMaxAbsTol \geq& \quad |\Lambda^{\text{lin}} - \Lambda^{\text{orig}}| \\
    \text{ and, if }|\Lambda^{\text{orig}}| > \varepsilon^\text{zero}, & \\
    && \paramMaxRelTol \geq& \quad \frac{|\Lambda^{\text{lin}} - \Lambda^{\text{orig}}|}{|\Lambda^{\text{orig}}|}.
\end{align*}
Here, the value of the linearized friction function with respect to the variable values of the given solution is denoted by $\Lambda^{\text{lin}}$ and the corresponding non-linear friction function value by $\Lambda^{\text{orig}}$.
Furthermore, we can only check for the relative tolerance if the non-linear friction value is not zero.
We ensure this by comparing the term against a corresponding epsilon value $\varepsilon^\text{zero}$, which we define to be equal to $0.001\,\text{bar}$.

For the initial linearization determined in function \texttt{applyInitialLinearization} of Algorithm~\ref{algo:standardSMIP} in Line~\ref{line:SMIP_linearizationInit}, we use the constant velocity approximation introduced by \citep{Hen2018}, which was also used as initial linearization in \citep{HopHenLenKoc2021} and is given as
\begin{align*}
    \alpha^0_{x,a,t} &= 0 \quad \alpha^p_{x,a,t} = 0 \quad \alpha^q_{x,a,t} = \frac{\fricI{a} \lenI{a}}{4 \areaI{a} \diamI{a}} \absVelo{x,a,0} \\
    \text{ with } \quad \absVelo{x,a,0} &:= \frac{\sGasConst \tempI{a} \zFactI{a}}{\areaI{a}}\frac{q_0}{\pressI{x,0}}, \quad\text{ and }\quad q_0 := \max\{|\mFlowI{x,a,0}|, q^\text{min}\}.
\end{align*}
Here, \absVelo{x,a,0} denotes the absolute value of the gas velocity based on the initial state quantities using a minimal absolute flow value of $q^\text{min}$ to avoid $\alpha^q_{x,a,t}$ being zero.
The advantage of this linearization is that it is symmetric in $q$ and also zero for $q=0$.

For all subsequent iterations of Algorithm~\ref{algo:standardSMIP}, we determine the linearization based on the pressure values $\press'$ and flow values $\mFlow'$ of a previous MIP solution by function \texttt{linearizeBasedOnPreviousSolution} in Line~\ref{line:SMIP_linearizationUpdate}.
As linearization, we use the tangential plane on the non-linear function in the point ($\press'$, $\mFlow'$), which is the first-order Taylor expansion at this point and is given as
\begin{align*}
    \alpha^p_{x,a,t} &= c_a \cdot \frac{-|\mFlowPrevI{x,a,t}| \mFlowPrevI{x,a,t}}{(\pressPrevI{x,t})^2} \\
    \alpha^q_{x,a,t} &= c_a \cdot \frac{2|\mFlowPrevI{x,a,t}|}{\pressPrevI{x,t}} \\
    \alpha^0_{x,a,t} &= c_a \cdot \frac{|\mFlowPrevI{x,a,t}| \mFlowPrevI{x,a,t}}{\pressPrevI{x,t}} - \alpha^p_{x,a,t}\cdot\pressPrevI{x,t} - \alpha^q_{x,a,t}\cdot\mFlowPrevI{x,a,t}.
\end{align*}

\subsubsection{Comments on infeasibility}
Algorithm~\ref{algo:standardSMIP} struggles with infeasibility in two different ways:
First, the linearized MIP problems of each iteration might be infeasible even if the overall MINLP problem \MINLP is not.
Second, the algorithm has no functionality to detect the overall infeasibility of the general problem \MINLP.

Despite those structural weaknesses, we still expect the algorithm to perform well on real-world-based instances, like those used in our computational experiments in Section~\ref{sec:computational_experiments}.
On these, infeasibility usually does not occur, as real-world gas transport networks are, in general, very flexible and allow for a broad range of different pressure, flow, and inflow values.
In addition, the slack values can adjust problematic inflow and pressure demands.

The computational results we present in Section~\ref{sec:computational_experiments} confirm our expectations, as we found a feasible solution for each instance of the problem \MINLP and rarely encountered infeasibility when solving the linearized MIP models, see Table~\ref{tab:heuristicResultStats}.

\subsection{Reduced time horizon heuristics}\label{sec:reduceHorizonHeuristics}
To further accelerate the solving process,
we introduce additional heuristics that find binary assignments by solving a series of models on time horizons of reduced size, where the size of a time horizon is defined as the number of future time steps.
The heuristics are based on two well-known ideas for time-expanded problems: The rolling horizon and the aggregated horizon.
While rolling horizon heuristics rely on iteratively solving the model for a subset of time steps that are gradually shifted to the future, the aggregated horizon combines certain time steps and afterwards completes the aggregated solution to one for the whole time horizon.

Note that we maintain the general sequential algorithm structure of Algorithm~\ref{algo:standardSMIP} for all of the heuristics below.

\subsubsection{Rolling horizon heuristic}
\begin{algorithm}[htp]
  \SetKwInOut{Parameters}{Parameters}
  \DontPrintSemicolon
  \KwData{An instance \inst of problem \MINLP}
  \Parameters{%
    Size of the rolling horizon \paramRolHorizon\\
    Maximum number of iteration \paramMaxIt\\
    Absolute convergence tolerance \paramMaxAbsTol\\
    Relative convergence tolerance \paramMaxRelTol
  }
  \KwResult{A list of valid binary assignments for \inst}
  % vars
  \SetKwData{mip}{mip}
  \SetKwData{initialTimeStep}{initialTimeStep}
  \SetKwData{savedSolutionStates}{savedSolutionStates}
  \SetKwData{startState}{startState}
  \SetKwData{lastTimeStep}{lastTimeStep}
  \SetKwData{reducedMip}{reducedMip}
  \SetKwData{foundPrimalSolution}{foundPrimalSolution}
  \SetKwData{solutionPointer}{solutionPointer}
  % functions
  \SetKwFunction{solveByRollingHorizonHeuristic}{solveByRollingHorizonHeuristic}
  \SetKwFunction{createList}{createList}
  \SetKwFunction{append}{append}
  \SetKwFunction{getStartState}{getStartState}
  \SetKwFunction{getLast}{getLast}
  \SetKwFunction{createModelForReducedTimeHorizon}{createModelForReducedTimeHorizon}
  \SetKwFunction{solve}{solve}
  \SetKwFunction{getFirstFutureState}{getFirstFutureState}
  \SetKwFunction{appendAll}{appendAll}
  \SetKwFunction{getAllFutureState}{getAllFutureState}
  \SetKwFunction{combineStatesToOverallSolution}{combineStatesToOverallSolution}
  \SetKwFunction{SMIP}{SMIP}
  % function declaration
  \SetKwProg{MajorFunc}{Function}{ is}{}
  % algo start
  \MajorFunc{\solveByRollingHorizonHeuristic{\mip}}{
    \initialTimeStep $\gets 0$\;
    \savedSolutionStates $\gets$ \createList{}\;
    \savedSolutionStates.\append{\getStartState{\inst}}\;
    \While{$\initialTimeStep+\paramRolHorizon \leq \nTimesteps$}{
      \startState  $\gets$ \savedSolutionStates.\getLast{}\;
      \lastTimeStep $\gets \initialTimeStep+\paramRolHorizon$\;
      \tcp{Create MIP model with reduced time horizon}
      \reducedMip $\gets$ \createModelForReducedTimeHorizon{\inst, \startState, \initialTimeStep, \lastTimeStep} \label{line:RH_createReducedModel}\;
      \foundPrimalSolution, \solutionPointer $\gets$ \solve{\reducedMip}\;
      \eIf{\foundPrimalSolution}{
        \tcp{Found feasible solution for the reduce time horizon mip\\check if this is the last iteration}
        \eIf{$\initialTimeStep+\paramRolHorizon < \nTimesteps$}{
          \tcp{This is not the last iteration\\$\Rightarrow$ save state of first future time step}
          \savedSolutionStates.\append{\getFirstFutureState{\solutionPointer}} \label{line:RH_saveFirstFutureTimeStep}\;
        }{
          \tcp{This is the last iteration!\\$\Rightarrow$ Save all future solution states}
          \savedSolutionStates.\appendAll{\getAllFutureState{\solutionPointer}} \label{line:RH_saveAllFutureTimeStep}\;
        }
        \tcp{Prepare the next iteration}
        $\initialTimeStep \gets \initialTimeStep +1$\;
      }{
        \tcp{Reduced time horizon mip was not successful\\return foundPrimalSolution=False}
        \KwRet{False, nullptr} \label{line:RH_failure}\;
      }
    }
    \tcp{All iterations have been successful!\\$\Rightarrow$ Return final solution}
    \KwRet{True, \combineStatesToOverallSolution{\savedSolutionStates}} \label{line:RH_combineSolutionStates}\;
  }
  \tcp{Use the \texttt{SMIP} algorithm, but replace the `solve' function}
  \KwRet{\SMIP{\inst, \paramMaxIt, \paramMaxAbsTol, \paramMaxRelTol, $\solve=\solveByRollingHorizonHeuristic$}}\;
  \caption{Sequential rolling horizon heuristic algorithm (\texttt{RH})}
  \label{algo:rollingHorizon}
\end{algorithm}
In our rolling horizon heuristic \texttt{RH} given in Algorithm~\ref{algo:rollingHorizon}, we only change a single line in contrast to the \texttt{SMIP} Algorithm~\ref{algo:standardSMIP}:
Instead of solving the linearized MIP for the whole time horizon in Line~\ref{line:SMIP_mipSolve}, we find a solution by using a rolling horizon approach of fixed length \paramRolHorizon.
We start by building a reduced model in Line~\ref{line:RH_createReducedModel} that only considers the first \paramRolHorizon future time steps by not adding any variables or constraints for the remaining ones.
After solving the reduced model, we save all variable values regarding the first of the considered future time steps for our overall solution in Line~\ref{line:RH_saveFirstFutureTimeStep}.
We call this set of values the \emph{solution state} of this time step and use it as the initially fixed state when building the reduced model in the next iteration.
Note that we still use the overall initial time step 0 to determine the model's constant parameters, like, for example, the fixed temperature $T_a$ of a pipe $a\in\setPipes$.
Finally, for the next iteration, we also include the earliest not already included future time step into our reduced time horizon, which thereby again has a length of \paramRolHorizon.
We repeat this procedure until we reach the iteration containing the very last future time step \nTimesteps.
After saving all the solution states of this iteration in Line~\ref{line:RH_saveAllFutureTimeStep}, we can create a solution for the entire time horizon by combining all the collected solution states in Line~\ref{line:RH_combineSolutionStates}.

As a minor detail, we note that for all variables without a fixed initial state for time $t=0$, for example, for the simple state variables \activeState, we use the model as defined in Section~\ref{sec:model} for the first iteration and then add the there determined values for $t=0$ to the corresponding solution state.
For all subsequent iterations, we fix the values determined in the previous iteration for the initial state and can therefore ignore all the constraints added to the model for $t=0$.

\subsubsection{Aggregated horizon heuristic}
For the aggregated horizon heuristic \texttt{AH} stated as Algorithm~\ref{algo:aggregatedHorizon}, we create a time horizon of reduced size \paramAggHorizon by combining certain time steps, resulting in longer but fewer time steps in the overall horizon.
We represent this by a function $\vartheta: \{0,1,\dots,\paramAggHorizon\} \rightarrow \setTimesteps=\{0,1,\dots,\nTimesteps\}$, which maps the time steps of the aggregated time horizon to time steps of the original time horizon and is created as variable \texttt{aggToOrigMap} by calling the function \texttt{mapAggregatedToOriginalTimeSteps} in Line~\ref{line:AH_createAggregatedTimeSteps}.
A time interval between two subsequent time steps $t$ and $t+1$ in the aggregated horizon represents the combination of all original time intervals between the time points $\vartheta(t)$ and $\vartheta(t+1)$.
To combine a similar amount of original time steps for each time interval, we determine $\vartheta$ as
\begin{align*}
  \bar{c} &:= \lceil \frac{\nTimesteps}{\paramAggHorizon} \rceil  \qquad \qquad \qquad \text{the larger number of original time intervals to combine} \\
  \ubar{c} &:= \lfloor \frac{\nTimesteps}{\paramAggHorizon} \rfloor  \qquad \qquad \qquad \text{the smaller number of original time intervals to combine} \\
  \bar{n} &:= \nTimesteps \mod{} \paramAggHorizon \qquad \text{the amount of larger aggregated time intervals}\\
  \vartheta(t) &:=
  \begin{cases}
      t \cdot \bar{c} & \text{ if } t \leq \bar{n} \\
      \bar{c}\cdot\bar{n} + t \cdot \ubar{c} & \text{ if } \bar{n} < t
  \end{cases},
\end{align*}
and give an example as
\begin{equation*}
  \nTimesteps = 7 \qquad\&\qquad \paramAggHorizon = 4 \qquad\implies\qquad
  \begin{array}{l|rrrrr}
    t            & 0 & 1 & 2 & 3 & 4\\
    \vartheta(t) & 0 & 2 & 4 & 6 & 7
  \end{array}.
\end{equation*}
Note that $\vartheta(0)=0$ and $\vartheta(\paramAggHorizon)=\nTimesteps$ hold.
If \paramAggHorizon is a divisor of \nTimesteps, each aggregated time interval combines $\bar{c}=\ubar{c}$ original time intervals.

\begin{algorithm}[tbh]
  \SetKwInOut{Parameters}{Parameters}
  \DontPrintSemicolon
  \KwData{An instance \inst of problem \MINLP}
  \Parameters{%
    Size of the aggregated horizon \paramAggHorizon\\
    Maximum number of iteration for the aggregated horizon MIP \paramMaxItAgg\\
    Maximum number of iteration for the full horizon MIP \paramMaxItFull\\
    Absolute convergence tolerance \paramMaxAbsTol\\
    Relative convergence tolerance \paramMaxRelTol
  }
  \KwResult{A list of valid binary assignments for \inst}
  % vars
  \SetKwData{aggToOrigMap}{aggToOrigMap}
  \SetKwData{aggregatedBAs}{aggregatedBAs}
  \SetKwData{varTypesToFix}{varTypesToFix}
  \SetKwData{baSetFullHorizon}{baSetFullHorizon}
  \SetKwData{aggBA}{aggBA}
  \SetKwData{baSetForAggBA}{baSetForAggBA}
  \SetKwData{mip}{mip}
  % functions
  \SetKwFunction{mapAggregatedToOriginalTimeSteps}{mapAggregatedToOriginalTimeSteps}
  \SetKwFunction{createTimeAggregatedProblemInstance}{createTimeAggregatedProblemInstance}
  \SetKwFunction{SMIP}{SMIP}
  \SetKwFunction{removeDuplicateBAsWithRespectToVarsToFix}{removeDuplicateBAsWithRespectToVarsToFix}
  \SetKwFunction{createList}{createList}
  \SetKwFunction{createProblemInstanceWithPartiallyFixedBinaries}{createProblemInstanceWithPartiallyFixedBinaries}
  \SetKwFunction{applyInitialLinearization}{applyInitialLinearization}
  \SetKwFunction{initializeLinearizationBasedOnAggSolution}{initializeLinearizationBasedOnAggSolution}
  \SetKwFunction{solve}{solve}
  \SetKwFunction{solveMIPUsingPartialStartSolution}{solveMIPUsingPartialStartSolution}
  \SetKwFunction{addAll}{addAll}
  % algo start
  \aggToOrigMap $\gets$ \mapAggregatedToOriginalTimeSteps{\inst,\paramAggHorizon} \label{line:AH_createAggregatedTimeSteps}\;
  \tcp{Create problem on aggregated time horizon}
  $\instAgg \gets \createTimeAggregatedProblemInstance{\inst, \aggToOrigMap}$ \label{line:AH_createAggregatedProblemInstance}\;
  \tcp{Solve aggregated problem via \texttt{SMIP}}
  \aggregatedBAs $\gets$ \SMIP{\instAgg, \paramMaxItAgg, \paramMaxAbsTol, \paramMaxRelTol} \label{line:AH_solveAggregatedProblem}\;
  \varTypesToFix $\gets \{$ \\
  \nonl\enskip$\modeOp{a},\modeAc{a},\modeCl{a} \enskip \forall a\in\setRegulatorArcs, \quad \modeOp{a} \enskip \forall a\in\setValvesOrig,$\\
  \nonl\enskip $\activeArcI{a} \enskip \forall a\in\setArtificialLinks, \quad \activeCfgI{c} \enskip \forall c\in\setAuxCompressorConfigsI{a} \enskip \forall a\in\setCompressorArcs, \quad \activeStateI{s} \enskip \forall s\in\setSimpleStates\}$ \label{line:AH_defineVariablesToFix}\;
  \aggregatedBAs $\gets$ \removeDuplicateBAsWithRespectToVarsToFix{\aggregatedBAs, \varTypesToFix} \label{line:AH_removeDuplicatedBasBasedOnFixedVars}\;
  \baSetFullHorizon $\gets$ \createList{}\;
  \tcp{Complete the time aggregated binary assignments to those for the full horizon}
  \For{$\aggBA \in \aggregatedBAs$}{
    \instFixed $\gets$ \createProblemInstanceWithPartiallyFixedBinaries{\inst, \aggBA, \varTypesToFix, \aggToOrigMap} \label{line:AH_fixAggregatedSolutionInFullHorizonInstance}\;
    \tcp{Solve \texttt{SMIP} with adjusted initial linearization and partial start solution}
    \baSetForAggBA $\gets$ \SMIP{\instFixed, \paramMaxItFull, \paramMaxAbsTol, \paramMaxRelTol,\\
    \nonl\enskip$\applyInitialLinearization=$\\
    \nonl\enskip\enskip\initializeLinearizationBasedOnAggSolution{\instFixed, \aggBA},\\
    \nonl\enskip$\solve=\solveMIPUsingPartialStartSolution{\mip, \aggBA}$}
    \label{line:AH_sMIPCallWithFixedVars}\;
    \baSetFullHorizon.\addAll{\baSetForAggBA}\;
  }
  \tcp{Return all full horizon binary assignment found be completing aggregated ones}
  \KwRet{\baSetFullHorizon}\;
  \caption{Sequential aggregated horizon heuristic algorithm (\texttt{AH})}
  \label{algo:aggregatedHorizon}
\end{algorithm}

Using the function $\vartheta$, we can now derive a time aggregated problem instance \instAgg from the original problem instance \inst by calling the function \texttt{createTimeAggregatedProblemInstance} in Line~\ref{line:AH_createAggregatedProblemInstance} of Algorithm~\ref{algo:aggregatedHorizon}.
Since all the elements sets of \inst and \instAgg coincide except of the time horizon, we are only missing a description of the time-dependent instance parameters.
In a first step, we determine the inflow demands.
As they represent the average inflow rate over a time interval, we determine the time aggregated inflow demands as the average of original inflow demands, weighted by the length of the corresponding time intervals, i.e., for the aggregated time $t\in\{1,\dots,\paramAggHorizon\}$ and a boundary node $v\in\setBoundaryNodes$ we define
\begin{equation*}
  \demandInflow{v,t} := \frac{\sum_{i=\vartheta(t-1)+1}^{\vartheta(t)} \left(\granularityI{i}-\granularityI{i-1}\right)\cdot\demandInflow{v,i} }{\granularityI{\vartheta(t)}-\granularityI{\vartheta(t-1)}}.
\end{equation*}
For all the remaining time-dependent parameters, like, for example, the function \granularityI{t} representing the time difference between $t$ and the initial time step in seconds, the pressure demands at the entries, or the general pressure and flow bounds, we use for each aggregated time point $t$ the parameter value of the corresponding original time point $\vartheta(t)$ in \inst and thereby complete the aggregated instance \instAgg.

Using the \texttt{SMIP} Algorithm~\ref{algo:standardSMIP}, we determine a set of valid binary assignments for the aggregated problem instance.
Note that we use the parameter \paramMaxItAgg to specify the maximum number of iterations to use in \texttt{SMIP} to solve \instAgg.
In the following, we explain the steps to complete these to full binary assignments for the original problem instance \inst.
The general idea is to fix certain variables in the model based on the aggregated binary assignment and thereby create the problem instance \instFixed, which we then again solve by a variant of the \texttt{SMIP} Algorithm~\ref{algo:standardSMIP} to find the binary values for the rest of the original time horizon.
Our set of variables to fix for creating \instFixed is given in Line~\ref{line:AH_defineVariablesToFix} of Algorithm~\ref{algo:aggregatedHorizon}.

Since we do not fix all the variables, the aggregated binary assignments created by solving \instAgg might be identical with respect to the variable fixations.
In order to avoid unnecessary computation, we remove duplicated binary assignments by calling the function \texttt{removeDuplicateBAsWithRespectToVarsToFix} in Line~\ref{line:AH_removeDuplicatedBasBasedOnFixedVars} of Algorithm~\ref{algo:aggregatedHorizon}.
Now, we create for each remaining binary assignment a variant \instFixed of the original problem instance \inst by calling function \texttt{createProblemInstanceWithPartiallyFixedBinaries} in Line~\ref{line:AH_fixAggregatedSolutionInFullHorizonInstance}.
When denoting the binary assignment by $B$ and the value of a binary variable $x$ in $B$ by $\varsigma(x,B)$, we perform the following fixations:
\begin{align*}
    &\forall t^\text{agg}\in\{0,1,\dots,\paramAggHorizon\} \text{ and variables }x\in\text{varTypesToFix}:
    & x_{\vartheta(t^\text{agg})} &= \varsigma(x_{t^\text{agg}},B) \\
    &\forall t^\text{agg}\in\{1,\dots,\paramAggHorizon\} \text{ and variables }x\in\text{varTypesToFix}\\
    &\quad\text{with } \varsigma(x_{t^\text{agg}-1},B) = \varsigma(x_{t^\text{agg}},B), \\
    &\quad \forall t\in\{\vartheta(t^\text{agg}-1)+1,\dots,\vartheta(t^\text{agg})-1\}:
    & x_t &= \varsigma(x_{t^\text{agg}},B)
\end{align*}
In words, we do two types of fixations:
First, we fix all binary variables of the corresponding types whose time steps are part of the aggregated time horizon to their corresponding solution values defined by $B$.
Second and for all time steps in between those time steps from the aggregated horizon, we check if a variable has the same solution value with respect to $B$ at the two surrounding time steps from the aggregated horizon.
If it does, we also fix them to exactly that solution value.
Otherwise, we do not fix its value.
This means that we allow a variable value change between two aggregated time points to happen at any point in time in the corresponding original time interval.

After creating the instance \instFixed, we again solve it by the \texttt{SMIP} Algorithm~\ref{algo:standardSMIP} using the parameter \paramMaxItFull to specify its maximum number of iterations.
However, we make two adjustments to the algorithm.
First, we replace the \texttt{applyInitialLinearization} function with a newly defined function \texttt{initializeLinearizationBasedOnAggSolution}.
In there, we still use the constant velocity approximation but determine the coefficient $\alpha^q_{x,a,t}$ for each time $t$, pipe $a=(\ell,r)$ and end node $x\in\{\ell,r\}$ based on the pressures and flows of the aggregated solution that created the binary assignment.
If $t=\vartheta(t^\text{agg})$ for some $t^\text{agg}$ in the aggregated time horizon, we use the pressures and flows of that aggregated time step.
Otherwise, we linearly interpolate the values.
As a second adjustment, we provide a partial starting solution to each attempt to solve a linearized MIP, from which the solver tries to build a complete solution for the problem.
We represent this adjustment by replacing the original \texttt{solve} function with the function \texttt{solveMIPUsingPartialStartSolution}.
As the partial start solution, we provide the values for all binary variables \activeFlowDirI{f}, which determine the activity of a flow direction $f$, and \modeOp{a}, which represent the mode for all fence node valves $a$.
The values for each entity and time step are created based on the given binary assignment using the same algorithm as in function \texttt{createProblemInstanceWithPartiallyFixedBinaries}.
Since both variables do not have a meaningful value for the initial time step, we use for all original time steps, which exist before the first time step with a corresponding aggregated time step, the variable value of that aggregated time step.

As a result of calling the adjusted \texttt{SMIP} in Line~\ref{line:AH_sMIPCallWithFixedVars}, we obtain a set of binary assignments for the original problem instance \inst.
The union of these for all aggregated binary assignments is returned as the result of Algorithm~\ref{algo:aggregatedHorizon}.
Note that the binary assignments obtained by completing different aggregated binary assignments are also different due to the guaranteed variable fixations on those original time steps with a corresponding aggregated time step.

\subsubsection{Aggregated horizon heuristic using rolling horizon}
As the final heuristic based on reduced time horizons, we propose a combination of the previous two algorithms and is given as Algorithm~\ref{algo:aggregatedRollingHorizon}.
We refer to it as \texttt{ARH}.

\begin{algorithm}[ht]
  \SetKwInOut{Parameters}{Parameters}
  \DontPrintSemicolon
  \KwData{An instance \inst of problem \MINLP}
  \Parameters{%
    Size of the aggregated horizon \paramAggHorizon\\
    Size of the rolling horizon \paramRolHorizon\\
    Maximum number of iteration for the aggregated horizon MIP \paramMaxItAgg\\
    Maximum number of iteration for the full horizon MIP \paramMaxItFull\\
    Absolute convergence tolerance \paramMaxAbsTol\\
    Relative convergence tolerance \paramMaxRelTol
  }
  \KwResult{A list of valid binary assignments for \inst}
  % functions
  \SetKwFunction{adjustAhByUsingRollingHorizonToSolveAggregatedMIP}{adjustAhByUsingRollingHorizonToSolveAggregatedMIP}
  \SetKwFunction{ahUsingRh}{ahUsingRh}
  % algo start
  \ahUsingRh $\gets$ \adjustAhByUsingRollingHorizonToSolveAggregatedMIP{\paramRolHorizon}\;
  \KwRet{\ahUsingRh{\inst,\paramAggHorizon, \paramMaxItAgg, \paramMaxItFull, \paramMaxAbsTol, \paramMaxRelTol}}\;
  \caption{Sequential aggregated and rolling horizon heuristic algorithm (\texttt{ARH})}
  \label{algo:aggregatedRollingHorizon}
\end{algorithm}

The heuristic is nearly equivalent to the aggregated horizon heuristic given as Algorithm~\ref{algo:aggregatedHorizon}.
The only difference is that we use in Line~\ref{line:AH_solveAggregatedProblem} the rolling horizon Algorithm~\ref{algo:rollingHorizon} instead of the sequential mixed-integer programming Algorithm~\ref{algo:standardSMIP} to solve the problem instance \instAgg on the aggregated time horizon.
We represent this in Algorithm~\ref{algo:aggregatedRollingHorizon} by the function
\texttt{adjustAhByUsingRollingHorizonToSolveAggregatedMIP}, which returns a copy of Algorithm~\ref{algo:aggregatedHorizon} adjusted by the mentioned change.
This adjusted algorithm is then used to solve the given problem instance \inst.
Note that we are given the parameter \paramRolHorizon, which specifies the length of the rolling horizon.

\subsection{A dynamic node limit}\label{sec:dynamic_node_limit}
When testing the above-given heuristics, we observed that during the solving process of the linearized MIP models, we find good primal solutions fast and then spend much time proving their optimality.
This is typical behavior for branch-and-bound-based MIP solvers according to \citep{BerHenKoc2018}.
However, for our heuristics, we are only interested in the optimal (or at least a good) primal solution for each MIP.
The corresponding dual bound is not relevant to us.
Therefore, there is the potential to save time if we can end the solving process at the point at which the solver has found a good solution and changed its focus to improving the dual bound.

To implement the described behavior, we define a maximum limit $L$ on the number of processed branch-and-bound nodes, which depends on the so far found feasible solutions and can be applied to the MIP solving processes in the heuristics.
The general idea is to increase the node limit as long as we find new best primal solutions that significantly improve the primal bound.
If this is no longer the case, we hit the dynamic node limit and stop the run with the current best solution.

At the start of the solving process, we set the limit to $L=\infty$.
Whenever we find a new best feasible solution with objective function value $O$ at node $N$, we update the node limit dynamically as follows:
Let $O^\text{last}$ be the objective function value of the last feasible solution that updated the node limit and let $O^\text{last}=\infty$ if no feasible solution has been found so far.
Then update can be stated as
\begin{align*}
    \text{If } \frac{O^\text{last} - O}{O^\text{last}} > \paramImprovement : \qquad
    L &= \max(\paramNodeLimitMin, \paramNodeLimitRel \cdot N) \\
      O^\text{last} &= O.
\end{align*}
In words, we update the node limit $L$ whenever there is a significant relative improvement of at least $\paramImprovement\in(0,1)$ of the primal bound compared to its value at the last node limit update.
The new limit has a minimal value of \paramNodeLimitMin and is otherwise set to \paramNodeLimitRel times the current number of processed branch-and-bound nodes.
The minimal node limit avoids a premature termination of the solver and should cover the early stage of the solving process, in which the relative increments in the number of nodes are still relatively high.

\section{Computational Experiments} \label{sec:computational_experiments}
The performance of the presented algorithms is tested by conducting a series of computational experiments.
The corresponding set of instances is based on a real-world network with corresponding initial states as well as inflow and pressure patterns at the boundary nodes, which are provided by our project partner OGE~\citep{OGE}.
From this data we choose 40 especially challenging instances, by finding those with the highest amount of demand changes.
After giving a description of the corresponding network features and describing the instance generation process, we specify the setup of the experiments as well as the used parameters.
Then, we present the results of the two conducted experiments:
First, we compare the results of the general sequential mixed-integer programming (\texttt{SMIP}) Algorithm~\ref{algo:standardSMIP} used in Algorithm~\ref{algo:baseAlgo} against those of a black-box global MINLP solver.
Here, we focus purely on the solution quality to investigate if the algorithm is able to yield solutions close to the proven optimum.
In a second experiment, we evaluate the general performance of the reduced time horizon heuristics presented in Section~\ref{sec:reduceHorizonHeuristics} while using the general \texttt{SMIP} approach as a reference.
\subsection{Instance generation}\label{sec:instances}
The instances for testing our algorithms are based on a section of the real-world network of our project partner OGE~\citep{OGE}, which was also used in \citep{HopHenLenKoc2021}.
It has been modified by replacing the original network station topology with the artificial model that was manually created by industry experts.
Furthermore, the network outside of network stations is aggregated using the methods presented in~\citep{Len2021}.

We choose the demand scenarios defining our instances from a historical time period of 12 consecutive months.
During this time, the network topology, the manually created station model, and the parameters of the applied network aggregation changed multiple times.
Hence, we are given a slightly different network topology for each instance.
An overview of the network characteristics is given in Table~\ref{tab:networkStatistics}, in which the minimal and maximal amounts for each element type are given.
Most noteworthy, the number of network stations changed from 7 to 10, as three regulators arrangements in the network were replaced by small stations.
We note that compared to the instances used in \citep{HopHenLenKoc2021}, a less aggressive network aggregation was used, resulting in more elements outside of network stations.
\begin{table}[htb]
    \setlength{\tabcolsep}{5pt} % space between columns, default value = 6pt
    \centering
    \begin{tabular}{lrrrrrrrrrrrr}
        % nNodes nEntries nExits nArcs nPipes nOrigValves nControlValves
        % nAuxShortCuts nAuxRegulators nAuxBiRegulators nAuxCompressors nStations
        & $|\setVertices|$ & $|\setEntries|$ & $|\setExits|$ & $|\setArcs|$
        & $|\setPipes|$ & $|\setValvesOrig|$ & $|\setRegulators|$ & $|\setShortcuts|$
        & $|\setRegulatorArcs|$ & $|\setBiRegulatorArcs|$ & $|\setCompressorArcs|$ & $|\setNetworkStations|$%
        \medskip\\
        min & 375 & 12 & 114 & 421 & 303 & 0 & 11 & 29 & 20 & 2 & 16 & 7 \\
        max & 477 & 13 & 162 & 531 & 391 & 1 & 13 & 42 & 28 & 3 & 19 & 10
    \end{tabular}
    \caption{Network topology statistics. Since the network and manually created station model change over time, we give for each quantity the corresponding minimal and maximal amounts.}
    \label{tab:networkStatistics}
\end{table}

In Table~\ref{tab:stationStatistics}, we list parameters regarding the network stations.
If at least one of these values changed over time, we give two sets of values for a station representing the corresponding minimal and maximal amounts per quantity.
The network stations A to G are the same as in \citep{HopHenLenKoc2021}, while the newly introduced stations replacing the regulators are named H, I, and J.
\definecolor{myGray}{gray}{0.9}
\begin{table}[htb]
    \centering
    \begin{tabular}{lrrrrrrrrr}
        % name nFenceGroupNodes nAuxShortcuts nAuxRegulators nAuxBiRegulators nAuxCompressors % nStates nFlowDirs nConfigs nConflicts
        $i\in\setNetworkStations$ & $|\setFenceNodesNaviStation{i}|$ &
        $|\setShortcutsNaviStation{i}|$ & $|\setRegulatorArcsNaviStation{i}|$ &
        $|\setBiRegulatorArcsNaviStation{i}|$ & $|\setCompressorArcsNaviStation{i}|$ &
        $|\setSimpleStatesI{i}|$ & $|\setFlowDirectionsI{i}|$ &
        $\sum_{a\in\setCompressorArcsNaviStation{i}} |\setAuxCompressorConfigsI{a}|$ & $|\setAuxCompressorConfigConflictsI{i}|$%
        \medskip\\
        \rowcolor{myGray}
        A & 2 & 1 & 0 & 0 & 2 & 5 & 3 & 7 & 0 \\
        B$_\text{min}$ & 2 & 1 & 0 & 0 & 3 & 5 & 2 & 4 & 0 \\
        B$_\text{max}$ & 2 & 1 & 0 & 0 & 3 & 6 & 3 & 4 & 0 \\
        \rowcolor{myGray}
        C$_\text{min}$ & 4 & 2 & 6 & 0 & 1 & 3 & 4 & 1 & 0 \\
        \rowcolor{myGray}
        C$_\text{max}$ & 6 & 4 & 6 & 0 & 1 & 9 & 6 & 1 & 0 \\
        D$_\text{min}$ & 3 & 2 & 1 & 0 & 2 & 5 & 7 & 2 & 0 \\
        D$_\text{max}$ & 3 & 6 & 2 & 0 & 3 & 7 & 8 & 3 & 0 \\
        \rowcolor{myGray}
        E$_\text{min}$ & 5 & 1 & 5 & 0 & 1 & 9 & 11 & 2 & 0 \\
        \rowcolor{myGray}
        E$_\text{max}$ & 6 & 8 & 6 & 1 & 2 & 19 & 13 & 3 & 0 \\
        F$_\text{min}$ & 6 & 6 & 1 & 2 & 2 & 7 & 3 & 6 & 0 \\
        F$_\text{max}$ & 6 & 6 & 1 & 2 & 3 & 11 & 6 & 12 & 6 \\
        \rowcolor{myGray}
        G$_\text{min}$ & 10 & 15 & 6 & 0 & 5 & 25 & 12 & 18 & 17 \\
        \rowcolor{myGray}
        G$_\text{max}$ & 10 & 16 & 10 & 0 & 5 & 33 & 19 & 18 & 17 \\
        H & 2 & 0 & 1 & 0 & 0 & 1 & 1 & 0 & 0 \\
        \rowcolor{myGray}
        I & 3 & 1 & 1 & 0 & 0 & 1 & 1 & 0 & 0 \\
        J & 2 & 1 & 2 & 0 & 0 & 3 & 2 & 0 & 0
    \end{tabular}
    \caption{Network station statistics. Since the network and manually created station model may change over time, we give for each station and each quantity the corresponding minimal and maximal counts. If there are no changes in a station, we only give one set of values.}
    \label{tab:stationStatistics}
\end{table}

\paragraph{Artificial compressor model size}
\begin{table}[ht]
\centering
\begin{tabular}{lccc}
        & \# continuous variables & \# binary variables & \# constraints \\
  base  & 56 -- 68 & 61 -- 72 & 236 -- 278 \\
  new   & 14 -- 17 & 54 -- 65 & 151 -- 184
\end{tabular}
\caption{Artificial compressor model sizes of the base and our new model based on the used instance set. Numbers vary over time. The corresponding minimal and maximal numbers originate from the same point in the overall time period.}
\label{tab:compressorModelComparison}
\end{table}

Based on the specific network station characteristics in the given time period, we compare the actual sizes of the artificial compressors arc models used in the base model of \citep{HopHenLenKoc2021} and our newly introduced model from Section~\ref{sec:auxCompressor}.
We count the number of continuous and binary variables as well as the number of linear inequality constraints needed for each discrete time step.
For the variables, we included the activation and flow variable per arc as well as all additional variables needed specifically for artificial compressors.
The constraints we count for our model are those defined in Equations~\eqref{eq:artificial_compressor_model}, \eqref{eq:artificial_compressor_config_model}, and \eqref{eq:compressor_conf_conflict}, while we count all the constraints (31)-(42) from \citep{HopHenLenKoc2021} for the base model.
Since we count inequality constraints, equations are counted twice.
The result is given in Table~\ref{tab:compressorModelComparison}.
We see that our model is smaller regarding all three categories, as it uses a slightly smaller number of binary variables and considerably fewer continuous variables and inequality constraints.

\paragraph{Challenging demand scenarios}
To properly examine if the presented algorithms are suitable for application in time-critical industry environments, we search for particularly challenging instances.
In \citep{HopHenZitGot2020}, where the gas flow model of \citep{HopHenLenKoc2021} was used to conduct a case study on hydrogen transport, the authors suggest that the difficulty of the problem corresponds to the number of element control changes needed to fulfill the given demand scenario.
As the number of changes is not known a priori, we instead look for demand scenarios featuring large demand changes for single boundary nodes since we expect these changes to induce inevitable changes in some network elements' control.
From a cooperation with our project partner OGE, we attain historical network demands in hourly resolution over a time period of 12 months, covering the majority of the year 2020.
The time period is nearly consecutive, having only three periods of 2, 6, and 11 missing days due to problems during the data creation.

To find demand scenarios with large demand changes from within this time period, we first specify the desired time horizon length for our instances to be 12 and 24 hours.
For both, we then enumerate all time horizons of the corresponding length and sum the differences between the initial and end time point of the two quantities used as future demands in Section~\ref{sec:futureDemands}: The boundary node inflow values, which are given as norm volumetric flow $\nvFlow$ in $1000\text{m}^3/\text{h}$, and the pressure in bar for entries with positive inflow.
The norm volumetric flow is defined as $\nvFlow = \mFlow / \nDens$, where \nDens is the norm density depending on the gas mixture.
We then create a final score $\Delta$ for the amount of changes of a time horizon by adding the two sums, using the quantities' values for the given units, and weighting the pressure changes by a factor of 100.
Finally, we choose from each of the two lists of time horizons a set of 20 demand scenarios by iteratively selecting from the remaining time horizons the one with the highest $\Delta$ value, which does not overlap with any of the already selected ones.
Overlapping is hereby defined as sharing at least one full hour in the time horizon.
If one time horizon ends at time point $t$ and a second one starts at $t$, they do not overlap.

\begin{figure}[bth]
  \centering

  \begin{subfigure}[c]{0.48\textwidth}
  % trim options: <left> <lower> <right> <upper>
  \includegraphics[trim={0.0cm 2.5cm 3.0cm 4.0cm},clip,width=0.9\textwidth]{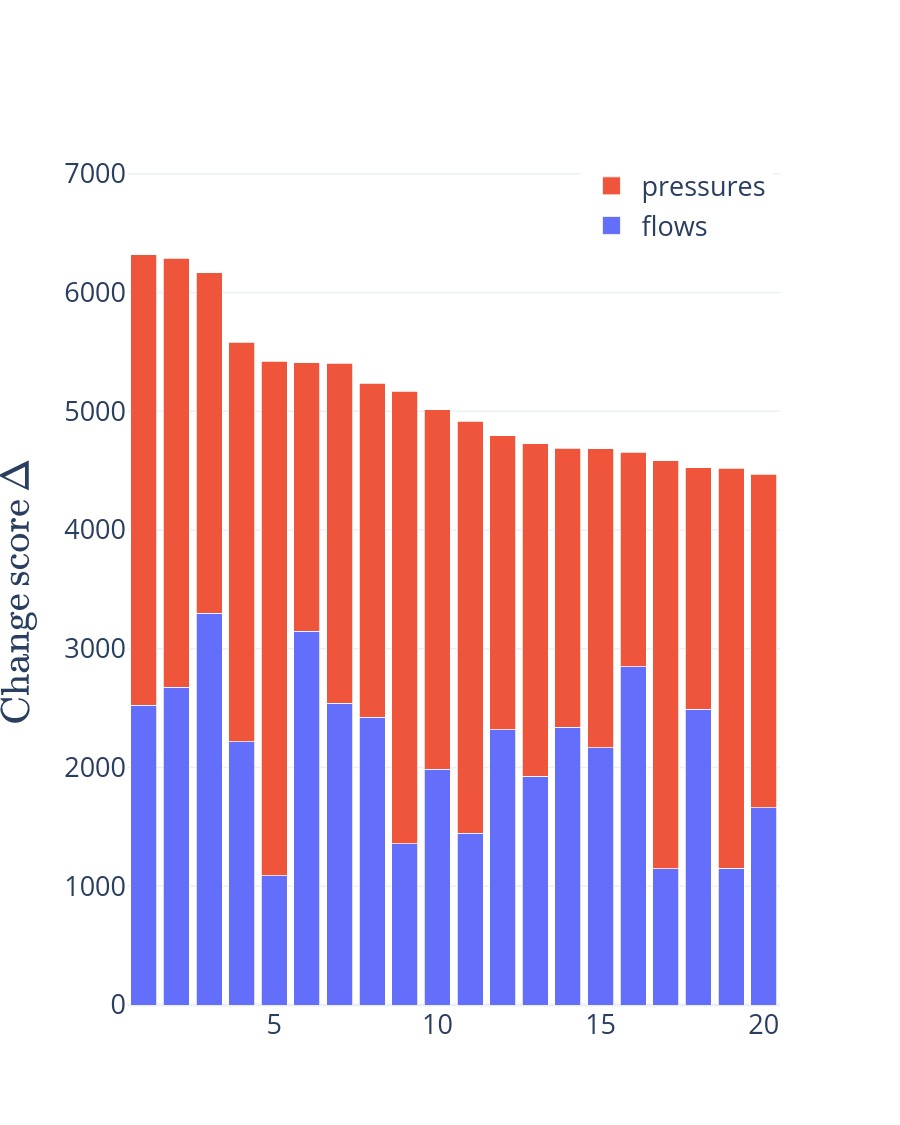}
  \subcaption{12h horizon}
  \label{fig:bestInstances12h}
  \end{subfigure}
  \begin{subfigure}[c]{0.48\textwidth}
  % trim options: <left> <lower> <right> <upper>
  \includegraphics[trim={0.0cm 2.5cm 3.0cm 4.0cm},clip,width=0.9\textwidth]{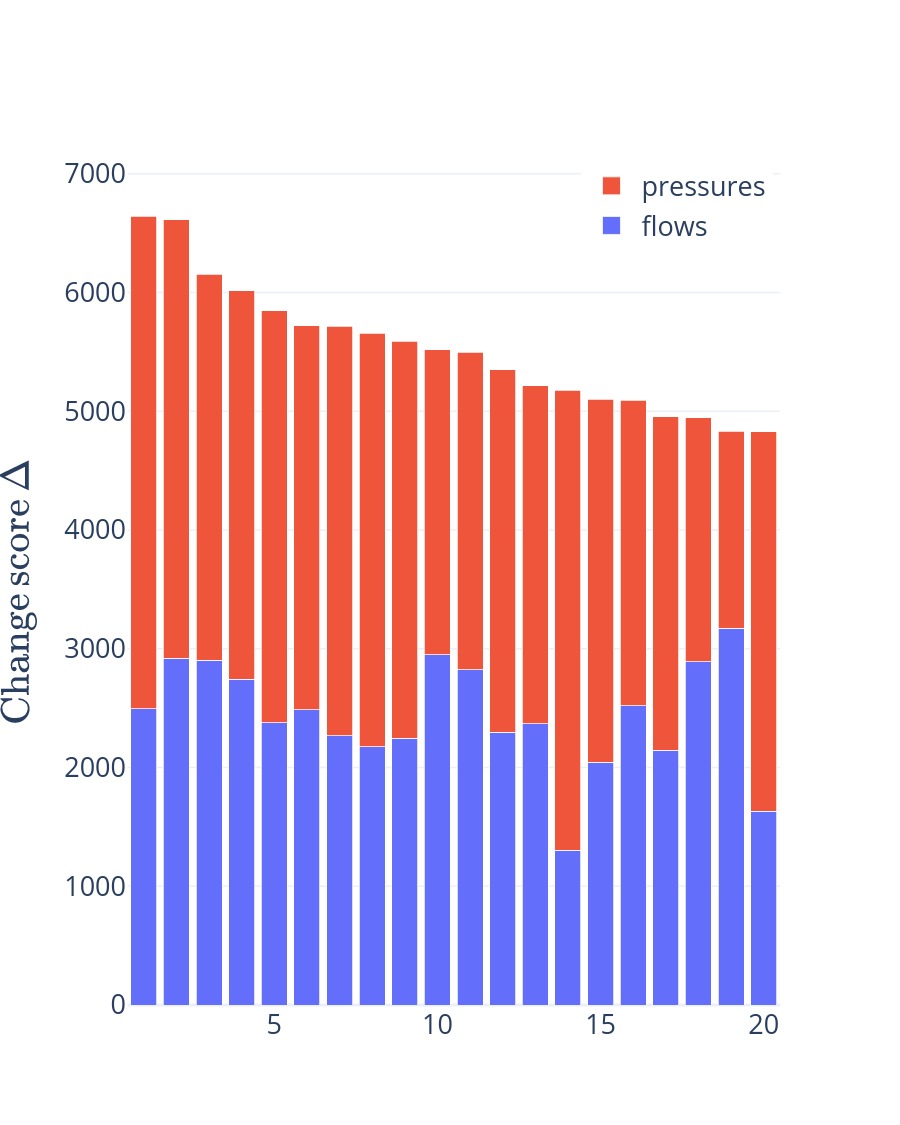}
  \subcaption{24h horizon}
  \label{fig:bestInstances24h}
  \end{subfigure}
  \caption{Change scores of the 20 non-overlapping demand scenarios with the highest change score $\Delta$ for each of the two lengths of the time horizon. Each bar represents one selected instance and is divided according to the pressure and flow share of $\Delta$, where the pressure is scaled by 100. The instances are sorted by change score. Figures created with \citep{Plo2015}.}
  \label{fig:bestGeneratedInstances}
\end{figure}

An overview of the $\Delta$ values for each of the 20 demands scenarios for each time horizon is given in Figure~\ref{fig:bestGeneratedInstances}.
The largest found score value is about 6600, which would be equivalent to a total amount of inflow changes of 6600\,1000m$^3/$h or a total amount of pressure changes of 66\,bar at entry nodes over the whole network.
The 40 demand scenarios, together with the corresponding network topology and initial network state at the start of the time horizon, define the 40 instances we use as our test set for the computational experiments.
When using a time horizon whose time steps are larger than 1 hour, we aggregate the hourly demands in the same way as in the function \texttt{createTimeAggregatedProblemInstance} used in the aggregated horizon heuristic presented in Algorithm~\ref{algo:aggregatedHorizon}.

\subsection{Setup and parameters}
In this section, we give an overview of the computational hardware setup and solver choice, the applied model adjustments to improve its numerical properties, and specify the parameter used for the computational experiments.
\subsubsection{Computational setup}
The computations were executed on a cluster using 4 cores and 40\,GB of RAM of a machine composed of two \emph{Intel Xeon Gold 5122} running at 3.60\,GHz.
The MIP problems are solved by \textsc{Gurobi} in version 9.5.1~\citep{Gur2021}, for which we change the general solution strategy to focus on finding primal solutions by setting the \texttt{MIPFocus} parameter to 1.
As the NLP solver, we used \textsc{Ipopt} in version 3.14.3~\citep{WacBie2006}.
Finally, we used \textsc{Scip} in version 8.0.0~\citep{BesBesCheChm2021} as the global MINLP solver.
For \textsc{Scip}, we change the feasibility tolerance for constraints \texttt{numerics/feastol} to $10^{-3}$, which is even higher than the corresponding default value of $10^{-4}$ for \textsc{Ipopt}, to improve the acceptance rate for start solutions.
We accessed both \textsc{Ipopt} and \textsc{Scip} via \textsc{Gams} in version 38.1.0~\citep{GAM2022}.
\subsubsection{Coefficient scaling and rounding}
To improve the numeric properties of the models, we use similar ranges for the continuous variables by scaling them such that pressure variables are used in bar and flow values in kg/s.
Furthermore, we apply the following modification to the pipes' Continuity equation~\eqref{eq:fd_discrete_continuity}, their Momentum equation~\eqref{eq:fd_discrete_momentum} in the original or the linearized variant, and the linearized maximum power constraint of artificial configurations \eqref{eq:maxPowerOfConfig}:
First, we set all coefficients to zero, whose absolute value divided by the largest absolute value coefficient is smaller than $10^7$.
Afterwards, we scale all coefficients of the corresponding constraint such that the smallest non-zero coefficient has an absolute value of 1.
\subsubsection{Parameters}
In the following, we list the parameters used in the model, in the heuristic algorithms, and as resources limits for the different solvers.
\paragraph{Model parameters}
% slack
For the model, we specify the maximum slack parameters to be 2\,bar for the pressure slack.
In the case of the inflow, we restrict the possible maximal deviation from the demands to $50\%$.
Hence, we use, instead of a fixed upper bound \slackQUB on the inflow slack, an upper bound specific to each boundary node and time step of
\begin{equation*}
    \slackQUB_{v,t} = 0.5\cdot|\demandInflow{v,t}|
\end{equation*}
for all boundary nodes $v\in\setBoundaryNodes$ and future time steps $t\in\setTimestepsNoZero$.
Furthermore, we set the minimum flow value used to determine the absolute velocity in the initial linearization of the friction term to $q^\text{min} = 100 \cdot 1000\text{m}^3/\text{h}$.

% objective function coefficients
We define the different weights used in the objective function defined in Equation~\eqref{eq:objective} as
\begin{align*}
 \costSlackQ &= 100 \cdot 3.6 / \nDens & \costSlackP &= 1000.0 & \costUnitRunning &= \costArtificialLinkChange = 50.0 \\
 \costVaModeChg &= \costRgModeChg = 500.0 & \costRgPLChg &= \costRgPRChg = \costFenceNodePChg = 10.0 & \costRgQChg &= \costFenceNodeQChg = 1.0 \cdot 3.6 / \nDens \\
 && \costSimpleStateChangeI{s}&\in[0,4000].
\end{align*}
Note that the multipliers for the flow change and the hourly flow slack are given with respect to norm volumetric flow \nvFlow in $1000\text{m}^3/\text{h}$.
Hence, we multiply them by a factor of $3.6 / \nDens$ to get the corresponding mass flow equivalent.
The costs for turning on a simple state of a network station are manually created by the industry experts at OGE, which is why we can only give the range of values for them.
\paragraph{Algorithm parameters}
Each heuristic algorithm defined in Section~\ref{sec:mip_based_heuristics} is given a set of parameters to adjust them.
First, we specify the converge tolerances for the friction term linearization as
\begin{align*}
  \paramMaxRelTol &= 0.001 & \paramMaxAbsTol &= 0.01\,\text{bar}.
\end{align*}
Next, the maximum number of sequential mixed-integer programming iterations we use is
\begin{align*}
  \paramMaxIt &= \paramMaxItAgg = 5 & \paramMaxItFull &= 3.
\end{align*}

Finally, we do not fix the parameters defining the reduced time horizon size globally but define different variants of the heuristics algorithms based on these.
Hence, we denote the rolling horizon heuristic of Algorithm~\ref{algo:rollingHorizon} with a rolling horizon size of $\paramRolHorizon=X$ by \texttt{RXH}, the aggregated horizon heuristic of Algorithm~\ref{algo:aggregatedHorizon} with aggregated horizon size of $\paramAggHorizon=X$ by \texttt{AXH}, and the aggregated and rolling horizon heuristic of Algorithm~\ref{algo:aggregatedRollingHorizon} using an aggregated horizon of length $\paramAggHorizon=X$ and a rolling horizon of length $\paramRolHorizon=Y$ by \texttt{AXRYH}.
\paragraph{Resource limits}
As final parameters, we define the resource limits of the different solving processes, which are given in Table~\ref{tab:resourceLimits}.
There we distinguish five different model variants: MIP with the dynamic node limit defined in Section~\ref{sec:dynamic_node_limit}, MIP without this node limit, an NLP run to complete binary assignments, a \emph{fast} NLP variant having a time limit of 5\,minutes, and finally, the MINLP run.
The given parameters for the dynamic node limit specify that the smallest node limit to stop is 1000 and that a significant primal solution upgrade featuring an objective function value improvement of at least 2\% increases the node limit to 150\% of the current number of processed nodes.
Note that all variants are additionally restricted by the general memory limit of 40\,GB.

\begin{table}[ht]
\centering
\begin{tabular}{lrrrrrrr}
  Variant           & Rel. Gap  & Time limit        &
  \paramNodeLimitMin & \paramNodeLimitRel & \paramImprovement \\
  MIP w/node limit  &       0.0 & $\infty$          & 1000 & 1.5 & 0.02 \\
  MIP wo/node limit & $10^{-2}$ & $86400\,\text{s}$ & -- & -- & -- \\
  NLP               &        -- & $\infty$          & -- & -- & --  \\
  NLP fast          &        -- & $300\,\text{s}$   & -- & -- & --  \\
  MINLP             & $10^{-6}$ & $600000\,\text{s}$& -- & -- & --
\end{tabular}
\caption{Optimality conditions and resource limits for the different model variants. We distinguish solving a MIP with and without the dynamic node limit defined in Section~\ref{sec:dynamic_node_limit}, solving an NLP with and without a time limit, and solving an MINLP.}
\label{tab:resourceLimits}
\end{table}
\subsection{Solution quality of general \texttt{SMIP}}
In our first experiment, we check the quality of solutions obtained by using the general sequential mixed-integer programming \texttt{SMIP} Algorithm~\ref{algo:standardSMIP} to determine binary assignments in Algorithm~\ref{algo:baseAlgo} and complete these to a full MINLP solution by an NLP solver.
The lower bound to compare the solutions against is obtained via a global MINLP solver.
Due to the size of our network and the general complexity of the problem, we were only able to obtain reasonable lower bounds for the smallest time horizon of size 1.
Depending on the actual instance, the corresponding MINLP model has the following sizes (rounded to 2 significant digits):
\begin{align*}
    \text{\# variables}&\in[2300,2900] & \text{\# binary variables}&\in[690,870] \\
    \text{\# constraints}&\in[3100,3800] & \text{\# non-linear constraints}&\in[300,390].
\end{align*}
Since the overall run time of the solving process has no priority for this experiment, we use for the single MIP models the parameter variant without the heuristic node limit, and for the NLP the general variant without a time limit, see Table~\ref{tab:resourceLimits}.
Also, we use the solution obtained from the \texttt{SMIP} algorithm as a starting solution for the MINLP.

\begin{figure}[thbp]
    \centering
    % trim options: <left> <lower> <right> <upper>
    \includegraphics[trim={0.0cm 1.5cm 2.5cm 4.0cm},clip,width=0.95\textwidth]{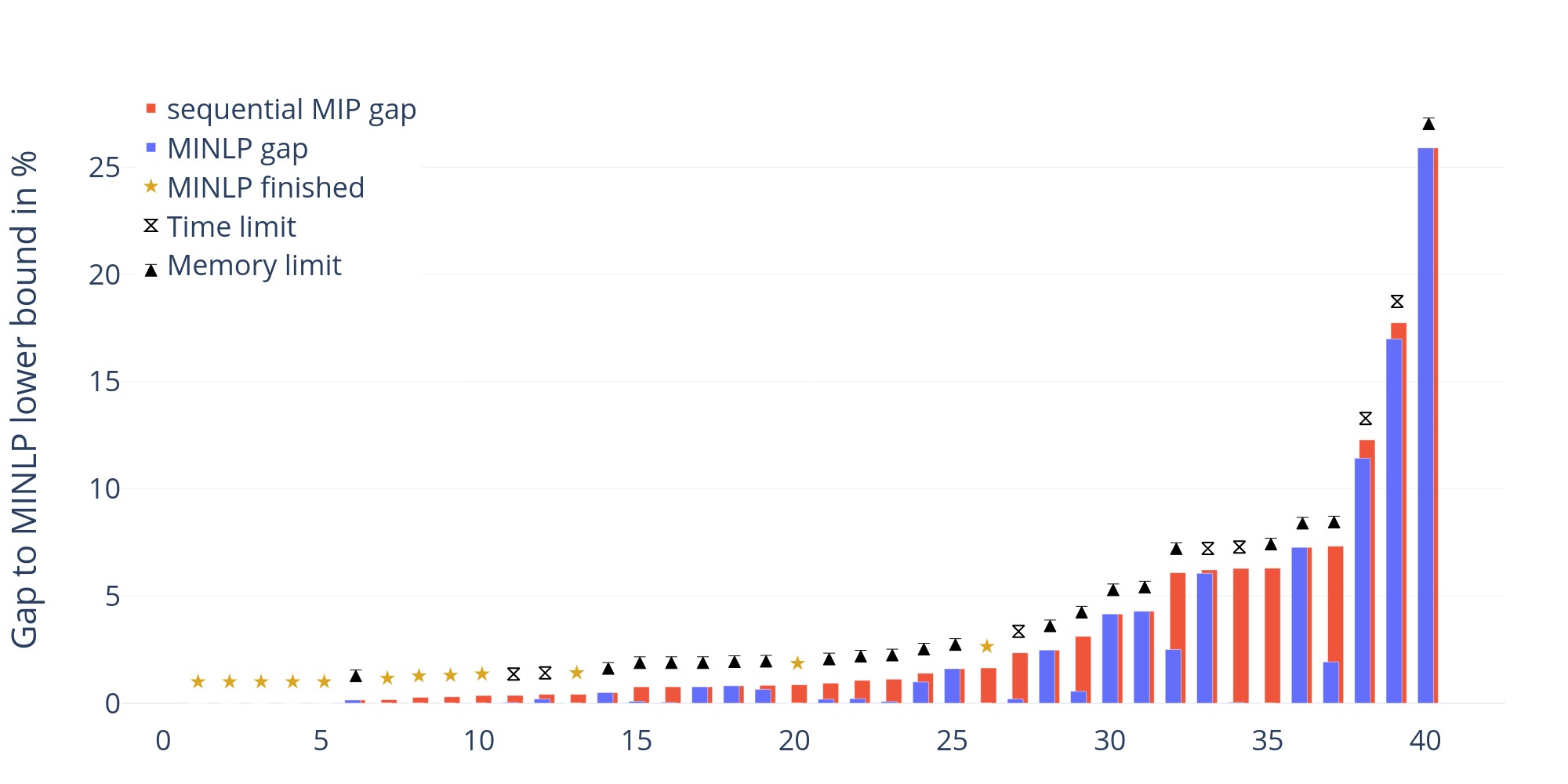}
    \caption{Comparison of the gap of the solutions created by the \texttt{SMIP} Algorithm~\ref{algo:standardSMIP} used in Algorithm~\ref{algo:baseAlgo} with respect to the lower bound found in the MINLP run and the gap of the MINLP itself. Each pair of bars represents one instance, sorted by the gap of the \texttt{SMIP} solution. The values are given in percent. For the instances marked by a star, the MINLP proved optimality. For those instances marked by an hourglass, the MINLP time limit was hit, while for instances marked by a black triangle pointing upwards to a horizontal bar, the MINLP hit the 40\,GB memory limit. Figure created with \citep{Plo2015}.}
    \label{fig:minlpVsSMipGap}
\end{figure}

The result is given in Figure~\ref{fig:minlpVsSMipGap}.
First, we notice that for only 12 out of the 40 instances, the MINLP was solved to optimality.
In 7 cases, the time limit of roughly one week was hit, while for the remaining 21 instances, the memory limit of 40\,GB was exceeded.
As a consequence, we do not always have the best possible lower bound available.
However, the average gap of 3.2\% of the \texttt{SMIP} solutions to that lower bound is still rather small and even closer to the average gap of 2.3\% of the MINLP itself.
Using the best found primal solutions of the MINLP as a reference for the highest possible lower bound values, the average gap for the \texttt{SMIP} solutions is 0.9\%.
Regarding run time, Algorithm~\ref{algo:baseAlgo} using the \texttt{SMIP} always finishes in less than 30 seconds, while the geometric mean run time of the MINLP is larger than a day.
In summary, the \texttt{SMIP} solutions are very close to the global optimum and can be computed in only a fraction of the global solver's execution time.
\subsection{Performance of reduced time horizon heuristics}
In our main experiment of the paper, we evaluate our proposed algorithm of using a reduced-time-horizon heuristic of Section~\ref{sec:reduceHorizonHeuristics} for creating binary assignments to be completed in Algorithm~\ref{algo:baseAlgo}.
As a reference approach, we use the \texttt{SMIP} for the binary assignment creation, which was shown to find close-to-optimal solutions in the previous section.
Not having to compare against a global MINLP solver enables us to solve instances with a time horizon consisting of 12 equally large time steps of 1-hour length for the 12-hour horizon and 2-hour length for the 24-hour horizon.
Depending on the actual instance, the MINLP model, which is never formulated and tackled as a whole model here, has the following sizes (rounded to 2 significant digits):
\begin{align*}
    \text{\# variables}&\in[25000,32000] & \text{\# binary variables}&\in[5600,7000] \\
    \text{\# constraints}&\in[35000,44000] & \text{\# non-linear constraints}&\in[3600,4700].
\end{align*}
After conducting some preliminary experiments, we determined for each of the reduced horizon heuristics the smallest time horizon lengths, which still give reasonably good results, and the largest time horizon lengths, at which the run times are still reasonably short.
As a result, we choose the following variants to test for each of the heuristics:
The rolling horizon heuristics from \texttt{R1H} to \texttt{R4H}, the aggregated time horizon heuristics \texttt{A2H}, \texttt{A3H}, \texttt{A4H}, and \texttt{A6H}, and for the last two aggregated horizon lengths also the heuristic using the rolling horizon heuristic for some of the sub-models with horizon lengths of 1 to 3, resulting in the variants \texttt{A4R1H}, \texttt{A4R2H}, \texttt{A4R3H}, \texttt{A6R1H}, \texttt{A6R2H}, and \texttt{A6R3H}.
For all variants, we solve the MIP models using the dynamic branch-and-bound node limit setting and solve the NLP using a time limit, see Table~\ref{tab:resourceLimits}.
For the reference \texttt{SMIP} approach, we ensure the highest possible solution quality by not using the node limit for MIP or the time limit for NLP.

\begin{figure}[tp]
  \centering

  \begin{subfigure}[c]{\textwidth}
  % trim options: <left> <lower> <right> <upper>
  \includegraphics[trim={0.0cm 0.0cm 3.0cm 3.0cm},clip,width=\linewidth]{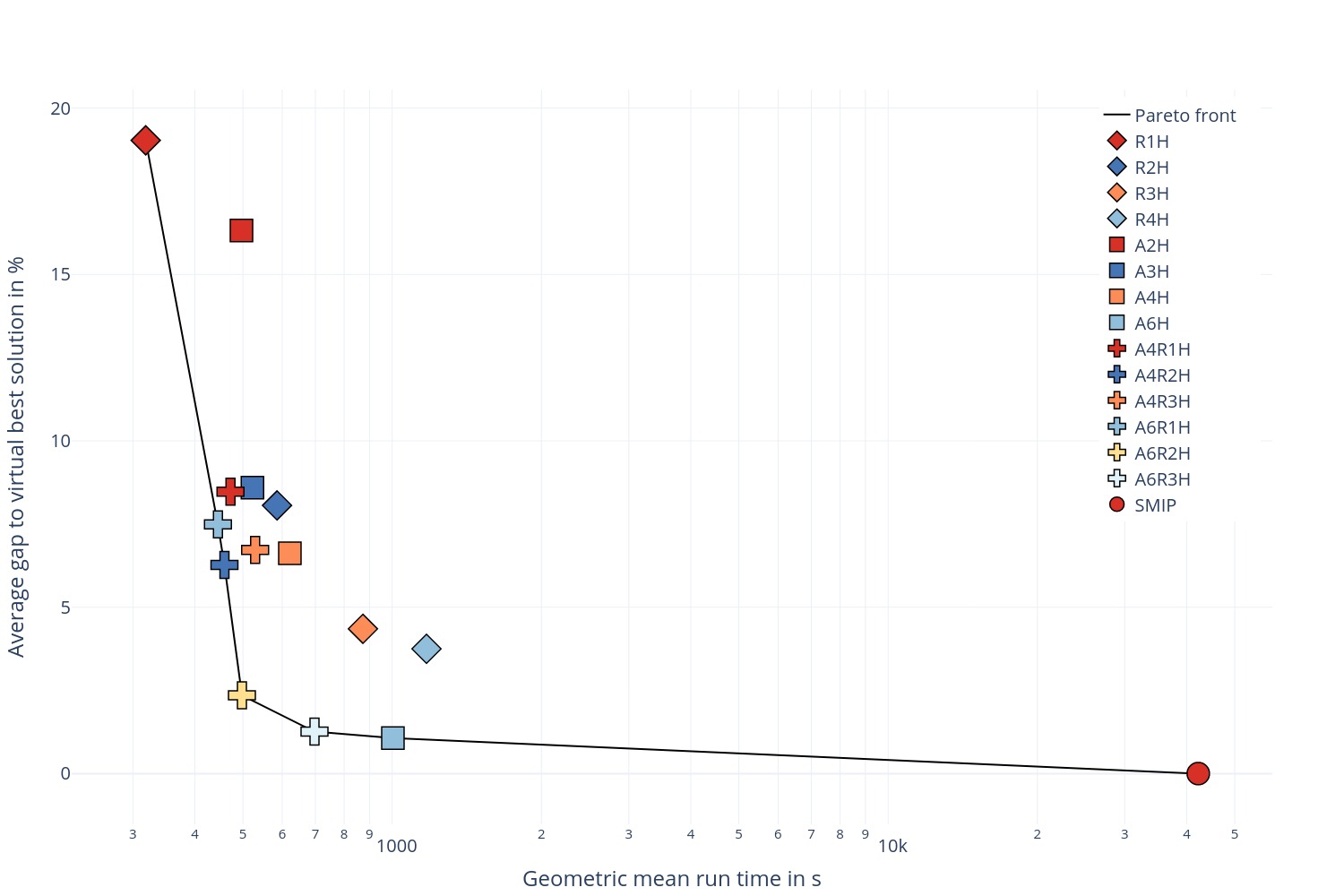}
  \subcaption{Individual heuristics}
  \label{fig:runTimeVsGap_geo_single}
  \end{subfigure}

  \begin{subfigure}[c]{\textwidth}
  % trim options: <left> <lower> <right> <upper>
  \includegraphics[trim={0.0cm 0.0cm 3.0cm 3.0cm},clip,width=\linewidth]{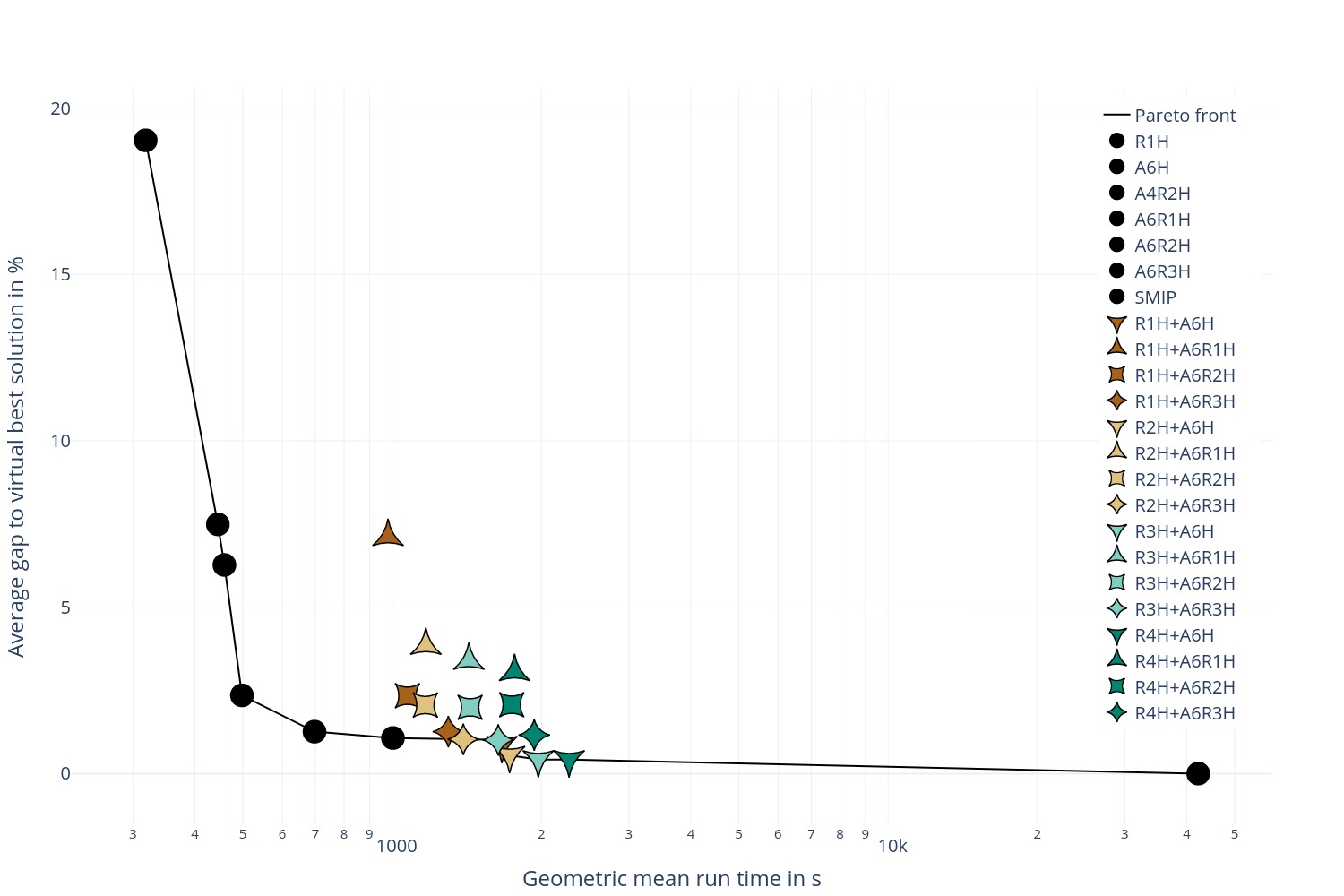}
  \subcaption{Combine heuristics and Pareto points of the individual heuristics}
  \label{fig:runTimeVsGap_geo_combi}
  \end{subfigure}
  \caption{Position of the heuristic variants in a geometric-mean-run-time vs. gap-to-virtual-best-solution plot. If a heuristic does not find a solution for some instance, we count the corresponding gap as 100\%. The Pareto efficient points are linked by a black line, forming the Pareto front. Figures created with \citep{Plo2015}.}
  \label{fig:runTimeVsGap_geo}
\end{figure}

\begin{figure}[tp]
  \centering
  \begin{subfigure}[c]{\textwidth}
  \includegraphics[trim={0.0cm 0.0cm 3.0cm 3.0cm},clip,width=\linewidth]{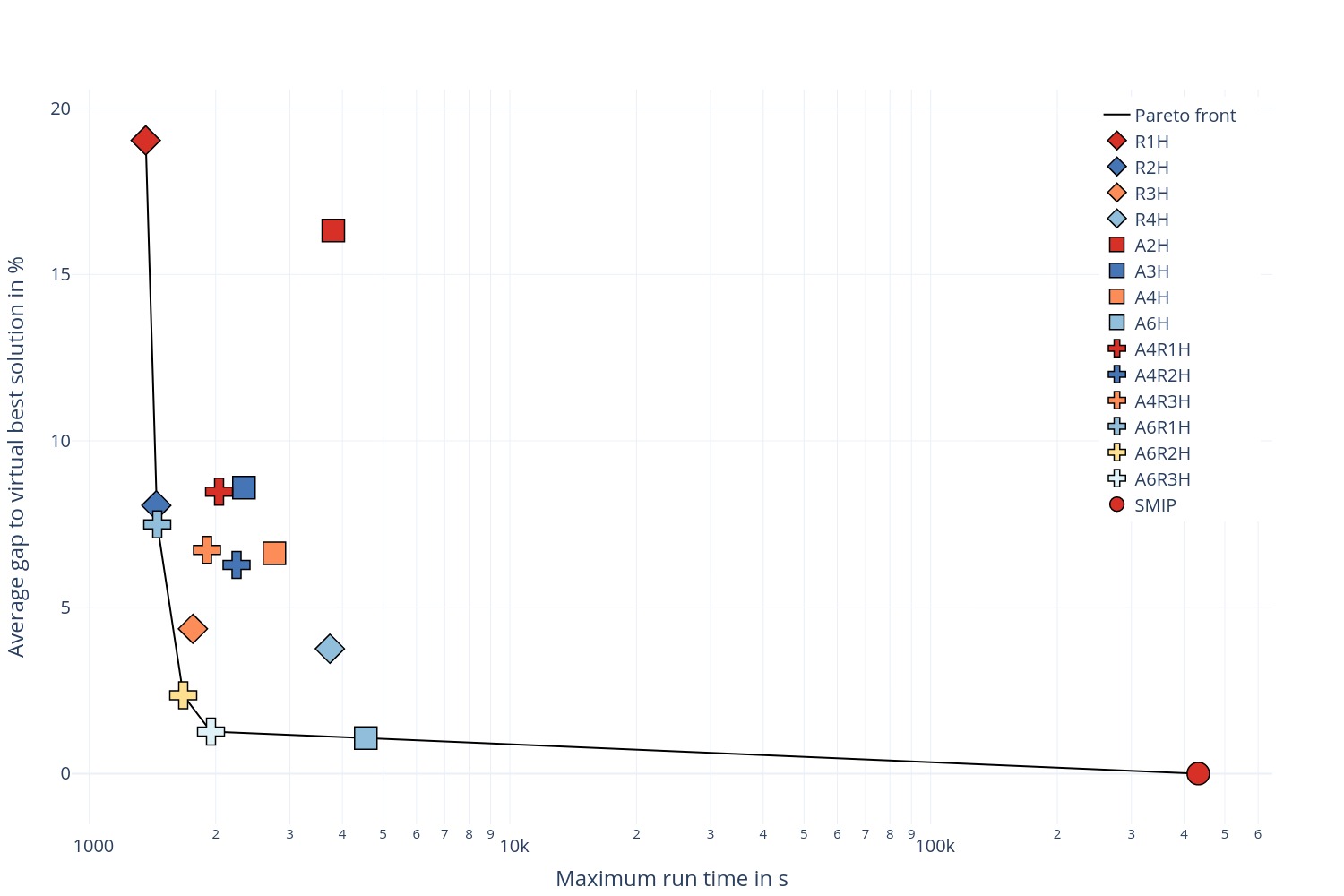}
  \subcaption{Individual heuristics}
  \label{fig:runTimeVsGap_max_single}
  \end{subfigure}

  \begin{subfigure}[c]{\textwidth}
  \includegraphics[trim={0.0cm 0.0cm 3.0cm 3.0cm},clip,width=\linewidth]{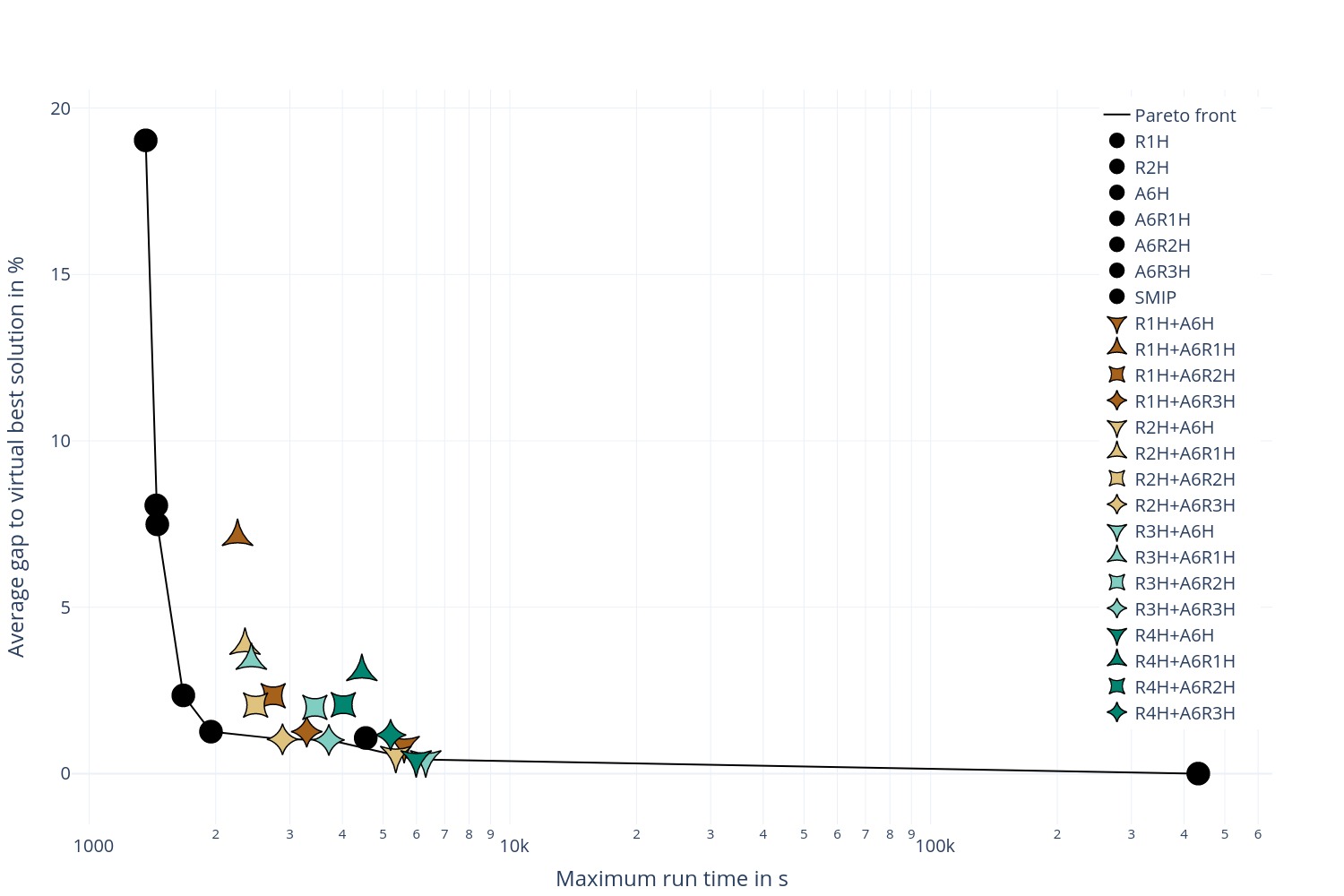}
  \subcaption{Combine heuristics and Pareto points of the individual heuristics}
  \label{fig:runTimeVsGap_max_combi}
  \end{subfigure}
  \caption{Position of the heuristic variants in a maximum-run-time vs. gap-to-virtual-best-solution plot. If a heuristic does not find a solution for some instance, we count the corresponding gap as 100\%. The Pareto efficient points are linked by a black line, forming the Pareto front. Figures created with \citep{Plo2015}.}
  \label{fig:runTimeVsGap_max}
\end{figure}

\begin{table}[thb]
\centering
\begin{tabular}{lrcrrrrcrrrrcrr}
 &&& \multicolumn{4}{c}{\# avg MIP runs} && \multicolumn{4}{c}{\# avg NLP runs} && \multicolumn{2}{c}{mean NLP time}\\
 & fail &  & total & opt & DNL & inf &  & total & conv & time & fail &  & by NLP & by inst\medskip \\
\rowcolor{myGray}
\texttt{R1H} & 4 &  & 46.73 & 46.62 & 0.00 & 0.10 &  & 3.38 & 3.12 & 0.05 & 0.20 &  & 161.6 & 307.4 \\
\texttt{R2H} & 1 &  & 46.08 & 45.35 & 0.70 & 0.03 &  & 3.67 & 3.58 & 0.03 & 0.07 &  & 132.8 & 418.6 \\
\rowcolor{myGray}
\texttt{R3H} & 0 &  & 44.00 & 41.42 & 2.58 & 0.00 &  & 3.88 & 3.75 & 0.05 & 0.07 &  & 130.3 & 474.9 \\
\texttt{R4H} & 0 &  & 38.25 & 30.80 & 7.45 & 0.00 &  & 3.77 & 3.70 & 0.03 & 0.05 &  & 120.8 & 430.9 \\
\rowcolor{myGray}
\texttt{A2H} & 1 &  & 10.40 & 8.82 & 1.48 & 0.10 &  & 4.55 & 4.45 & 0.05 & 0.05 &  & 86.6 & 258.3 \\
\texttt{A3H} & 1 &  & 11.22 & 8.82 & 2.33 & 0.07 &  & 4.65 & 4.65 & 0.00 & 0.00 &  & 78.3 & 246.6 \\
\rowcolor{myGray}
\texttt{A4H} & 0 &  & 11.68 & 8.90 & 2.67 & 0.10 &  & 4.55 & 4.55 & 0.00 & 0.00 &  & 77.5 & 278.0 \\
\texttt{A6H} & 0 &  & 12.62 & 9.15 & 3.42 & 0.05 &  & 4.90 & 4.85 & 0.00 & 0.05 &  & 76.6 & 315.6 \\
\rowcolor{myGray}
\texttt{A4R1H} & 0 &  & 24.55 & 24.48 & 0.07 & 0.00 &  & 4.78 & 4.78 & 0.00 & 0.00 &  & 71.4 & 278.6 \\
\texttt{A4R2H} & 0 &  & 19.02 & 18.62 & 0.40 & 0.00 &  & 4.47 & 4.42 & 0.00 & 0.05 &  & 75.5 & 245.8 \\
\rowcolor{myGray}
\texttt{A4R3H} & 0 &  & 15.03 & 12.57 & 2.42 & 0.03 &  & 4.47 & 4.45 & 0.03 & 0.00 &  & 75.9 & 260.5 \\
\texttt{A6R1H} & 1 &  & 33.77 & 33.75 & 0.00 & 0.03 &  & 5.35 & 5.33 & 0.00 & 0.03 &  & 74.4 & 285.5 \\
\rowcolor{myGray}
\texttt{A6R2H} & 0 &  & 28.12 & 27.38 & 0.72 & 0.03 &  & 4.75 & 4.75 & 0.00 & 0.00 &  & 70.3 & 266.8 \\
\texttt{A6R3H} & 0 &  & 25.23 & 22.52 & 2.70 & 0.00 &  & 5.08 & 5.05 & 0.00 & 0.03 &  & 76.9 & 324.8
\end{tabular}
\caption{Statistics for the used heuristic variants.
The first two columns state the name of the heuristic and the number of instances for which it fails to find a feasible solution.
The next four columns state the average number of MIP models solved per instance, followed by the corresponding average number of those MIPs finished with an optimal solution, finished due to reaching the dynamic node limit, or finished with proven infeasibility of the corresponding model.
Next, we have a similar block of four columns for the average number of NLP models solved per instance.
First, the total average is stated (which is equal to the average number of binary assignments found), then the average number of converged NLPs, the average number of NLP finishing with a feasible solution at the given time limit, and finally, the average number of failed attempts to find a feasible solution within the time limit.
The final two columns depict the geometric mean run time for the NLP in seconds, first for a single NLP run and second for the summed run time of all NLP runs per instance.
If a heuristic did not produce any binary assignments for an instance, the NLP run time was set to 1\,second.}
\label{tab:heuristicResultStats}
\end{table}

We present the results averaged over all instances for each heuristic variant in two pairs of figures.
Each figure shows the position of each variant with respect to the run time and the average gap to the virtual best solution.
In the Figures~\ref{fig:runTimeVsGap_geo}, the run time is determined as the geometric mean run time of all the instances and in the Figures~\ref{fig:runTimeVsGap_max}, we use the maximum run time.
If a heuristic failed to find a solution for a specific instance, we used a gap value of 100\%.
In Table~\ref{tab:heuristicResultStats}, we gave an overview of the number of fails as well as the average results for the single MIP solves and NLP solves per heuristic variant.
In the figures, we connected those heuristics being Pareto efficient by a black line, where Pareto efficient heuristics are those not having another heuristic with a smaller average gap and, at the same time, a shorter run time value.

In addition to the individual heuristics displayed in the upper Figures~\ref{fig:runTimeVsGap_geo_single} and~\ref{fig:runTimeVsGap_max_single}, we also present in the lower Figures~\ref{fig:runTimeVsGap_geo_combi} and~\ref{fig:runTimeVsGap_max_combi} results for the \emph{combination} of two different heuristic variants.
For each combination and instance, we use the best solution found in any of the two heuristics and add up the overall run time of both approaches.
As a reference, we also included all Pareto efficient individual heuristics from the corresponding upper pictures.
For the combination, we always choose a pair of a pure rolling horizon heuristic and a heuristic using an aggregated time horizon.
We do this since we expect the two fundamentally different approaches to have a higher probability complementing each each other in finding high-quality solutions.
For the combination, we chose all four rolling horizon heuristics.
For the time-aggregation-based heuristics, we used those that are part of the Pareto front for the geometric mean and the maximum run time figures, i.e., \texttt{A6H}, \texttt{A6R1H}, \texttt{A6R2H}, and \texttt{A6R3H}.
Note that the combined heuristics in the figures have the same color for using the same rolling horizon heuristics and the same symbol for using the same time-aggregation-based heuristic.

The figures show that the reference \texttt{SMIP} approach is the overall best approach in terms of solution quality with an average gap of 0.0\%, as it yields the best solution for all of the instances.
However, it is much slower than every of the reduce time horizon heuristics both on average and regarding the maximum run time.
This is mainly due to the fact that, on average, 1.5 out of the at most 5 solved MIP models reached the 1-day time limit before proofing optimality.
For the reduced time horizon heuristics, the solution quality is, in general, very good:
The average gap to the virtual best solution is below 20\% for all considered variants and even below 10\% for all but two.
The best average gap apart from the general \texttt{SMIP} approach was achieved by the \texttt{A6H} heuristic with a gap of 1.07\% and 17/40 instances, for which the best solution has been found.
Two very close variants are \texttt{A6R2H} and \texttt{A6R3H}, which are close in terms of the average gap but have better run time values.
Hence, the summary for the individual variants is to choose an aggregated horizon heuristic with a long time horizon and either solve it directly for the best results or internally use the rolling horizon approach for shorter run times at the cost of a slightly worse solution quality.
When comparing the geometric mean and the maximum run time, the pictures are very similar, with only two heuristic variants being in only one of the corresponding Pareto fronts.
The maximum run time values are between 2 to 5 times higher for all the heuristics except \texttt{A2H}, where the factor is over 7.
For the two variants \texttt{A6R2H} and \texttt{A6R3H} with the lowest run times of those with an average gap below 3\%, the maximum run time is about half an hour, which is very fast for these types of problems. % we can not really proof this statement ...

Using two heuristics in combination improves the overall solution quality as expected, with 6 of the 16 combinations having an average gap better than the best individual solution and all of the combinations involving \texttt{A6H} reaching an average gap below 1\%.
However, this comes at the cost of rather large run time increases.
Only the fastest of the combined heuristics has a better geometric run time than the one of the slowest Pareto efficient reduced time horizon heuristic.
However, in that larger run time segment, several Pareto efficient combinations of the \texttt{A6R3H} and \texttt{A6H} heuristic combined with the rolling horizon heuristics can be found for both run time metrics.
The best of them decrease the gap to the \texttt{SMIP} variant to a minimum of 0.42\% for the two combinations of \texttt{R3H} and \texttt{R4H} with \texttt{A6H} while still being much faster than \texttt{SMIP} itself.

The reason for the overall bad run time for the combination of two heuristics is the rather long solve time for the NLP model, as displayed in the last two columns of Table~\ref{tab:heuristicResultStats}.
For heuristics like \texttt{R1H} and \texttt{R2H} that determine binary assignments very fast, the geometric mean NLP time per instance is equal to 96.5\% and 71.5\% of the geometric mean overall run time, respectively.
Hence, testing more binary assignments for yielding good overall solutions to the problem is apparently a considerable time investment.

To compare our algorithm against a procedure similar to the one proposed in \citep{HopHenLenKoc2021}, we also tested the \texttt{SMIP} approach with only a single iteration, i.e., $\paramMaxIt = 1$.
This is close to the algorithm proposed there, as the linearization used for the model determining the binary solution values is only updated in the rare case of failure in the subsequent procedure searching for a non-linear solution for the found binary assignment.
The \texttt{SMIP} heuristic with $\paramMaxIt = 1$ has an average gap of 42.89\% to the virtual best solution and is, therefore, worse than all of our proposed heuristic variants.
For 14 out of the overall 40 instances, it reached the time limit of 1 day.
This leads to a geometric mean run time of 9541.0 seconds, which is again considerably larger than those of all our heuristics variants.
Hence, we can conclude that our proposed algorithms are clearly superior.

\section{Outlook}\label{sec:outlook}
In this article, we presented several improvements for solving the control optimization problem on transient gas networks in a time-critical environment.
First, we introduced a new model for the artificial arcs representing the compressors, which is configuration-based and closer to the original compression capabilities while using fewer variables and constraints.
Furthermore, we propose a new algorithm that finds better solutions faster.
It first determines a set of solution values for all binary variables, which we call a binary assignment, and afterwards solves a corresponding non-linear program to complete it to a full solution of the problem.
The binary assignments are obtained based on sequential mixed-integer programming, two heuristics based on reduced time horizons, and a specialized dynamic node limit.
Our approach was tested on particularly challenging real-world instances featuring a large amount of supply and demand changes.
The results confirm that we are able to find solutions close to the global optimum while being very fast compared to black-box MINLP solvers and comparable approaches from the literature.

We see two main reasons for the success of our heuristic approach for the given problem.
First, we recognized that from the multitude of possible future control recommendations, only a few combinations fit a given scenario.
This is especially true for the given objective function that mainly punishes control modes and network state changes.
Second, the network's pipes act as gas storage and can therefore compensate for a local shortage or surplus of gas for a short amount of time.
For this reason, changing the control of an element a little too early or too late is often feasible and only causes a small objective function value penalty.
This helps for both of the heuristics:
The rolling horizon heuristic results in too late control changes due to not having a full view of the future.
In the case of the aggregated horizon, the time points in which control decisions can be made are limited.
Hence, they may be a little late or a little early.
Due to the pipe storage, solutions with these control time offsets are likely to be similar to one of the few reasonable control recommendations while only having a slightly worse objective function value.

Future research on the overall topic can continue in many directions.
Model-wise, we could further refine the compressor model to guarantee that the corresponding feasible region is an actual relaxation of the original feasible region.
Another possibility to improve the accuracy would be to also update the linearization of the power bound during the sequential mixed-integer programming and include the original non-linear formula in the final non-linear programming model.
However, as we have seen in the results, the run time for solving the NLP models is already quite high.
Therefore, it might be beneficial to just complete the most promising binary assignments to full solutions, for example, by comparing their corresponding linearized MIP solution value, or to only complete one of a set of similar binary assignments.
Regarding the instance set used in our computational experiments, it would be interesting to examine if heterogeneous step lengths in the time horizon lead to more complex solutions or to evaluate our algorithm in the context of pure hydrogen networks similar to those presented in~\citep{HopHenZitGot2020}.

\section*{Acknowledgements}
The work for this article has been conducted in the Research Campus MODAL, funded by the German Federal Ministry of Education and Research (BMBF) (fund numbers 05M14ZAM \& 05M20ZBM).

%
% ---- Bibliography ----
%
\typeout{} % this solves the "no citations found"-problem, no idea why, see https://www.overleaf.com/learn/latex/Questions/BibTeX_isn%27t_working%3B_my_%5Ccite_are_showing_up_as_question_marks_(%3F)
\bibliography{bibliography}

\newpage
\appendix
% change figure&table count to A.1, A.2 ... in the appendix
\renewcommand\thefigure{\thesection.\arabic{figure}}
\setcounter{figure}{0}
\renewcommand\thetable{\thesection.\arabic{table}}
\setcounter{table}{0}
\section{Appendix}

\begin{table}[htp]
\centering
\newcommand{\realPlus}{\hfill\makebox[\widthof{$\in\{0,1\}$}][l]{$\in\setRealsNonNeg$}}
\newcommand{\realAll}{\hfill\makebox[\widthof{$\in\{0,1\}$}][l]{$\in\setReals$}}
\newcommand{\binary}{\hfill$\in\{0,1\}$}

{
\def\arraystretch{1.50}% 1 is the default
\begin{tabular}{lll}
Variable & Meaning & Unit \\
\hline
\pressI{v,t}             \realPlus & pressure at node $v\in\setVertices$ & bar \\
\mFlowI{v,a,t}           \realAll  & flow into or out of pipe $a=(l,r)\in\setPipes$ at end node $v\in\{l,r\}$     & kg/s \\
\mFlowI{a,t}             \realAll  & flow over arc $a\in\setArcs\setminus\setPipes$ & kg/s \\
\inflowI{v,t}            \realAll  & inflow into boundary node $v\in\setBoundaryNodes$ & kg/s \\
\modeOp{a,t}             \binary   & selection of open mode for $a\in\setValves\cup\setRegulators$ & 1 \\
\modeCl{a,t}             \binary   & selection of closed mode for $a\in\setRegulators$ & 1 \\
\modeAc{a,t}             \binary   & selection of active mode for $a\in\setRegulators$ & 1 \\
\activeArcI{a,t}         \binary   & activity of artificial arc $a\in\setArtificialLinks$ & 1 \\
\activeFwdI{a,t}         \binary   & selection of forward direction for $a\in\setBiRegulatorArcs$ & 1 \\
\activeBwdI{a,t}         \binary   & selection of backward direction for $a\in\setBiRegulatorArcs$ & 1 \\
\mFlowFwdI{a,t}          \realPlus & forward flow of $a\in\setBiRegulatorArcs$ & kg/s \\
\mFlowBwdI{a,t}          \realPlus & backward flow of $a\in\setBiRegulatorArcs$ & kg/s \\
\activeCfgI{c,t}         \binary   & selection of configuration $c\in\setAuxCompressorConfigsI{a}$ for $a\in\setCompressorArcs$ & 1 \\
\activeStateI{s,t}       \binary   & selection of simple state $s\in\setSimpleStates$ & 1 \\
\activeFlowDirI{f,t}     \binary   & selection of flow direction $f\in\setFlowDirections$ & 1 \\
\slackPPos{v,t}          \realPlus & positive pressure slack for boundary node $v\in\setBoundaryNodes$ & bar \\
\slackPNeg{v,t}          \realPlus & negative pressure slack for boundary node $v\in\setBoundaryNodes$ & bar \\
\slackQPos{v,t}          \realPlus & positive inflow slack for boundary node $v\in\setBoundaryNodes$ & kg/s \\
\slackQNeg{v,t}          \realPlus & negative inflow slack for boundary node $v\in\setBoundaryNodes$ & kg/s \\
\valveModeChg{a,t}       \binary   & mode change for original valve $a\in\setValvesOrig$ & 1 \\
\regModeChg{a,t}         \binary   & mode change for regulator $a\in\setRegulators$ & 1 \\
\rgPLChg{a,t}           \realPlus & change of incoming pressure of active $a\in\setRegulators$ & bar \\
\rgPRChg{a,t}           \realPlus & change of outgoing pressure of active $a\in\setRegulators$ & bar \\
\rgQChg{a,t}            \realPlus & change of flow over of active $a\in\setRegulators$ & kg/s \\
\auxArcActivityChg{a,t} \binary   & activity change for artificial arc $a\in\setArtificialLinks$ & 1 \\
\simpleStateTurnOn{s,t} \binary   & activation of simple state $s\in\setSimpleStates$ & 1 \\
\fnPChg{a,t}            \realPlus & change of pressure of fence node $v\in\setFenceNodesNaviStation{i}$ for $i\in\setNetworkStations$ & bar \\
\fnQChg{a,t}            \realPlus & change of flow over valve of fence node $v\in\setFenceNodesNaviStation{i}$ for $i\in\setNetworkStations$ & kg/s
\end{tabular}
}
\caption{List of all used variables, specifying their domain, meaning and unit. Note that $1\,\text{bar}=10^5\,\text{Pa}$.}
\label{tab:Variables}
\end{table}

\end{document}